\documentclass[10pt]{article}
\usepackage{amsfonts,color}
\usepackage{amssymb}
 \usepackage{footnote}
\usepackage{mathtools}
\usepackage{amsthm,wasysym}
\usepackage{amsmath}
\usepackage[english]{babel}
\usepackage{bbm}
\usepackage{colonequals}
\usepackage{mathrsfs}

\usepackage{graphicx}
\usepackage{amsfonts}
\usepackage{hyperref}
\usepackage{subcaption}
\usepackage{nicefrac}
\usepackage{gensymb}
\usepackage{enumerate}
\usepackage[nameinlink,capitalize]{cleveref}
\usepackage{cite}
\usepackage[nottoc,numbib]{tocbibind}
\usepackage{pgf,tikz,pgfplots}
\pgfplotsset{compat=1.14}
\usepackage{mathrsfs}
\usetikzlibrary{arrows}
\usetikzlibrary{arrows.meta}
\usepackage{overpic}

\usepackage{chngcntr}

\counterwithin*{equation}{section}

\DeclareMathOperator{\at}{\bigg\vert}
\newcommand{\vb}[1]{\mathbf{#1}}
\newcommand{\bm}[1]{\boldsymbol{#1}}

\newcommand{\sym}{\mathrm{sym}\,}

\newcommand{\skw}{\mathrm{skew}\,}

\newcommand{\tr}{\mathrm{tr}\,}
\newcommand{\id}{\mathrm{id}\,}
\newcommand{\cof}{\mathrm{cof}\,}
\newcommand{\jump}[1]{\ensuremath{[\![#1]\!]} }

\newcommand{\ptr}{\mathrm{tr}\,}
\newcommand{\ttr}{\mathrm{tr}_{n}^{\perp}}
\newcommand{\rtr}{\mathrm{tr}_{t}^{\parallel}}
\newcommand{\ntr}{\mathrm{tr}_{n}^{\parallel}}

\newcommand{\dd}{\mathrm{d}}
\newcommand{\D}{\mathrm{D}}
\newcommand{\di}{\mathrm{div} \,}
\newcommand{\Di}{\mathrm{Div}\,}
\newcommand{\drot}{\mathrm{div}\bm{R}}
\newcommand{\rot}[1]{\mathrm{div}(\bm{R}\,{#1})}
\newcommand{\rog}{\bm{R} \nabla}
\newcommand{\curl}{\mathrm{curl}\,}
\newcommand{\Curl}{\mathrm{Curl}\,}

\newcommand{\Ned}{\mathcal{N}_{I}}
\newcommand{\RT}{\mathcal{RT}}
\newcommand{\BDM}{\mathcal{BDM}}
\newcommand{\Nedtwo}{\mathcal{N}_{II}}
\newcommand{\Lag}{\mathcal{L}}
\newcommand{\Ber}{\mathcal{B}}
\newcommand{\Po}{\mathit{P}}
\newcommand{\Le}{\mathit{L}^2}
\newcommand{\Lez}{\mathit{L}^2_0}
\newcommand{\Hhalf}{\mathit{H}^{1/2}}
\newcommand{\Hnhalf}{\mathit{H}^{-1/2}}
\newcommand{\Hone}{\mathit{H}^1}
\newcommand{\Honez}{\mathit{H}_0^1}
\newcommand{\Hd}[1]{\mathit{H}(\mathrm{div}{#1})}
\newcommand{\Hdz}[1]{\mathit{H}_0(\mathrm{div}{#1})}
\newcommand{\Hr}[1]{\mathit{H}(\mathrm{curl}{#1})}

\newcommand{\Hc}[1]{\mathit{H}(\mathrm{curl}{#1})}

\newcommand{\Hcz}[1]{\mathit{H}_0(\mathrm{curl}{#1})}

\newcommand{\body}{V}
\newcommand{\surf}{A}
\newcommand{\curv}{s}

\newcommand{\R}{\mathbb{R}}
\newcommand{\U}{\mathit{U}}

\newcommand{\X}{\mathit{X}}
\newcommand{\C}{\mathit{C}}

\newcommand{\tem}{\mathcal{T}}
\newcommand{\ver}{\mathcal{V}}
\newcommand{\edge}{\mathcal{E}}
\newcommand{\face}{\mathcal{F}}
\newcommand{\cell}{\mathcal{C}}

\newcommand{\tv}{\bm{\ell}}

\newcommand{\lame}{\lambda_{\mathrm{e}}}
\newcommand{\lammi}{\lambda_{\mathrm{micro}}}

\newcommand{\mue}{\mu_{\mathrm{e}}}
\newcommand{\muc}{\mu_{\mathrm{c}}}
\newcommand{\mumi}{\mu_{\mathrm{micro}}}
\newcommand{\muma}{\mu_{\mathrm{macro}}}
\newcommand{\Lc}{L_\mathrm{c}}
\newcommand{\Ce}{\mathbb{C}_{\mathrm{e}}}
\newcommand{\Cc}{\mathbb{C}_{\mathrm{c}}}
\newcommand{\Cmic}{\mathbb{C}_{\mathrm{micro}}}

\newcommand{\Pm}{\bm{P}}
\newcommand{\ud}{\vb{u}}

\newcommand{\RN}[1]{%
	\textup{\uppercase\expandafter{\romannumeral#1}}%
}

\allowdisplaybreaks[1]

\makeindex
                    
\usepackage{chngcntr}              

\newtheoremstyle{break}
{\topsep}{\topsep}%
{\itshape}{}%
{\bfseries}{}%
{\newline}{}%
\theoremstyle{break}

\newtheorem{theorem}{Theorem}

\newtheorem{remark}{Remark}

\newtheorem{definition}{Definition}

\counterwithin{theorem}{section}
\counterwithin{lemma}{section}
\counterwithin{remark}{section}
\counterwithin{observation}{section}
\counterwithin{definition}{section}

   \makeatletter
\let\@fnsymbol\@arabic
\makeatother

\crefname{Problem}{Problem.}{Problem.}

\title{Polytopal templates for the formulation of semi-continuous vectorial finite elements of arbitrary order}
\author{\normalsize{Adam Sky}\thanks{Corresponding author: Adam Sky, Institute of Structural Mechanics, Statics and Dynamics, Technische Universit\"at Dortmund, August-Schmidt-Str. 8, 44227 Dortmund, Germany, email: adam.sky@tu-dortmund.de}
	\quad and \quad
	\normalsize{Ingo Muench}\thanks{Ingo Muench, Institute of Structural Mechanics, Statics and Dynamics, Technische Universit\"at Dortmund, August-Schmidt-Str. 8, 44227 Dortmund, Germany, email: ingo.muench@tu-dortmund.de}
}

\begin{document}

\maketitle

\begin{abstract}
The Hilbert spaces $\Hc{}$ and $\Hd{}$ are needed for variational problems formulated in the context of the de Rham complex in order to guarantee well-posedness. Consequently, the construction of conforming subspaces is a crucial step in the formulation of viable numerical solutions. Alternatively to the standard definition of a finite element as per Ciarlet, given by the triplet of a domain, a polynomial space and degrees of freedom, this work aims to introduce a novel, simple method of directly constructing semi-continuous vectorial base functions on the reference element via polytopal templates and an underlying $\Hone$-conforming polynomial subspace.
The base functions are then mapped from the reference element to the element in the physical domain via consistent Piola transformations. 
The method is defined in such a way, that the underlying $\Hone$-conforming subspace can be chosen independently, thus allowing for constructions of arbitrary polynomial order. The base functions arise by multiplication of the basis with template vectors defined for each polytope of the reference element. We prove a unisolvent construction of N\'ed\'elec elements of the first and second type, Brezzi-Douglas-Marini elements, and Raviart-Thomas elements. An application for the method is demonstrated with two examples in the relaxed micromorphic model.

\vspace*{0.25cm}


\end{abstract}

\section{Introduction}
In many variational problems, well-posedness necessitates the use of either the $\Hc{}$ or $\Hd{}$ Hilbert spaces. Some classical examples are Maxwell's equations \cite{Mon03,MONK1993101,Lee} and mixed Poisson problems \cite{Banz}. More recent examples are the tangential-displacement-normal-normal-stress (TDNNS) method in elasticity \cite{NEUNTEUFEL2021113857,Pechstein2018}, and the relaxed micromorphic model \cite{Neff2014,Neff2015,Reg}. 
Other examples are curl based plasticity models \cite{Cordero2010,Ebobisse2018}. 
Commonly, the application of the $\Hc{}$ or $\Hd{}$ spaces arises in problems associated with the de Rham complex \cite{Demkowicz2000,PaulyDeRham}.
In fact, as shown in \cite{Arnold2021}, complexes are a powerful and general tool for mixed variational formulations. The latter is also demonstrated in \cite{PaulyEl} and \cite{PaulyDiv,Botti,Hu} where the elasticity complex and the $\di\Di$-complex are explored in the context of mixed formulations of linear elasticity and the biharmonic equation, respectively.

Since analytical solutions to partial differential equations are rarely possible for general domain geometries or boundary conditions, the application of numerical schemes with conforming subspaces is required. 
Unlike in the classical Hilbert space $\Hone$, an element of the $\Hc{}$-space is only required to be tangentially continuous \cite{Solin}. Analogously, elements of the $\Hd{}$-space are only required to be normal-continuous \cite{Solin}. 
As such, the formulation of $\Hc{}$- and $\Hd{}$-conforming finite elements is more complex. The pioneering works \cite{Nedelec1980} and \cite{Ned2} introduced the N\'ed\'elec elements of the first and second types, which represent polynomial subspaces with the minimal regularity requirements of the $\Hc{}$-space $\Nedtwo^p \subset \Ned^p \subset \Hc{}$. In \cite{BDM} and \cite{Raviart} the authors introduced the Brezzi-Douglas-Marini and Raviart-Thomas elements, that allow to construct polynomial subspaces for the $\Hd{}$-space $\BDM^p \subset \RT^p \subset \Hd{}$, such that the elements exhibit the minimal regularity needed in $\Hd{}$. 
The elements are given in the classical element definition as per Ciarlet \cite{Ciarl78}, and allow for application on general grids.

An alternative methodology to construction of a basis directly on the grid is to build base functions on the reference elements and map them to the physical elements on the grid by consistent transformations. The construction of low order vectorial finite elements is demonstrated in \cite{Anjam2015,Sky2021,SKY2022115298,SkyOn,Schroder2022}. The formulation of higher order elements on the basis of Legendre polynomials can be found in \cite{Zaglmayr2006,Joachim2005,Solin}. Further, in \cite{AINSWORTH2018178,Ainsworth2015}, the authors present a higher order construction based on Bernstein polynomials. 
We note that in general, the mapping alone does not suffice in order to assert a consistent transformation and some additional algorithmic is required in order to avoid the orientation problem \cite{Sky2021,SKY2022115298,Fuentes2015,Anjam2015,Zaglmayr2006,Ainsworth2003}.

The aim of this work is to establish a method of defining $\Hc{}$ and $\Hd{}$ base functions on the reference element, such that the underlying polynomial basis can be chosen independently. As such, the method allows to directly construct conforming finite elements by using for example, Lagrange, Legendre, Jacobi or Bernstein polynomials. 
This goal is achieved by defining a template on the reference element, which can be subsequently used in conjunction with an $\Hone$-conforming polynomial basis of one's choice, in order to span a semi-continuous finite element space. 
The template is composed of vector sets associated with the polytopes of the reference element. Consequently, we dub the methodology "polytopal templates". In this work we consider subspaces for the Hilbert spaces $\Hc{}$ and $\Hd{}$.

This paper is structured as follows. First, we introduce the classical Hilbert spaces and their corresponding differential and trace operators. Next, we derive two-dimensional polytopal templates for the construction of N\'ed\'elec elements of the first and second type, Brezzi-Douglas-Marini elements, and Raviart-Thomas elements on the reference triangle. The methodology is subsequently utilized to derive polytopal templates on the reference tetrahedron for N\'ed\'elec elements of the second type and Brezzi-Douglas-Marini elements. We demonstrate the application of elements using the relaxed micromorphic model with one example in antiplane shear and one three-dimensional example.
Lastly, we present our conclusions and outlook. 

\begin{figure}
    	\centering
    	\definecolor{asl}{rgb}{0.4980392156862745,0.,1.}
    	\definecolor{asb}{rgb}{0.,0.4,0.6}
    	\begin{tikzpicture}[line cap=round,line join=round,>=triangle 45,x=1.0cm,y=1.0cm]
    		\clip(-4,-0.5) rectangle (12.5,4.5);
    		
    		\def\x{-8.5}
    		\def\y{4}
    		
    		\draw [color=asb,line width=1.pt] (5+\x,-3+\y) -- (6.5+\x,-3+\y) -- (6.5+\x,-1.5+\y) -- (5+\x, -1.5+\y) -- (5+\x, -3+\y);
    		\fill[opacity=0.1, asb] (5+\x,-3+\y) -- (6.5+\x,-3+\y) -- (6.5+\x,-1.5+\y) -- (5+\x, -1.5+\y) -- cycle;
    		
    		\draw [-to,color=black,line width=1.pt] (5+\x,-3+\y) -- (7+\x,-3+\y);
    		\draw [-to,color=black,line width=1.pt] (5+\x,-3+\y) -- (5+\x,-1+\y);
    		\draw (7+\x,-3+\y) node[color=black,anchor=west] {$\xi$};
    		\draw (5+\x,-1+\y) node[color=black,anchor=south] {$\eta$};
    		\draw (5.75+\x,-2.25+\y) node[color=asb] {$\Omega$};
    		
    		\draw [-to,color=asl,line width=1.pt] (5.75+\x,-1.5+\y) -- (5.75+\x,-1.+\y);
    		\draw (5.75+\x,-1.+\y) node[color=asl,anchor=south] {$\bm{\nu}$};
    		
    		\draw [-to,color=asl,line width=1.pt] (5.5+\x,-1.5+\y) -- (6+\x,-1.5+\y);
    		\draw (5.9+\x,-1.5+\y) node[color=asl,anchor=south west] {$\bm{\tau}$};
    		
    		\draw [-Triangle,color=black,line width=1.pt] (-1.5,2) -- (0.5,2);
    		\draw (-1.5,2) node[color=black,anchor=south west] {$\vb{x}:\Omega \to \body$};
    		
    		\fill [asb, opacity=0.1] plot [smooth cycle] coordinates {(1,3) (3,4) (7, 2) (10,3) (12,1) (10,0) (5,1) (2,1)};
    		
    		\begin{scope}
    			\clip(6,-0.5) rectangle (12.5,4.5);
    			\draw [asl, densely dashed] plot [smooth cycle] coordinates {(1,3) (3,4) (7, 2) (10,3) (12,1) (10,0) (5,1) (2,1)};
    		\end{scope}
    	    \begin{scope}
    	    	\clip(0,-0.5) rectangle (6,4.5);
    	    	\draw [asb] plot [smooth cycle] coordinates {(1,3) (3,4) (7, 2) (10,3) (12,1) (10,0) (5,1) (2,1)};
    	    \end{scope}
    		
    		\draw [-to,color=black,line width=1.pt] (0,0) -- (1,0);
    		\draw [-to,color=black,line width=1.pt] (0,0) -- (0,1);
    		\draw (1,0) node[color=black,anchor=west] {$x$};
    		\draw (0,1) node[color=black,anchor=south] {$y$};
    		
    		\draw [-to,color=asl,line width=1.pt] (11.35,2) -- (12,2.55);
    		\draw (11.95,2.5) node[color=asl,anchor=south west] {$\vb{n}$};
    		\draw [-to,color=asl,line width=1.pt] (11.35-0.22,2+0.25) -- (11.35+0.25,2-0.28);
    		\draw (11.35+0.3,2-0.25) node[color=asl,anchor=west] {$\vb{t}$};
    		
    		\draw (6.5,1.15) node[color=asb,anchor=south] {$\body$};
    		\draw (4.3,3.5) node[color=asb,anchor=west] {$\surf_D$};
    		\draw (11.6,0.4) node[color=asl,anchor=north] {$\surf_N$};
    	\end{tikzpicture}
        \caption{A domain $\body$ with Dirichlet and Neumann boundaries mapped from some reference domain $\Omega$.}
    	\label{fig:domain_elast}
    \end{figure}
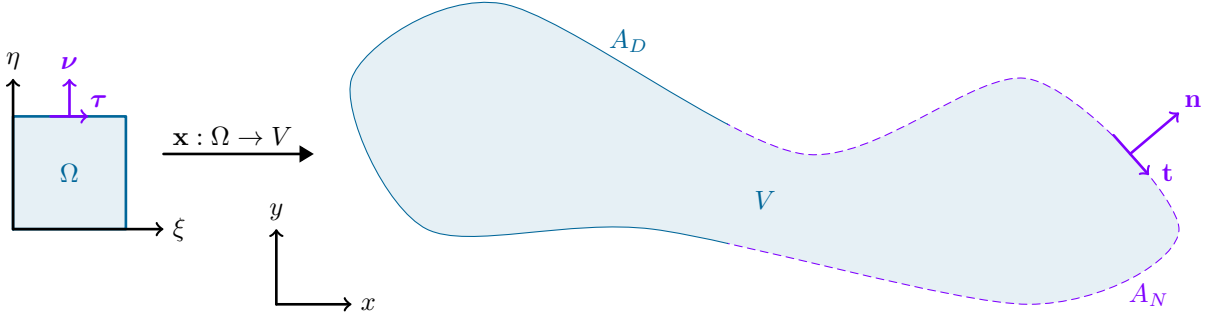

\hfill \break
The following definitions are employed throughout this work, see also \cref{fig:domain_elast}:
\begin{itemize}
    \item Vectors are indicated by bold letters. Non-bold letters represent scalars.
    \item In general, formulas are defined using the Cartesian basis, where the base vectors are denoted by $\vb{e}_1$, $\vb{e}_2$ and $\vb{e}_3$.
    \item Three-dimensional domains in the physical space are denoted with $\body \subset \R^3$. The corresponding reference domain is given by $\Omega$.
    \item Analogously, in two dimensions we employ $\surf \subset \R^2$ for the physical domain and $\Gamma$ for the reference domain. 
    \item Curves on the physical domain are denoted by $\curv$, whereas curves in the reference domain by $\mu$.
    \item The tangent and normal vectors in the physical domain are given by $\vb{t}$ and $\vb{n}$, respectively. Their counterparts in the reference domain are $\bm{\tau}$ for tangent vectors and $\bm{\nu}$ for normal vectors.
\end{itemize}

\section{Hilbert spaces and trace operators}
Hilbert spaces are the natural function spaces in the formulation of variational problems \cite{Mon03}.  
In preparation for the construction of conforming subspaces we introduce the classical  Hilbert spaces and their associated norms
\begin{subequations}
		\begin{align}
			\Hone(\body) &= \{ u \in \Le(\body) \; | \; \nabla u \in [\Le(\body)]^3 \} \, , & \|u \|^2_{\Hone(\body)} &= \|u \|^2_{\Le} + \| \nabla u \|^2_{\Le} \, , \\
			\Hc{, \body} &= \{ \vb{u} \in [\Le(\body)]^3 \; | \; \curl \vb{u} \in [\Le(\body)]^3 \} \, , & \|u \|^2_{\Hc{,\body}} &= \|\vb{u} \|^2_{\Le} + \| \curl \vb{u} \|^2_{\Le} \, , \\
			\Hd{, \body} &= \{ \vb{u} \in [\Le(\body)]^3 \; | \; \di \vb{u} \in \Le(\body) \} \, , & \|u \|^2_{\Hd{,\body}} &= \|\vb{u} \|^2_{\Le} + \| \di \vb{u} \|^2_{\Le} \, ,
		\end{align}
	\end{subequations}
which are based on the Lebesgue space
\begin{align}
		&\Le(\body) = \{ u : \body \to  \mathbb{R} \; | \; \| u \|_{\Le(\body)} < \infty  \} \, , && \| u \|^2_{\Le(\body)} = \langle u , \, u \rangle_{\Le(\body)} = \int_\body u^2 \, \dd \body \, .
	\end{align}
Note that on two-dimensional domains the differential operators are reduced to
\begin{align}
    &\nabla u = \begin{bmatrix} u_{,x} \\ u_{,y} \end{bmatrix} \, ,
    &&\rot{\ud} = u_{2,x} - u_{1,y} \, , &&  \rog u = \begin{bmatrix} u_{,y} \\ -u_{,x}  \end{bmatrix} \, , && \bm{R} = \begin{bmatrix}
    0 & 1 \\
    -1 & 0
    \end{bmatrix} \, , && \di \ud = u_{1,x} + u_{2,y} \, , 
\end{align}
such that two curl operators are introduced: one for two-dimensional vectors $\rot{\cdot}$, and one for scalars $\rog(\cdot)$. 
On contractible domains the Hilbert spaces are connected by the exact de Rham sequence \cite{Arnold2021,Demkowicz2000,PaulyDeRham,DEMKOWICZ2005267}, see \cref{fig:derham}.

Each Hilbert space is associated with a corresponding trace operator \cite{Hiptmair}. The trace of the function $u \in \Hone(\body)$ is defined by the linear bounded operator 
\begin{align}
				&\ptr u = u \at_{\partial \body} \in \Hhalf(\partial \body) \, , && \exists \, c > 0: \quad \| \ptr u \|_{\Hhalf(\partial \body)} \leq c \| u \|_{\Hone(\body)} \quad \forall \, u \in \Hone(\body) \, .
\end{align}
In other words, the trace restricts the function to the boundary of the domain. 
The trace of the $\Hc{,\body}$ space is given by the tangential components on the boundary. On a surface the tangent vector is not unique, and therefore, the tangential projection is defined using the normal vector and the cross product
\begin{align}
				&\ttr \vb{u} = \vb{n} \times \vb{u} \at_{\partial \body} \in [\Hnhalf(\partial \body)]^3 \, , &&\exists \, c > 0: \quad \|\ttr \vb{u}\|_{\Hnhalf(\partial \body)} \leq c \|\vb{u} \|_{\Hc{,\body}} \quad \forall \,\vb{u} \in \Hc{,\body} \, .
\end{align}
For two-dimensional domains $\surf \subset \R^2$ the tangent vector is unique and the trace operator reduces to
\begin{align}
				&\rtr \vb{u} = \langle \vb{t} ,\, \vb{u} \rangle \at_{\partial \surf} \in \Hnhalf(\partial \surf) \, .
\end{align}
Lastly, the trace of the $\Hd{,\body}$ space is defined by the normal projection at the boundary
\begin{align}
				&\ntr\vb{u} = \langle \vb{n} ,\, \vb{u} \rangle \at_{\partial \body} \in \Hnhalf(\partial \body) \, , && \exists \, c > 0: \quad \|\ntr\vb{u}\|_{\Hnhalf(\partial \body)} \leq c \|\vb{u} \|_{\Hd{,\body}} \quad \forall \,\vb{u} \in \Hd{,\body} \, .
\end{align}
In this work we define the trace space via
\begin{align}
    \Hhalf(\partial \body) = \{ v \in \Le(\partial \body) \; | \; \exists \, u \in \Hone(\body) :  \ptr u = v  \} \, , 
\end{align}
and $\Hnhalf(\partial \body)$ is its dual. A thorough treatment of fractional Sobolev spaces is found in \cite{DINEZZA2012521}.
The trace operators are used to define Hilbert spaces with boundary conditions. The exact de Rham sequence holds also on Hilbert spaces with vanishing traces, see \cref{fig:derham0}.
Further, the trace operators allow to identify finite elements of a specific Hilbert space via interface conditions. We state the interface theorem \cite{neunteufel2021mixed,Solin} and apply it to the construction of conforming subspaces.
\begin{theorem}[Interface conditions]
		
		A finite element space is a conforming subspace of a Hilbert space if and only if the jump of the trace of its elements vanishes for all arbitrarily defined interfaces $\Xi_{ij} = \body_i \cap \body_j \, , i \neq j$ where $\body = \body_i \cup \body_j \subset \R^3$ and $\Xi_{ij} \subset \R^{2}$ (and analogously for two-dimensional domains) 
		\begin{subequations}
			\begin{align}
				&u \in \Hone(\body) && \iff &  \jump{ \ptr u }\at_{\Xi_{ij}} &= 0 \quad \forall \, \Xi_{ij} = \body_i  \cap \body_j  \, , \\
				& \vb{u} \in \Hc{,\body} && \iff &  \jump{ \ttr \vb{u} }\at_{\Xi_{ij}} &= 0 \quad \forall \, \Xi_{ij} = \body_i  \cap \body_j  \, , \\
				& \vb{u} \in \Hd{,\body} && \iff &  \jump{ \ntr \vb{u} }\at_{\Xi_{ij}} &= 0 \quad \forall \, \Xi_{ij} = \body_i  \cap \body_j  \, .
			\end{align}
		\end{subequations}
	\end{theorem}

In the following sections the construction of arbitrary order $\Hc{}$- and $\Hd{}$-conforming subspaces is presented. The construction is based on a polytopal association of base functions.
	\begin{definition} [Polytopal base functions] \label{def:poly}
		Each base function is associated with its respective polytope and the underlying Hilbert space as follows:
		\begin{enumerate}
			\item A vertex base function has a  vanishing trace on all other vertices and non-neighbouring edges and faces.
			\item An edge base function has a vanishing trace on all other edges and non-neighbouring faces. 
			\item A face base function has a vanishing trace on all other faces.
			\item A cell base function has a vanishing trace on the entire boundary of the element.
		\end{enumerate}
		The definition is general and the respective trace may change according to the corresponding Hilbert space.
	\end{definition}
\begin{figure}
    \centering
    \begin{tikzpicture}[line cap=round,line join=round,>=triangle 45,x=1.0cm,y=1.0cm]
				\clip(0.,0.3) rectangle (13.3,1.0);
				\draw (0.,0.9) node[anchor=north west] {$\R$};
				\draw [-Triangle,line width=.5pt] (0.6,0.6) -- (2.1,0.6);
				\draw (0.9,1.1) node[anchor=north west] {$\id$};
				\draw (2.1,0.9) node[anchor=north west] {$\Hone(\body)$};
				\draw [-Triangle,line width=.5pt] (3.5,0.6) -- (5.,0.6);
				\draw (3.8,1.1) node[anchor=north west] {$\nabla$};
				\draw (5.,0.9) node[anchor=north west] {$\Hc{,\body}$};
				\draw [-Triangle,line width=.5pt] (7.1,0.6) -- (8.6,0.6);
				\draw (7.3,1.1) node[anchor=north west] {$\curl$};
				\draw (8.6,0.9) node[anchor=north west] {$\Hd{,\body}$};
				\draw [-Triangle,line width=.5pt] (10.5,0.6) -- (12.,0.6);
				\draw (10.8,1.1) node[anchor=north west] {$\di$};
				\draw (12.,0.9) node[anchor=north west] {$\Le(\body)$};
			\end{tikzpicture}
\begin{tikzpicture}[line cap=round,line join=round,>=triangle 45,x=1.0cm,y=1.0cm]
				\clip(0.,0.0) rectangle (10.0,1.5);
				\draw (0.,0.9) node[anchor=north west] {$\R$};
				\draw [-Triangle,line width=.5pt] (0.6,0.6) -- (2.1,0.6);
				\draw (0.9,1.1) node[anchor=north west] {$\id$};
				\draw (2.1,0.9) node[anchor=north west] {$\Hone(\surf)$};
				\draw [-Triangle,line width=.5pt] (3.5,0.6) -- (5.,0.6);
				\draw (3.8,1.1) node[anchor=north west] {$\nabla$};
				\draw (5.,0.9) node[anchor=north west] {$\Hc{,\surf}$};
				\draw [-Triangle,line width=.5pt] (7.,0.6) -- (8.5,0.6);
				\draw (7.1,1.1) node[anchor=north west] {$\drot$};
				\draw (8.5,0.9) node[anchor=north west] {$\Le(\surf)$};
			\end{tikzpicture}
\begin{tikzpicture}[line cap=round,line join=round,>=triangle 45,x=1.0cm,y=1.0cm]
				\clip(0.,0.0) rectangle (10.0,1.3);
				\draw (0.,0.9) node[anchor=north west] {$\R$};
				\draw [-Triangle,line width=.5pt] (0.6,0.6) -- (2.1,0.6);
				\draw (0.9,1.1) node[anchor=north west] {$\id$};
				\draw (2.1,0.9) node[anchor=north west] {$\Hone(\surf)$};
				\draw [-Triangle,line width=.5pt] (3.5,0.6) -- (5.,0.6);
				\draw (3.8,1.1) node[anchor=north west] {$\rog$};
				\draw (5.,0.9) node[anchor=north west] {$\Hd{,\surf}$};
				\draw [-Triangle,line width=.5pt] (6.8,0.6) -- (8.3,0.6);
				\draw (7.1,1.1) node[anchor=north west] {$\di$};
				\draw (8.3,0.9) node[anchor=north west] {$\Le(\surf)$};
			\end{tikzpicture}
    \caption{Classical de Rham exact sequences for three- and two-dimensional contractible domains.
			The range of each operator is exactly the kernel of the next operator in the sequence.}
    \label{fig:derham}
\end{figure}
\begin{figure}
    \centering
    \begin{tikzpicture}[line cap=round,line join=round,>=triangle 45,x=1.0cm,y=1.0cm]
				\clip(0.,0.3) rectangle (13.3,1.0);
				\draw (0.,0.9) node[anchor=north west] {$\R$};
				\draw [-Triangle,line width=.5pt] (0.6,0.6) -- (2.1,0.6);
				\draw (0.9,1.1) node[anchor=north west] {$\id$};
				\draw (2.1,0.9) node[anchor=north west] {$\Honez(\body)$};
				\draw [-Triangle,line width=.5pt] (3.5,0.6) -- (5.,0.6);
				\draw (3.8,1.1) node[anchor=north west] {$\nabla$};
				\draw (5.,0.9) node[anchor=north west] {$\Hcz{,\body}$};
				\draw [-Triangle,line width=.5pt] (7.1,0.6) -- (8.6,0.6);
				\draw (7.3,1.1) node[anchor=north west] {$\curl$};
				\draw (8.6,0.9) node[anchor=north west] {$\Hdz{,\body}$};
				\draw [-Triangle,line width=.5pt] (10.6,0.6) -- (12.1,0.6);
				\draw (10.9,1.1) node[anchor=north west] {$\di$};
				\draw (12.1,0.9) node[anchor=north west] {$\Lez(\body)$};
			\end{tikzpicture}
\begin{tikzpicture}[line cap=round,line join=round,>=triangle 45,x=1.0cm,y=1.0cm]
				\clip(0.,0.0) rectangle (10.0,1.5);
				\draw (0.,0.9) node[anchor=north west] {$\R$};
				\draw [-Triangle,line width=.5pt] (0.6,0.6) -- (2.1,0.6);
				\draw (0.9,1.1) node[anchor=north west] {$\id$};
				\draw (2.1,0.9) node[anchor=north west] {$\Honez(\surf)$};
				\draw [-Triangle,line width=.5pt] (3.5,0.6) -- (5.,0.6);
				\draw (3.8,1.1) node[anchor=north west] {$\nabla$};
				\draw (5.,0.9) node[anchor=north west] {$\Hcz{,\surf}$};
				\draw [-Triangle,line width=.5pt] (7.,0.6) -- (8.5,0.6);
				\draw (7.1,1.1) node[anchor=north west] {$\drot$};
				\draw (8.5,0.9) node[anchor=north west] {$\Lez(\surf)$};
			\end{tikzpicture}
\begin{tikzpicture}[line cap=round,line join=round,>=triangle 45,x=1.0cm,y=1.0cm]
				\clip(0.,0.0) rectangle (10.0,1.3);
				\draw (0.,0.9) node[anchor=north west] {$\R$};
				\draw [-Triangle,line width=.5pt] (0.6,0.6) -- (2.1,0.6);
				\draw (0.9,1.1) node[anchor=north west] {$\id$};
				\draw (2.1,0.9) node[anchor=north west] {$\Honez(\surf)$};
				\draw [-Triangle,line width=.5pt] (3.5,0.6) -- (5.,0.6);
				\draw (3.8,1.1) node[anchor=north west] {$\rog$};
				\draw (5.,0.9) node[anchor=north west] {$\Hdz{,\surf}$};
				\draw [-Triangle,line width=.5pt] (6.9,0.6) -- (8.4,0.6);
				\draw (7.2,1.1) node[anchor=north west] {$\di$};
				\draw (8.4,0.9) node[anchor=north west] {$\Lez(\surf)$};
			\end{tikzpicture}
    \caption{De Rham exact sequences for Hilbert spaces with vanishing traces. The Lebesgue zero-space is characterized by functions with a vanishing integral over the domain.}
    \label{fig:derham0}
\end{figure}

\section{Two-dimensional templates}
\begin{figure}
    \centering
    \definecolor{asl}{rgb}{0.4980392156862745,0.,1.}
		\definecolor{asb}{rgb}{0.,0.4,0.6}
		\begin{tikzpicture}[line cap=round,line join=round,>=triangle 45,x=1.0cm,y=1.0cm]
			\clip(-1.5,-1) rectangle (4.5,4.5);
			\draw (-0.5,-0.5) node[circle,fill=asb,inner sep=1.5pt] {};
			\draw (-0.5,4) node[circle,fill=asb,inner sep=1.5pt] {};
			\draw (4,-0.5) node[circle,fill=asb,inner sep=1.5pt] {};
			\draw [color=asb,line width=.6pt] (-0.5,0) -- (-0.5,3);
			\draw [color=asb,line width=.6pt] (0,-0.5) -- (3,-0.5);
			\draw [color=asb,line width=.6pt] (0.3,3.3) -- (3.3,0.3);
			\draw [dotted,color=asb,line width=.6pt] (0,0) -- (0,3) -- (3,0) -- (0,0);
			\fill[opacity=0.1, asb] (0,0) -- (0,3) -- (3,0) -- cycle;
			\draw (-0.5,-0.5) node[color=asb,anchor=north east] {$_{v_1}$};
			\draw (4,-0.5) node[color=asb,anchor=north west] {$_{v_3}$};
			\draw (-0.5,4) node[color=asb,anchor=south east] {$_{v_2}$};
			
			\draw (-0.58,1.5) node[color=asb,anchor=east] {$_{e_{12}}$};
			\draw (1.5,-.5) node[color=asb,anchor=north] {$_{e_{13}}$};
			\draw (1.94,1.94) node[color=asb,anchor=north west] {$_{e_{23}}$};
			
			\draw (0.65,0.7)
			node[color=asb,anchor=south west] {$_{c_{123}}$};
			
		\end{tikzpicture}
    \caption{Decomposition of the unit triangle into vertices, edges and the cell.}
    \label{fig:tri}
\end{figure}
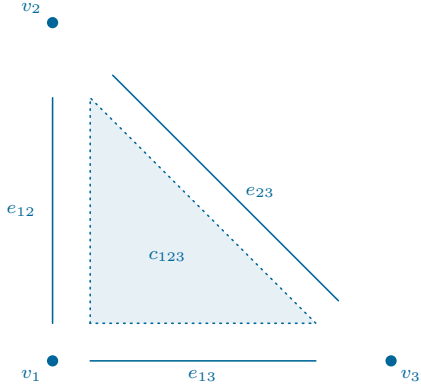

\begin{figure}
    	\centering
    	\begin{subfigure}{0.3\linewidth}
    		\centering
    		\includegraphics[width=0.7\linewidth]{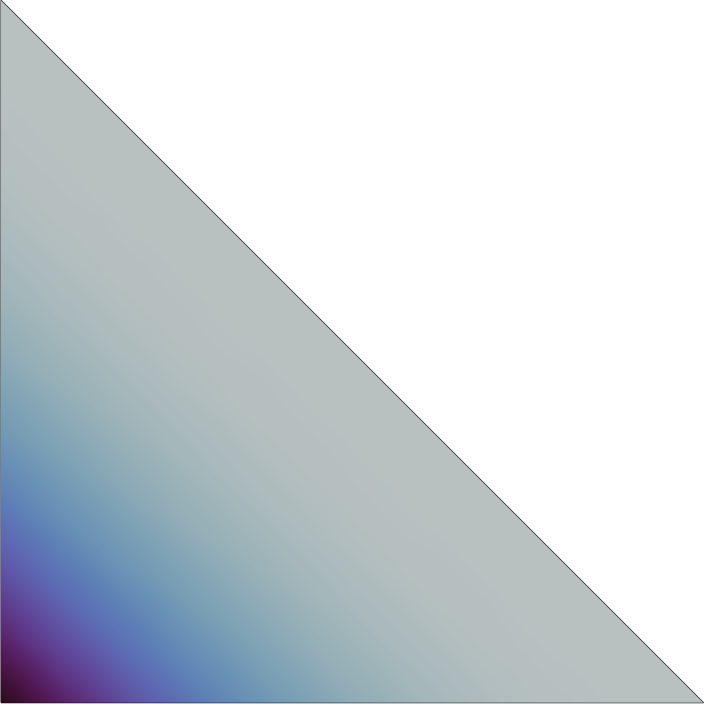}
    		\caption{}
    	\end{subfigure}
    	\begin{subfigure}{0.3\linewidth}
    		\centering
    		\includegraphics[width=0.7\linewidth]{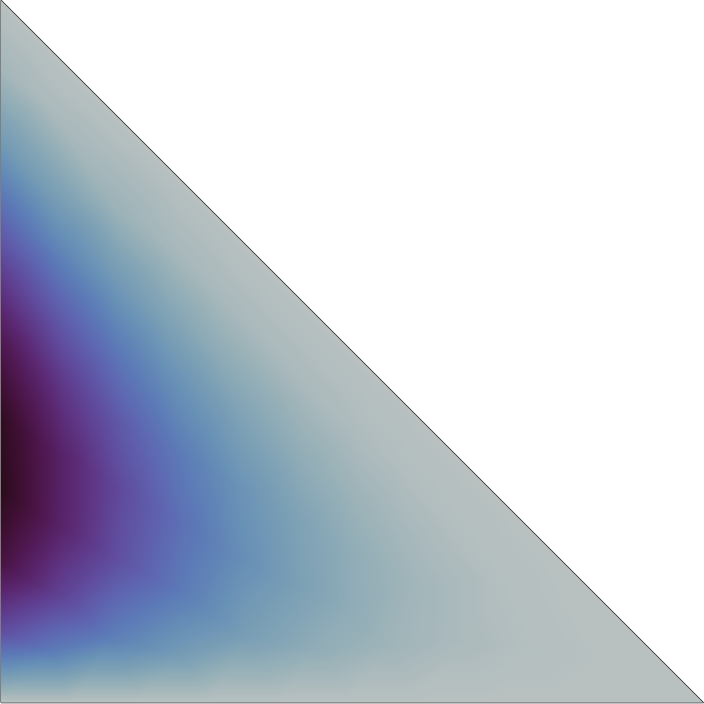}
    		\caption{}
    	\end{subfigure}
    	\begin{subfigure}{0.3\linewidth}
    		\centering
    		\includegraphics[width=0.7\linewidth]{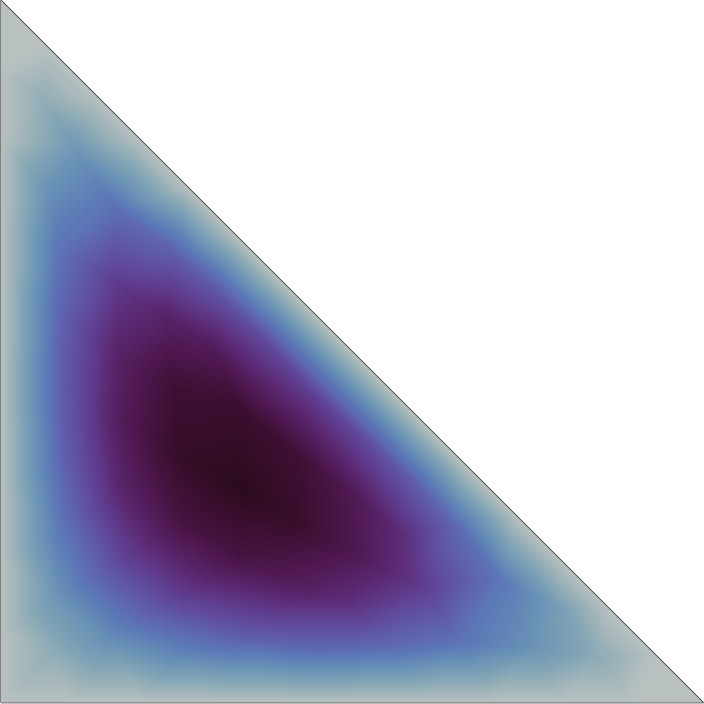}
    		\caption{}
    	\end{subfigure}
    	\caption{Vertex (a), edge (b) and cell (c) base functions on the reference triangle.
    	The transition from light to dark represents the increase in value of the function. The lowest value of each base function is zero and the maximum value may change according to the basis.}
    	\label{fig:base}
    \end{figure}

This section is dedicated to the introduction of polytopal templates on the reference triangle 
\begin{align}
    \Gamma = \{ (\xi, \eta) \in [0,1]^2 \; | \; \xi + \eta \leq 1 \} \, .
\end{align}
To that end, the triangle is decomposed into its base polytopes given by its vertices $\{v_1, v_2, v_3\}$, its edges $\{e_{12}, e_{13}, e_{23}\}$, and its interior cell $c_{123}$, see \cref{fig:tri}. 
Further, each polytope is associated with base functions belonging to an $\Hone$-conforming subspace $\U^p(\Gamma)$ with $\dim \U^p(\Gamma) = \dim \Po^p(\Gamma) = (p+2)(p+1)/2$. 
\begin{definition}[Triangle $\U^p(\Gamma)$-polytopal spaces]
Each polytope is associated with a space of base functions as follows:
\begin{itemize}
    \item Each vertex is associated with the space of its respective base function $\ver^p_i$. As such, there are three spaces in total $i \in \{1,2,3\}$ and each one is of dimension one, $\dim \ver^p_i = 1 \quad \forall \, i \in \{1,2,3\}$.
    \item For each edge there exists a space of edge functions $\edge^p_j$ with the multi-index $j \in \mathcal{J} = \{(1,2),(1,3),(2,3)\}$.
    The dimension of each edge space is given by $\dim \edge^p_j = p-1$.
    \item Lastly, the cell is equipped with the space of cell base functions $\cell_{123}^p$ with $\dim \cell_{123}^p = (p-2)(p-1)/2$.
\end{itemize}
The association with a respective polytope reflects \cref{def:poly} with the trace operator for $\Hone$-spaces.
\end{definition}
A depiction of vertex, edge, and cell base functions is given in \cref{fig:base}. Clearly, this is the standard definition of base functions for approximations in $\Hone$. Common examples of such bases are Lagrange, Legendre and Bernstein.

In the following we define templates on the reference triangle, such that their multiplications with corresponding base functions from $\U^p(\Gamma)$ generate vectorial base functions for either $\Hc{}$ or $\Hd{}$.

\subsection{N\'ed\'elec II}
\begin{figure}
		\centering
		\definecolor{asl}{rgb}{0.4980392156862745,0.,1.}
		\definecolor{asb}{rgb}{0.,0.4,0.6}
		\begin{tikzpicture}[line cap=round,line join=round,>=triangle 45,x=1.0cm,y=1.0cm]
			\clip(-2,-1.5) rectangle (12.5,4.5);
			\draw (-0.5,-0.5) node[circle,fill=asb,inner sep=1.5pt] {};
			\draw (-0.5,4) node[circle,fill=asb,inner sep=1.5pt] {};
			\draw (4,-0.5) node[circle,fill=asb,inner sep=1.5pt] {};
			\draw [color=asb,line width=.6pt] (-0.5,0) -- (-0.5,3);
			\draw [color=asb,line width=.6pt] (0,-0.5) -- (3,-0.5);
			\draw [color=asb,line width=.6pt] (0.3,3.3) -- (3.3,0.3);
			\draw [dotted,color=asb,line width=.6pt] (0,0) -- (0,3) -- (3,0) -- (0,0);
			\fill[opacity=0.1, asb] (0,0) -- (0,3) -- (3,0) -- cycle;
			\draw (-0.5,-0.5) node[color=asb,anchor=north east] {$_{v_1}$};
			\draw (4,-0.5) node[color=asb,anchor=north west] {$_{v_3}$};
			\draw (-0.5,4) node[color=asb,anchor=south east] {$_{v_2}$};
			
			\draw (-0.58,1.5) node[color=asb,anchor=west] {$_{e_{12}}$};
			\draw (1.5,-.5) node[color=asb,anchor=south] {$_{e_{13}}$};
			\draw (1.94,1.94) node[color=asb,anchor=north east] {$_{e_{23}}$};
			
			\draw [-to,color=asl,line width=1pt] (-0.75,0) -- (-0.75,0.75);
			\draw [-to,color=asl,line width=1pt] (-1.5,2.25) -- (-0.75,3);
			
			\draw [-to,color=asl,line width=1pt] (0,-0.75) -- (0.75,-0.75);
			\draw [-to,color=asl,line width=1pt] (2.25,-1.5) -- (3,-0.75);
			
			\draw [-to,color=asl,line width=1pt] (0.45,3.45) -- (1.2,3.45);
			\draw [to-,color=asl,line width=1pt] (3.45,0.45) -- (3.45,1.2);
			
			\draw [to-,color=asl,line width=1pt, dashdotted] (-1.5,1.5) -- (-0.75,1.5);
			\draw [-to,color=asl,line width=1pt, densely dashed] (-0.75,1.125) -- (-0.75,1.875);
			
			\draw [-to,color=asl,line width=1pt, dashdotted] (1.5,-1.5) -- (1.5,-0.75);
			\draw [-to,color=asl,line width=1pt, densely dashed] (1.125,-0.75) -- (1.875,-0.75);
			
			\draw [-to,color=asl,line width=1pt, dashdotted] (2,2) -- (2.75,2.75);
			\draw [-to,color=asl,line width=1pt, densely dashed] (1.7,2.3) -- (2.3,1.7);
			
			\draw [-to,color=asl,line width=1pt, dotted] (0.65,0.65) -- (1.35,0.65);
			\draw [-to,color=asl,line width=1pt, dotted] (0.65,0.65) -- (0.65,1.35);
			
			\draw (0.65,0.7)
			node[color=asb,anchor=south west] {$_{c_{123}}$};
			
			\draw [-to,color=asl,line width=1pt] (6,3) -- (7,3);
			\draw (7,3)
			node[color=asl,anchor=west] {Vertex-edge template vectors};
			\draw [-to,color=asl,line width=1pt, densely dashed] (6,2.25) -- (7,2.25);
			\draw (7,2.25)
			node[color=asl,anchor=west] {Edge template vectors};
			\draw [-to,color=asl,line width=1pt, dashdotted] (6,1.5) -- (7,1.5);
			\draw (7,1.5)
			node[color=asl,anchor=west] {Edge-cell template vectors};
			\draw [-to,color=asl,line width=1pt, dotted] (6,0.75) -- (7,0.75);
			\draw (7,0.75)
			node[color=asl,anchor=west] {Cell template vectors};
		\end{tikzpicture}
		\caption{Template vectors for the reference N\'ed\'elec triangle element of the second type on their corresponding polytope.}
		\label{fig:tri_nii}
	\end{figure}
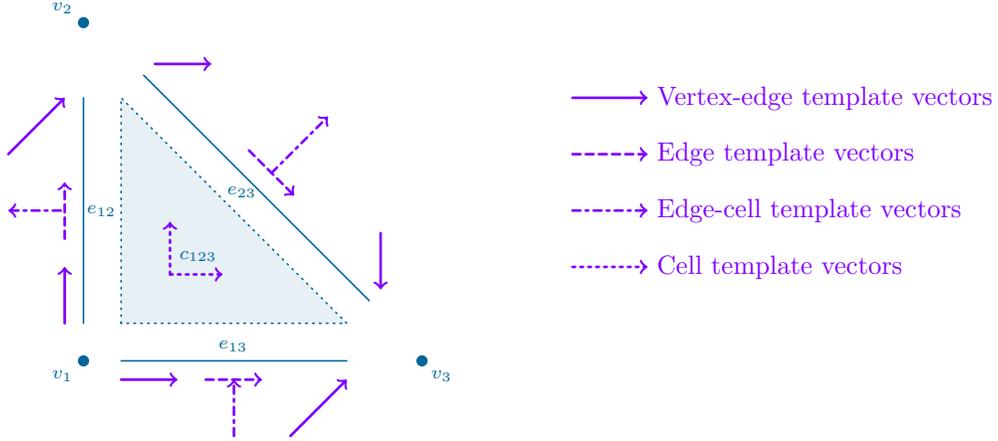

\begin{figure}
    	\centering
    	\begin{subfigure}{0.24\linewidth}
    		\centering
    		\includegraphics[width=0.95\linewidth]{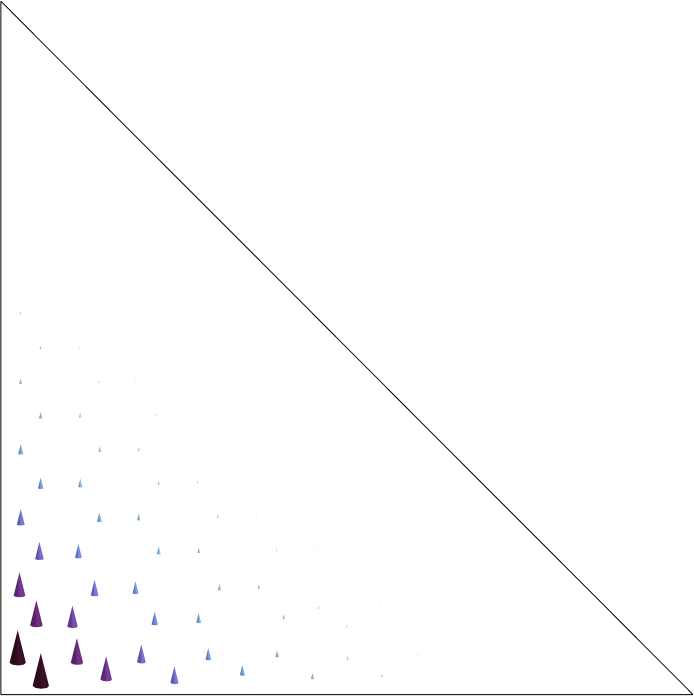}
    		\caption{}
    	\end{subfigure}
        \begin{subfigure}{0.24\linewidth}
        	\centering
        	\includegraphics[width=0.95\linewidth]{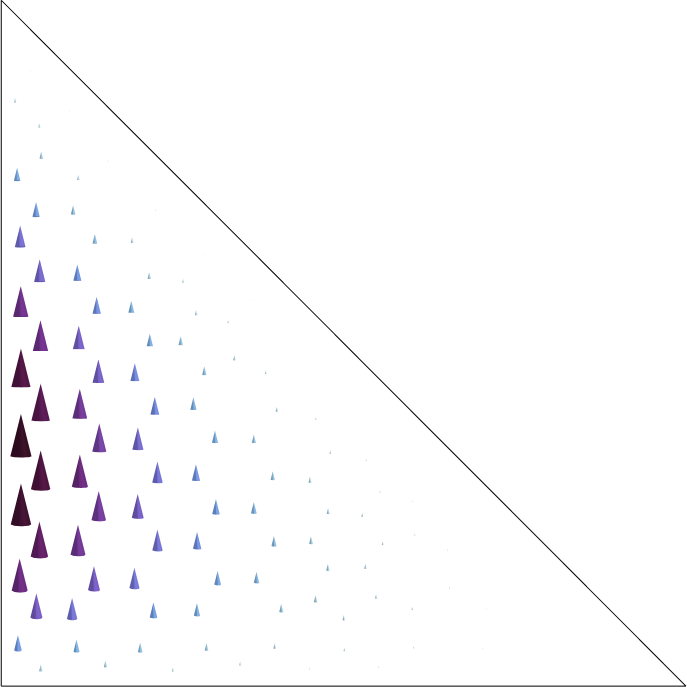}
        	\caption{}
        \end{subfigure}
        \begin{subfigure}{0.24\linewidth}
        	\centering
        	\includegraphics[width=0.95\linewidth]{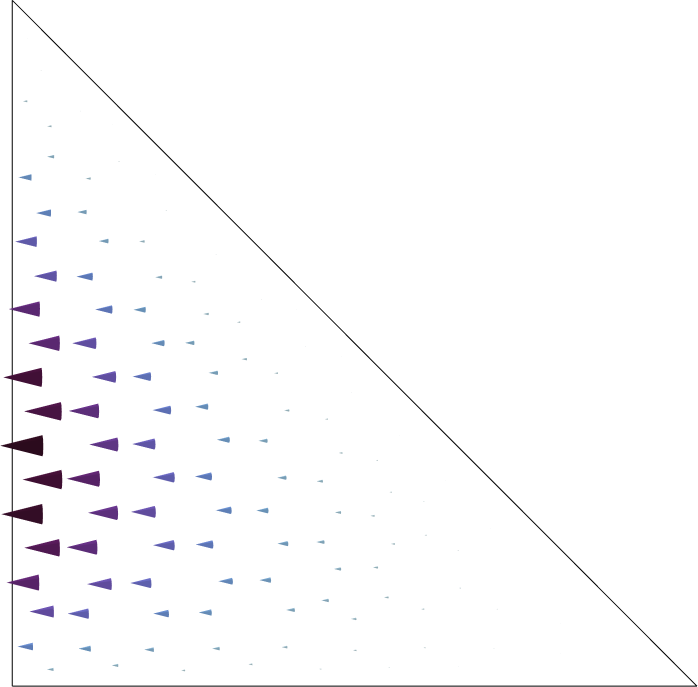}
        	\caption{}
        \end{subfigure}
        \begin{subfigure}{0.24\linewidth}
        	\centering
        	\includegraphics[width=0.95\linewidth]{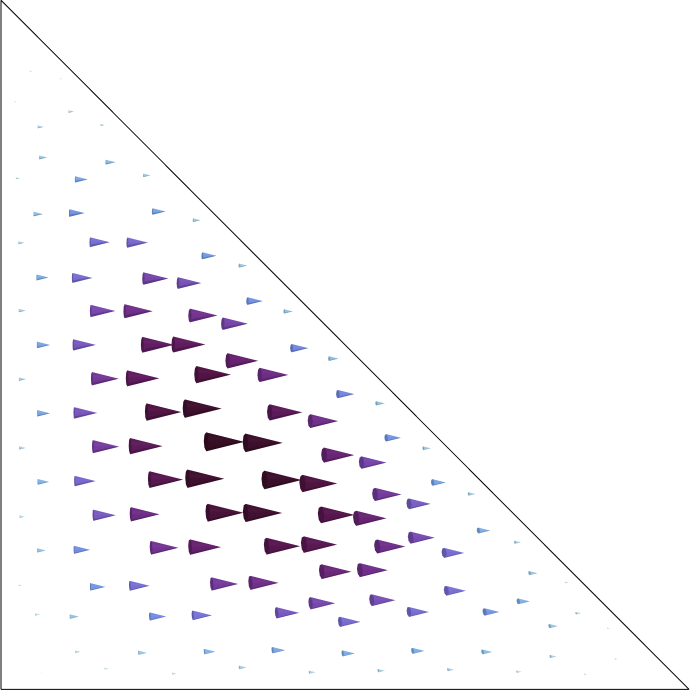}
        	\caption{}
        \end{subfigure}
    	\caption{Cubic vertex-edge (a), edge (b), edge-cell (c) and pure cell (d) vectorial base functions of the N\'ed\'elec element of the second type on the reference triangle.
    	The colour and size of the arrows represents the intensity of the field.}
    	\label{fig:nedpii}
    \end{figure}

In order to construct a template for the N\'ed\'elec element of the second type \cite{Ned2} we consider the decomposition of the reference triangle in \cref{fig:tri}. On the first vertex $v_1$ we define a vector with a projection of one on the tangent vector of the first edge $e_{12}$ and a zero projection on the second edge $e_{13}$. Next we define a vector with a projection of one on the tangent vector of the first edge $e_{12}$. Further, we construct a normal vector on the first edge $e_{12}$. Lastly, we define two unit vectors in the cell. 
The remaining vectors for their respective polytopes are computed by mapping the triangle $c_{123}$ to various permutations of $c_{ijk}$
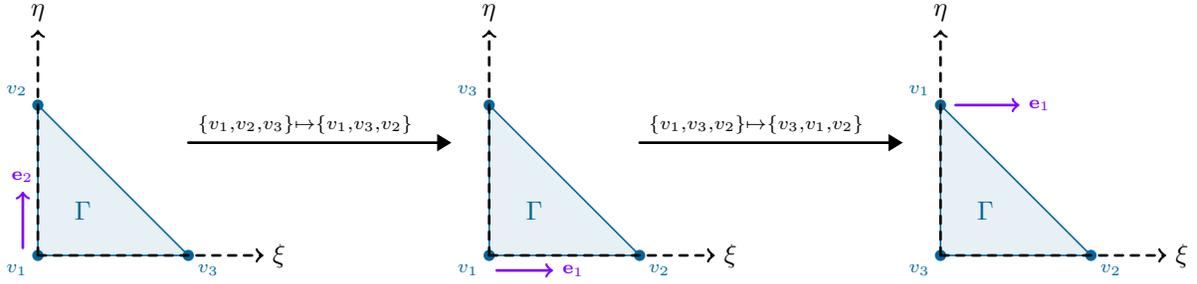
\begin{figure}
		\centering
		\definecolor{asl}{rgb}{0.4980392156862745,0.,1.}
\definecolor{asb}{rgb}{0.,0.4,0.6}
\begin{tikzpicture}[line cap=round,line join=round,>=triangle 45,x=1.0cm,y=1.0cm]
	\clip(-0.5,-0.5) rectangle (15.5,3.5);
	
	\draw (0,0) node[circle,fill=asb,inner sep=1.5pt] {};
	\draw (0,2) node[circle,fill=asb,inner sep=1.5pt] {};
	\draw (2,0) node[circle,fill=asb,inner sep=1.5pt] {};
	\draw [color=asb,line width=.6pt] (0,0) -- (0,2) -- (2,0) -- (0,0);
	\fill[opacity=0.1, asb] (0,0) -- (0,2) -- (2,0) -- cycle;
	\draw (0.6,0.6) node[color=asb] {$\Gamma$};
	\draw (0,0) node[color=asb,anchor=north east] {$_{v_1}$};
	\draw (2,0) node[color=asb,anchor=north west] {$_{v_3}$};
	\draw (0,2) node[color=asb,anchor=south east] {$_{v_2}$};
	\draw [-to,color=black,line width=1.pt, dashed] (0,0) -- (3,0);
	\draw [-to,color=black,line width=1.pt, dashed] (0,0) -- (0,3);
	\draw (3,0) node[color=black,anchor=west] {$\xi$};
	\draw (0,3) node[color=black,anchor=south] {$\eta$};
	
	\draw [-to,color=asl,line width=1pt] (-0.2,0.1) -- (-0.2,0.85);
	\draw (-0.2,0.85) node[color=asl,anchor=south] {$_{\vb{e}_2}$};
	
	\draw (6,0) node[circle,fill=asb,inner sep=1.5pt] {};
	\draw (6,2) node[circle,fill=asb,inner sep=1.5pt] {};
	\draw (8,0) node[circle,fill=asb,inner sep=1.5pt] {};
	\draw [color=asb,line width=.6pt] (6,0) -- (6,2) -- (8,0) -- (6,0);
	\fill[opacity=0.1, asb] (6,0) -- (6,2) -- (8,0) -- cycle;
	\draw (6.6,0.6) node[color=asb] {$\Gamma$};
	\draw (6,0) node[color=asb,anchor=north east] {$_{v_1}$};
	\draw (8,0) node[color=asb,anchor=north west] {$_{v_2}$};
	\draw (6,2) node[color=asb,anchor=south east] {$_{v_3}$};
	\draw [-to,color=black,line width=1.pt, dashed] (6,0) -- (9,0);
	\draw [-to,color=black,line width=1.pt, dashed] (6,0) -- (6,3);
	\draw (9,0) node[color=black,anchor=west] {$\xi$};
	\draw (6,3) node[color=black,anchor=south] {$\eta$};
	
	\draw [-to,color=asl,line width=1pt] (6.1, -0.2) -- (6.85,-0.2);
	\draw (6.85,-0.2) node[color=asl,anchor=west] {$_{\vb{e}_1}$};
	
	\draw (12,0) node[circle,fill=asb,inner sep=1.5pt] {};
	\draw (12,2) node[circle,fill=asb,inner sep=1.5pt] {};
	\draw (14,0) node[circle,fill=asb,inner sep=1.5pt] {};
	\draw [color=asb,line width=.6pt] (12,0) -- (12,2) -- (14,0) -- (12,0);
	\fill[opacity=0.1, asb] (12,0) -- (12,2) -- (14,0) -- cycle;
	\draw (12.6,0.6) node[color=asb] {$\Gamma$};
	\draw (12,0) node[color=asb,anchor=north east] {$_{v_3}$};
	\draw (14,0) node[color=asb,anchor=north west] {$_{v_2}$};
	\draw (12,2) node[color=asb,anchor=south east] {$_{v_1}$};
	\draw [-to,color=black,line width=1.pt, dashed] (12,0) -- (15,0);
	\draw [-to,color=black,line width=1.pt, dashed] (12,0) -- (12,3);
	\draw (15,0) node[color=black,anchor=west] {$\xi$};
	\draw (12,3) node[color=black,anchor=south] {$\eta$};
	
	\draw (2,1.5) node[color=black,anchor=south west] {$_{\{v_1,v_2,v_3\} \mapsto \{v_1,v_3,v_2\}}$};
	\draw [-Triangle,color=black,line width=1.pt] (2,1.5) -- (5.5,1.5);
	
	\draw [-to,color=asl,line width=1pt] (12.2,2) -- (13.05,2);
	\draw (13.05,2) node[color=asl,anchor=west] {$_{\vb{e}_1}$};
	
	\draw (8,1.5) node[color=black,anchor=south west] {$_{\{v_1,v_3,v_2\} \mapsto \{v_3,v_1,v_2\}}$};
	\draw [-Triangle,color=black,line width=1.pt] (8,1.5) -- (11.5,1.5);
\end{tikzpicture}
		\caption{Derivation of a template vector on the remaining edges via permutations of the reference triangle using covariant Piola mappings.}
		\label{fig:permut}
	\end{figure}
on the unit domain via covariant Piola transformations (see \cref{ap:a}) and adjusting the sign to ensure a positive projection on the tangent vector, see \cref{fig:permut}. The complete template is depicted in \cref{fig:tri_nii}.
\begin{remark}
		The polytopal set is not unique and depends on the starting definition on the first polytopes and the resulting mapping. Further, one can change pure edge-type template vectors by adding or subtracting normal vectors without influencing the tangential projection. 
	\end{remark}
The resulting template is given by the super-set of the sets for the respective polytopes
	\begin{align}
		\tem = \{\tem_1,\tem_2,\tem_3,\tem_{12},\tem_{13},\tem_{23},\tem_{123}\} \, , 
	\end{align}
	where the polytopal sets read
	\begin{align}
		\tem_1 &= \{\vb{e}_1,\vb{e}_2\} \, , & \tem_2 &= \{\vb{e}_1 + \vb{e}_2,\vb{e}_1\} \, , & \tem_3 &= \{\vb{e}_1 + \vb{e}_2,-\vb{e}_2\} \, , \notag \\
		\tem_{12} &= \{\vb{e}_2,-\vb{e}_1\} \, , & \tem_{13} &= \{\vb{e}_1,\vb{e}_2\} \, , & \tem_{23} &= \{ (1/2) (\vb{e}_1 - \vb{e}_2),\vb{e}_1 + \vb{e}_2\} \, , \notag \\
		\tem_{123} &= \{\vb{e}_1,\vb{e}_2\} \, .
	\end{align}
	The set can now be used in conjunction with an underlying $\U^p(\Gamma)$-space to construct a N\'ed\'elec element of the second type of order $p$
	\begin{align}
			&\Nedtwo^p = \left \{ \bigoplus_{i=1}^3 \ver^p_i \otimes \tem_i \right \} \oplus \left \{ \bigoplus_{j \in \mathcal{J}} \edge^p_{j} \otimes \tem_{j} \right \} \oplus  \{\cell^p_{123} \otimes \tem_{123}\} \, , && \mathcal{J} = \{(1,2),(1,3),(2,3)\} \, , 
		\end{align}
		where $\ver^p_i$ are the sets of the vertex base functions, $\edge^p_{j}$ are the sets of edge base functions and $\cell^p_{123}$ is the set of cell base functions.
	\begin{definition}[Triangle $\Nedtwo^p$ base functions]
	The base functions of the triangle N\'ed\'elec element of the second type are defined on their respective polytope as follows.
	\begin{itemize}
	    \item On each edge $e_{ij}$ with vertices $v_i$ and $v_j$, the base functions read
	    \begin{subequations}
	        \begin{align}
	        \text{Vertex-edge:}& &\bm{\vartheta}(\xi, \eta) &= n \tv \, , & n &\in \ver^p_i \, , &\tv &\in \left \{ \tv \in \tem_i \; | \; \rtr \tv \at_{e_{ij}}  \neq 0 \right \} \, , \\
	        &&\bm{\vartheta}(\xi, \eta) &= n \tv \, , & n &\in \ver^p_j \, , &\tv &\in \left \{ \tv \in \tem_j \; | \; \rtr \tv \at_{e_{ij}}  \neq 0 \right \} \, , \\
	        \text{Edge:} &&\bm{\vartheta}(\xi, \eta) &= n \tv \, , & n &\in \edge^p_{ij} \, , &\tv &\in \left \{ \tv \in \tem_{ij} \; | \; \rtr \tv \at_{e_{ij}} \neq 0  \right \} \, ,
	    \end{align}
	    \end{subequations}
	    such that $\tv$ is an element of the template sets, whose tangential trace does not vanish on the edge. 
	    \item The cell base functions are given by  
	    \begin{subequations}
	        \begin{align}
	        \text{Edge-cell:}& &\bm{\vartheta}(\xi, \eta) &= n \tv \, , & n &\in \edge^p_{12} \, , &\tv &\in \left \{ \tv \in \tem_{12} \; | \; \rtr \tv \at_{e_{12}} = 0  \right \} \, , \\
	        & &\bm{\vartheta}(\xi, \eta) &= n \tv \, , & n &\in \edge^p_{13} \, , &\tv &\in \left \{ \tv \in \tem_{13} \; | \; \rtr \tv \at_{e_{13}} = 0  \right \} \, , \\
	        & &\bm{\vartheta}(\xi, \eta) &= n \tv \, , & n &\in \edge^p_{23} \, , &\tv &\in \left \{ \tv \in \tem_{23} \; | \; \rtr \tv \at_{e_{23}} = 0  \right \} \, , \\
	        \text{Cell:}& &\bm{\vartheta}(\xi, \eta) &= n \tv \, , & n &\in \cell^p_{123} \, , &\tv &\in \tem_{123} \, ,
	    \end{align}
	    \end{subequations}
	    such that their tangential trace vanishes on all edges.
	\end{itemize}
	\end{definition}
	A selection of cubic $\Nedtwo^p(\Gamma)$-base functions given by the construction is depicted in \cref{fig:nedpii}.
	\begin{theorem} [Linear independence]
		The tensor product of the template with an $\Hone$-conforming polynomial basis $\U^p$ yields a unisolvent N\'ed\'elec element of the second type.
		\label{th:li_tri_nedii}
	\end{theorem}
	\begin{proof}
		Under the assertion that the underlying $\Hone$-conforming polynomial basis $\U^p$ is unisolvent, unisolvence of the $\Nedtwo^p$ basis follows automatically since each base function of $\U^p$ is multiplied with two linearly independent vectors, thus inheriting the linear independence of the basis on the vectorial level. Further, the resulting basis has the required dimensionality
		\begin{align}
			\dim [\U^p(\Gamma)]^2 = \dim [\Po^p(\Gamma)]^2 = \dim \Nedtwo^p(\Gamma) \, , 
		\end{align}
		of the N\'ed\'elec finite element space.
	\end{proof}
	\begin{theorem} [$\Hc{,\surf}$-conformity]
		The constructed element on the reference domain is $\Hc{}$-conforming under covariant Piola transformations of the base functions.
		\label{th:con_tri_nedii}
	\end{theorem}
	\begin{proof}
		By construction, the tangential projection of each non-cell base function on the tangential vector of its respective polytope is the underlying $\Hone$-conforming base function
		\begin{align}
			\langle \bm{\tau}, \, \bm{\vartheta}_i \rangle = n_i \, .
		\end{align} 
		Since the templates are constructed by permutations of the reference element, this characteristic is extended to every corresponding polytope. Lastly, the covariant Piola transformations uphold the tangential projections in the physical domain. 
		As such, the condition at the interfaces of neighbouring elements is reduced from $\jump{\rtr\ud}|_{\Xi_{ij}} = 0$ to $\jump{\tr \langle \vb{t} , \, \ud \rangle}|_{\Xi_{ij}} = 0$ for the tangential components. This property is upheld by the underlying $\Hone$-conforming subspace.
	\end{proof}

\subsection{Brezzi-Douglas-Marini}
\begin{figure}
		\centering
		\definecolor{asl}{rgb}{0.4980392156862745,0.,1.}
		\definecolor{asb}{rgb}{0.,0.4,0.6}
		\begin{tikzpicture}[line cap=round,line join=round,>=triangle 45,x=1.0cm,y=1.0cm]
			\clip(-2,-1.5) rectangle (12.5,4.5);
			\draw (-0.5,-0.5) node[circle,fill=asb,inner sep=1.5pt] {};
			\draw (-0.5,4) node[circle,fill=asb,inner sep=1.5pt] {};
			\draw (4,-0.5) node[circle,fill=asb,inner sep=1.5pt] {};
			\draw [color=asb,line width=.6pt] (-0.5,0) -- (-0.5,3);
			\draw [color=asb,line width=.6pt] (0,-0.5) -- (3,-0.5);
			\draw [color=asb,line width=.6pt] (0.3,3.3) -- (3.3,0.3);
			\draw [dotted,color=asb,line width=.6pt] (0,0) -- (0,3) -- (3,0) -- (0,0);
			\fill[opacity=0.1, asb] (0,0) -- (0,3) -- (3,0) -- cycle;
			\draw (-0.5,-0.5) node[color=asb,anchor=north east] {$_{v_1}$};
			\draw (4,-0.5) node[color=asb,anchor=north west] {$_{v_3}$};
			\draw (-0.5,4) node[color=asb,anchor=south east] {$_{v_2}$};
			
			\draw (-0.58,1.5) node[color=asb,anchor=west] {$_{e_{12}}$};
			\draw (1.5,-.5) node[color=asb,anchor=south] {$_{e_{13}}$};
			\draw (1.94,1.94) node[color=asb,anchor=north east] {$_{e_{23}}$};
			
			\draw [-to,color=asl,line width=1pt] (-1.5,0) -- (-0.75,0);
			\draw [-to,color=asl,line width=1pt] (-1.5,3.75) -- (-0.75,3);
			
			\draw [-to,color=asl,line width=1pt] (0,-0.75) -- (0,-1.5);
			\draw [to-,color=asl,line width=1pt] (3.75,-1.5) -- (3,-0.75);
			
			\draw [to-,color=asl,line width=1pt] (0.45,3.45) -- (0.45,4.2);
			\draw [to-,color=asl,line width=1pt] (3.45,0.45) -- (4.2,0.45);
			
			\draw [-to,color=asl,line width=1pt, densely dashed] (-1.5,1.5) -- (-0.75,1.5);
			\draw [-to,color=asl,line width=1pt, dashdotted] (-0.75,1.125) -- (-0.75,1.875);
			
			\draw [to-,color=asl,line width=1pt, densely dashed] (1.5,-1.5) -- (1.5,-0.75);
			\draw [-to,color=asl,line width=1pt, dashdotted] (1.125,-0.75) -- (1.875,-0.75);
			
			\draw [to-,color=asl,line width=1pt, densely dashed] (2,2) -- (2.75,2.75);
			\draw [-to,color=asl,line width=1pt, dashdotted] (1.7,2.3) -- (2.3,1.7);
			
			\draw [-to,color=asl,line width=1pt, dotted] (0.65,0.65) -- (1.35,0.65);
			\draw [-to,color=asl,line width=1pt, dotted] (0.65,0.65) -- (0.65,1.35);
			
			\draw (0.65,0.7)
			node[color=asb,anchor=south west] {$_{c_{123}}$};
			
			\draw [-to,color=asl,line width=1pt] (6,3) -- (7,3);
			\draw (7,3)
			node[color=asl,anchor=west] {Vertex-edge template vectors};
			\draw [-to,color=asl,line width=1pt, densely dashed] (6,2.25) -- (7,2.25);
			\draw (7,2.25)
			node[color=asl,anchor=west] {Edge template vectors};
			\draw [-to,color=asl,line width=1pt, dashdotted] (6,1.5) -- (7,1.5);
			\draw (7,1.5)
			node[color=asl,anchor=west] {Edge-cell template vectors};
			\draw [-to,color=asl,line width=1pt, dotted] (6,0.75) -- (7,0.75);
			\draw (7,0.75)
			node[color=asl,anchor=west] {Cell template vectors};
		\end{tikzpicture}
		\caption{Template vectors for the reference Brezzi-Douglas-Marini triangle element on their corresponding polytope.}
		\label{fig:tri_bdm}
	\end{figure}
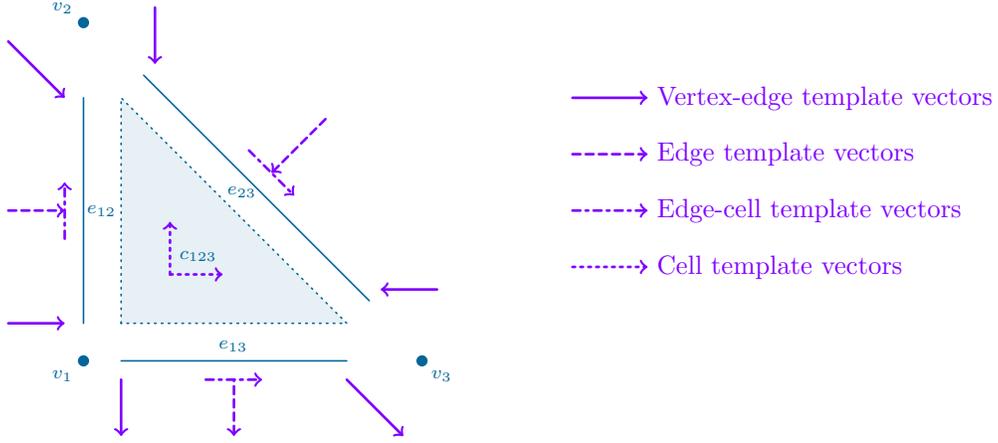

\begin{figure}
    	\centering
    	\begin{subfigure}{0.24\linewidth}
    		\centering
    		\includegraphics[width=0.95\linewidth]{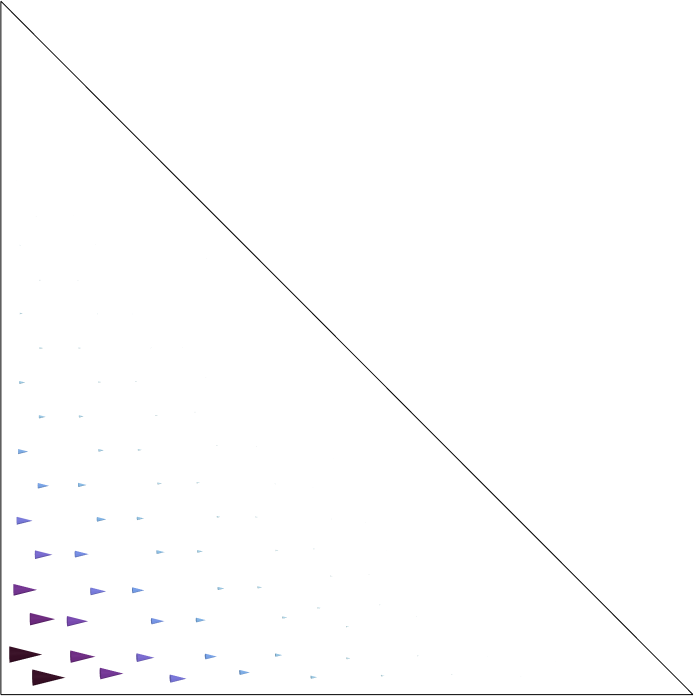}
    		\caption{}
    	\end{subfigure}
        \begin{subfigure}{0.24\linewidth}
        	\centering
        	\includegraphics[width=0.95\linewidth]{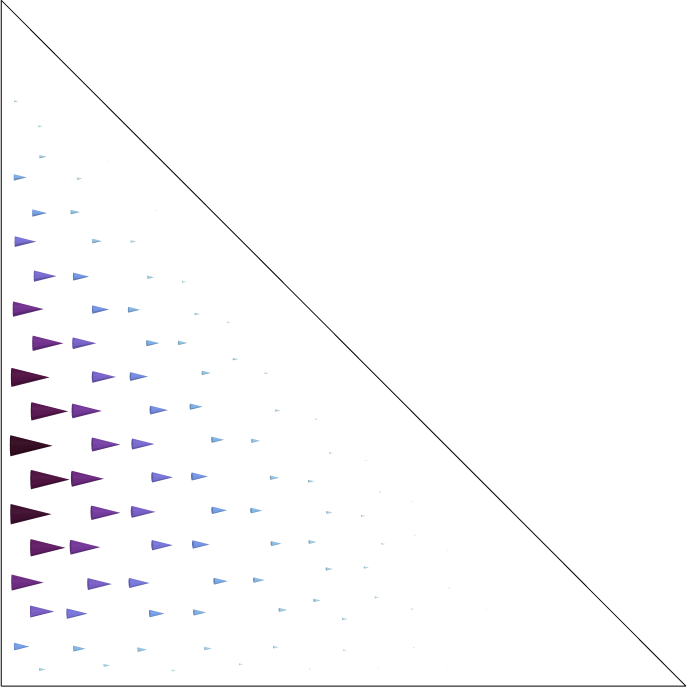}
        	\caption{}
        \end{subfigure}
        \begin{subfigure}{0.24\linewidth}
        	\centering
        	\includegraphics[width=0.95\linewidth]{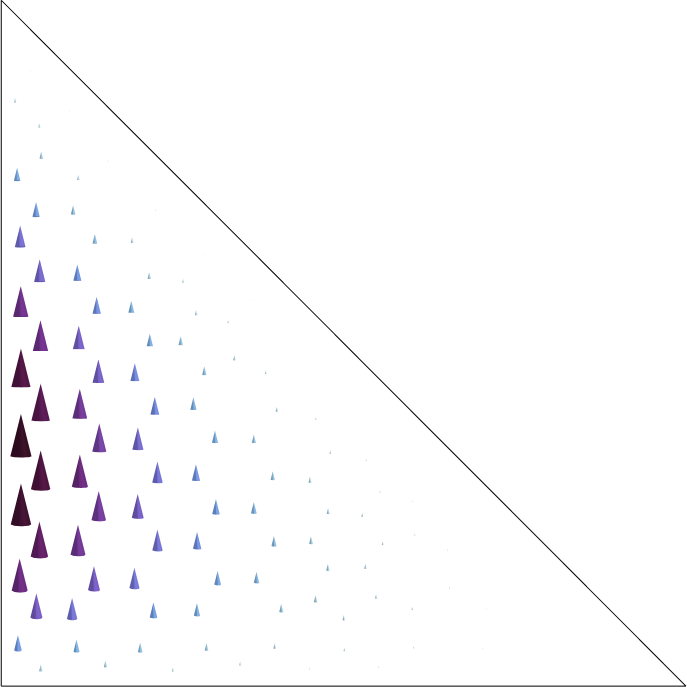}
        	\caption{}
        \end{subfigure}
        \begin{subfigure}{0.24\linewidth}
        	\centering
        	\includegraphics[width=0.95\linewidth]{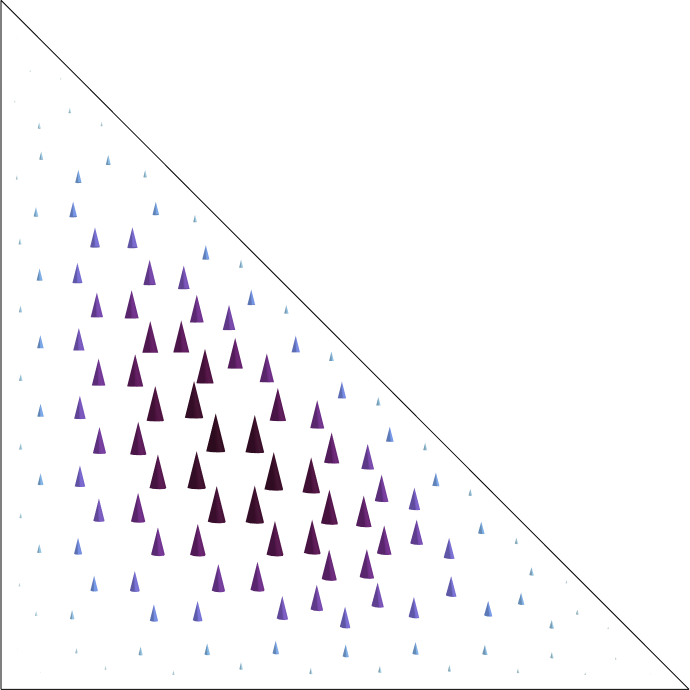}
        	\caption{}
        \end{subfigure}
    	\caption{Cubic vertex-edge (a), edge (b), edge-cell (c) and pure cell (d) base functions of the Brezzi-Douglas-Marini element on the reference triangle.}
    	\label{fig:bdm}
    \end{figure}

Next we construct the template for the Brezzi-Douglas-Marini element \cite{BDM} on the reference triangle. 
We define two normal vectors for the vertex $v_1$, one normal vector and one tangent vector on $e_{12}$ and the Cartesian basis in the cell $c_{123}$. By permutations of the reference triangle $c_{ijk}$ via contravariant Piola transformations (see \cref{ap:a}) we retrieve the remaining template vectors, compare with \cref{fig:permut}. The sign of the vectors is adjusted in order to ensure conformal normal projections on the edge normals.
Further, on the slanted edge, we modify the edge normal with the tangential vector to generate an orthogonal template.
The resulting template is depicted in \cref{fig:tri_bdm} and is given by the super-set of the sets belonging to the respective polytopes
	\begin{align}
		\tem = \{\tem_1,\tem_2,\tem_3,\tem_{12},\tem_{13},\tem_{23},\tem_{123}\} \, , 
	\end{align}
	where the polytopal sets read
	\begin{align}
		\tem_1 &= \{\vb{e}_1,-\vb{e}_2\} \, , & \tem_2 &= \{\vb{e}_1 - \vb{e}_2,-\vb{e}_2\} \, , & \tem_3 &= \{\vb{e}_1 - \vb{e}_2, -\vb{e}_1\} \, , \notag \\
		\tem_{12} &= \{\vb{e}_1,\vb{e}_2\} \, , & \tem_{13} &= \{-\vb{e}_2,\vb{e}_1\} \, , & \tem_{23} &= \{ -(1/2) (\vb{e}_1 + \vb{e}_2),\vb{e}_2 - \vb{e}_1\} \, , \notag \\
		\tem_{123} &= \{\vb{e}_1,\vb{e}_2\} \, .
	\end{align}
With the polytopal template at hand, we can construct the Brezzi-Douglas-Marini element by tensor products with an underlying $\U^p(\Gamma)$-basis
\begin{align}
			&\BDM^p = \left \{ \bigoplus_{i=1}^3 \ver^p_i \otimes \tem_i \right \} \oplus \left \{ \bigoplus_{j \in \mathcal{J}} \edge^p_{j} \otimes \tem_{j} \right \} \oplus  \{\cell^p_{123} \otimes \tem_{123}\} \, , && \mathcal{J} = \{(1,2),(1,3),(2,3)\} \, , 
		\end{align}
		where $\ver^p_i$ are the sets of vertex base functions, $\edge^p_{j}$ contain the edge base functions and $\cell^p_{123}$ is the set of cell base functions.
\begin{definition}[Triangle $\BDM^p$ base functions]
	The base functions of the triangle Brezzi-Douglas-Marini element are defined on their respective polytopes as follows.
	\begin{itemize}
	    \item On each edge $e_{ij}$ with vertices $v_i$ and $v_j$, the base functions read
	    \begin{subequations}
	        \begin{align}
	        \text{Vertex-edge:}& &\bm{\phi}(\xi, \eta) &= n \tv \, , & n &\in \ver^p_i \, , &\tv &\in \left \{ \tv \in \tem_i \; | \; \ntr \tv \at_{e_{ij}}  \neq 0 \right \} \, , \\
	        &&\bm{\phi}(\xi, \eta) &= n \tv \, , & n &\in \ver^p_j \, , &\tv &\in \left \{ \tv \in \tem_j \; | \; \ntr \tv \at_{e_{ij}}  \neq 0 \right \} \, , \\
	        \text{Edge:} &&\bm{\phi}(\xi, \eta) &= n \tv \, , & n &\in \edge^p_{ij} \, , &\tv &\in \left \{ \tv \in \tem_{ij} \; | \; \ntr \tv \at_{e_{ij}} \neq 0  \right \} \, ,
	    \end{align}
	    \end{subequations}
	    such that $\tv$ is an element of the template sets, whose normal trace does not vanish on the edge. 
	    \item The cell base functions are given by all permutations of 
	    \begin{subequations}
	        \begin{align}
	        \text{Edge-cell:}& &\bm{\phi}(\xi, \eta) &= n \tv \, , & n &\in \edge^p_{12} \, , &\tv &\in \left \{ \tv \in \tem_{12} \; | \; \ntr \tv \at_{e_{12}} = 0  \right \} \, , \\
	        & &\bm{\phi}(\xi, \eta) &= n \tv \, , & n &\in \edge^p_{13} \, , &\tv &\in \left \{ \tv \in \tem_{13} \; | \; \ntr \tv \at_{e_{13}} = 0  \right \} \, , \\
	        & &\bm{\phi}(\xi, \eta) &= n \tv \, , & n &\in \edge^p_{23} \, , &\tv &\in \left \{ \tv \in \tem_{23} \; | \; \ntr \tv \at_{e_{23}} = 0  \right \} \, , \\
	        \text{Cell:}& &\bm{\phi}(\xi, \eta) &= n \tv \, , & n &\in \cell^p_{123} \, , &\tv &\in \tem_{123} \, ,
	    \end{align}
	    \end{subequations}
	    such that their normal trace vanishes on all edges.
	\end{itemize}
	\end{definition}
	A depiction of cubic base functions is given in \cref{fig:bdm}.
\begin{theorem} [Linear independence]
		The tensor product of the template with an $\Hone$-conforming polynomial basis $\U^p$ yields a unisolvent Brezzi-Douglas-Marini element.
	\end{theorem}
	\begin{proof}
		The proof is analogous to the one for the $\Nedtwo^p$-space \cref{th:li_tri_nedii}, since the Brezzi-Douglas-Marini element can be constructed by rotating the base functions of a N\'ed\'elec element of the second type by 90 degrees. 
	\end{proof}
	\begin{theorem} [$\Hd{,\surf}$-conformity]
		
		The constructed element on the reference domain is $\Hd{}$-conforming under contravariant Piola transformations of the base functions.
	\end{theorem}
	\begin{proof}
		By design, the non-cell base functions exhibit the property 
		\begin{align}
		    \langle \bm{\nu} , \, \bm{\phi}_i \rangle = n_i  \, ,
		\end{align}
		on their corresponding edges due to their construction by permutations of the reference triangle with contravariant Piola transformations. The property is maintained by contravariant Piola transformation to the physical domain. Consequently, the interface condition between neighbouring elements $\jump{\ntr\ud}|_{\Xi_{ij}} = 0$ is given by $\jump{\tr \langle \vb{n} , \, \ud \rangle}|_{\Xi_{ij}} = 0$ for the normal components, which is satisfied by the underlying subspace $\U^p$.
	\end{proof}
	
\subsection{N\'ed\'elec I} \label{sec:ned}
The N\'ed\'elec element of the second type has the disadvantage that its curl is a polynomial space of a lower degree $\rot{\Nedtwo^p} = \Po^{p-1}$. Consequently, one loses one order of convergence in the curl terms. In order to ameliorate the convergence rate one can employ the N\'ed\'elec elements of the first type \cite{Nedelec1980}, which enhance the polynomial space with base functions orthogonal to the kernel $\bm{\vartheta}_i \in \ker^\perp(\drot)$ 
\begin{align}
    &\dim \Ned^p(\Gamma) = \dim ([\Po^p(\Gamma)]^2 \oplus \bm{R}\widetilde{\Po}^p(\Gamma) \bm{\xi}) = (p+3)(p+1) \, , && p \geq 0 \, ,
\end{align}
such that $\Nedtwo^p \subset \Ned^p$.
The polynomial space $[\Po^p(\Gamma)]^2 \oplus \bm{R}\widetilde{\Po}^p(\Gamma) \bm{\xi}$ is used to construct the N\'ed\'elec base functions on the reference triangle.
In order to enhance the space, one must be able to split the space between kernel functions and non-kernel functions while maintaining conformity. 
Here we follow the ideas presented in \cite{Zaglmayr2006,Joachim2005}, where one explicitly applies the operators in the exact polynomial sequences to construct the kernel of the next space. We complement the kernel space with our new approach and introduce a specific and intuitive polytopal template leading to non-kernel base functions. 
	
\begin{figure}
    	\centering
    	\begin{subfigure}{0.3\linewidth}
    		\centering
    		\includegraphics[width=0.7\linewidth]{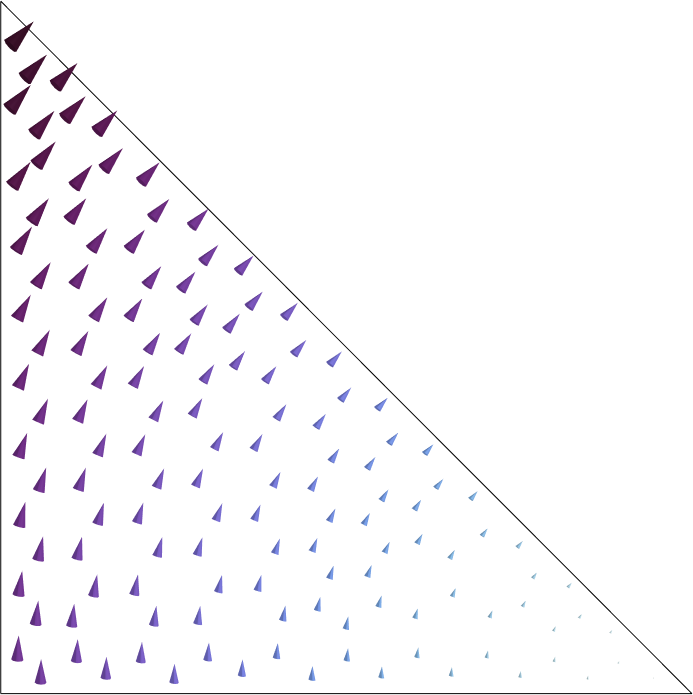}
    		\caption{}
    	\end{subfigure}
    	\begin{subfigure}{0.3\linewidth}
    		\centering
    		\includegraphics[width=0.7\linewidth]{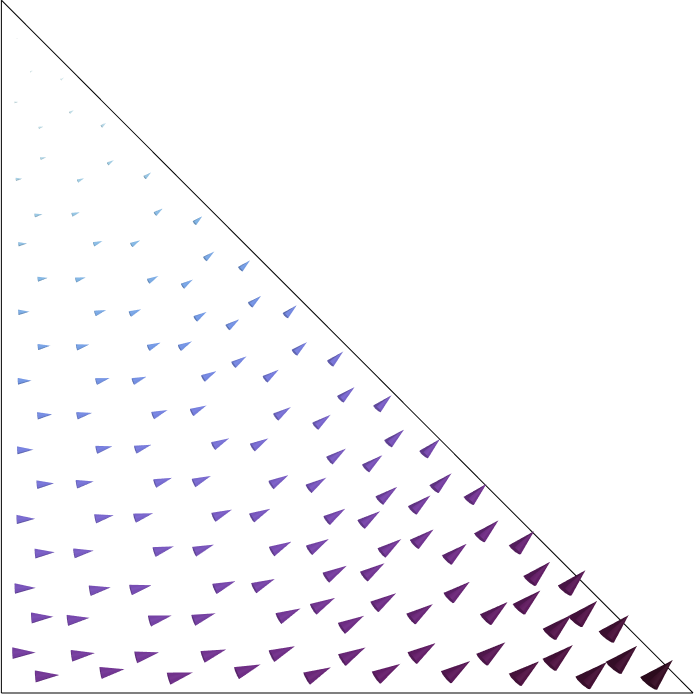}
    		\caption{}
    	\end{subfigure}
    	\begin{subfigure}{0.3\linewidth}
    		\centering
    		\includegraphics[width=0.7\linewidth]{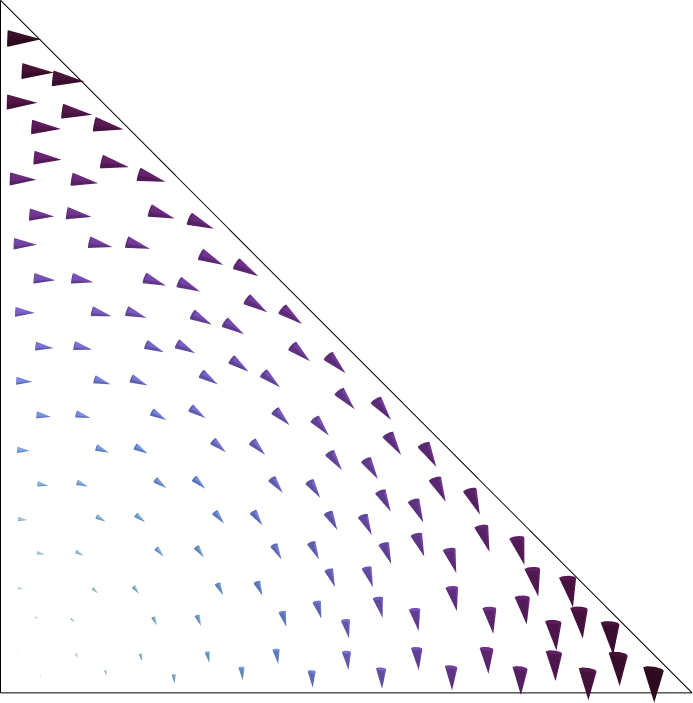}
    		\caption{}
    	\end{subfigure}
    	\caption{The base functions of the lowest order N\'ed\'elec element of the first type on the reference triangle belonging to the first (a), second (b) and third (c) edges.}
    	\label{fig:ned}
    \end{figure}

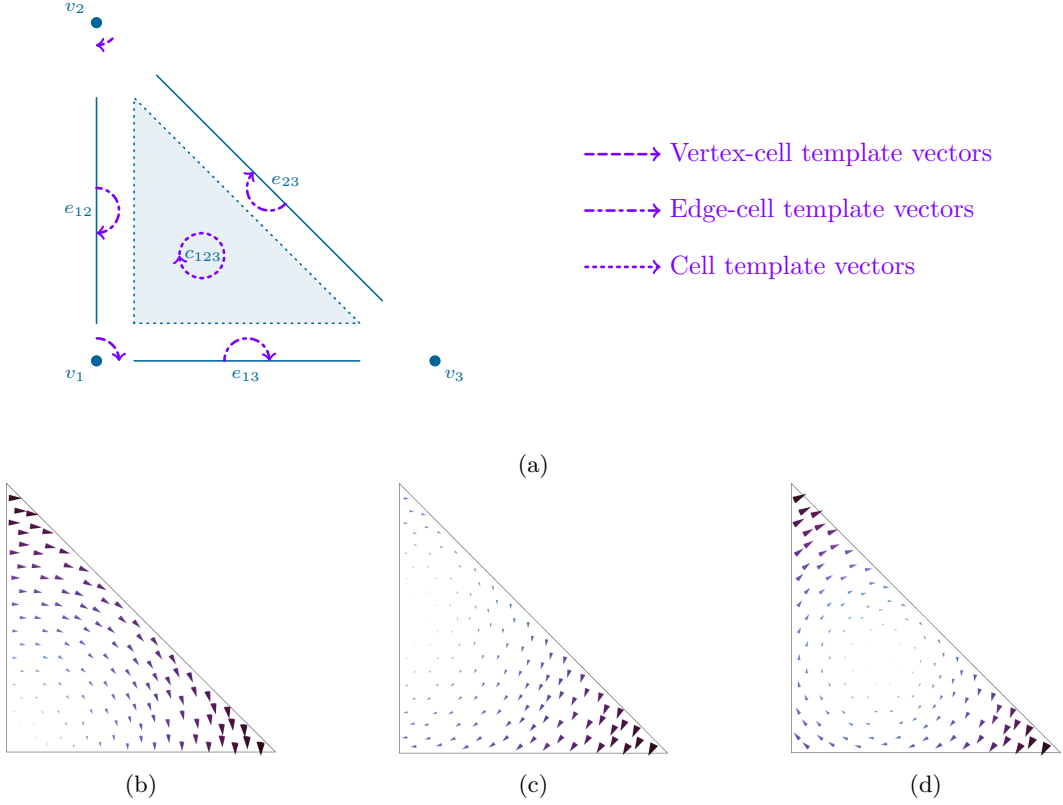
\begin{figure}
    	\centering
    	\begin{subfigure}{1\linewidth}
    		\centering
    		\definecolor{asl}{rgb}{0.4980392156862745,0.,1.}
    		\definecolor{asb}{rgb}{0.,0.4,0.6}
    		\begin{tikzpicture}[line cap=round,line join=round,>=triangle 45,x=1.0cm,y=1.0cm]
    			\clip(-2,-1.5) rectangle (12.5,4.5);
    			\draw (-0.5,-0.5) node[circle,fill=asb,inner sep=1.5pt] {};
    			\draw (-0.5,4) node[circle,fill=asb,inner sep=1.5pt] {};
    			\draw (4,-0.5) node[circle,fill=asb,inner sep=1.5pt] {};
    			\draw [color=asb,line width=.6pt] (-0.5,0) -- (-0.5,3);
    			\draw [color=asb,line width=.6pt] (0,-0.5) -- (3,-0.5);
    			\draw [color=asb,line width=.6pt] (0.3,3.3) -- (3.3,0.3);
    			\draw [dotted,color=asb,line width=.6pt] (0,0) -- (0,3) -- (3,0) -- (0,0);
    			\fill[opacity=0.1, asb] (0,0) -- (0,3) -- (3,0) -- cycle;
    			\draw (-0.5,-0.5) node[color=asb,anchor=north east] {$_{v_1}$};
    			
    			\draw [to-,asl,domain=180:270,line width=1.pt, densely dashed] plot ({0.3*cos(\x-180)-0.5}, {0.3*sin(\x-180)-0.5});
    			
    			\draw (4,-0.5) node[color=asb,anchor=north west] {$_{v_3}$};
    			\draw (-0.5,4) node[color=asb,anchor=south east] {$_{v_2}$};
    			
    			\draw [to-,asl,domain=90:135,line width=1.pt, densely dashed] plot ({0.3*cos(\x-180)-0.5}, {0.3*sin(\x-180)+4});
    			
    			\draw (-0.4,1.5) node[color=asb,anchor=east] {$_{e_{12}}$};
    			
    			\draw [to-,asl,domain=90:270,line width=1.pt, dashdotted] plot ({0.3*cos(\x-180)-0.5}, {0.3*sin(\x-180)+1.5});
    			
    			\draw (1.5,-.5) node[color=asb,anchor=north] {$_{e_{13}}$};
    			
    			\draw [to-,asl,domain=180:360,line width=1.pt, dashdotted] plot ({0.3*cos(\x-180)+1.5}, {0.3*sin(\x-180)-0.5});
    			
    			\draw (1.7,1.7) node[color=asb,anchor=south west] {$_{e_{23}}$};
    			
    			\draw [-to,asl,domain=135:-45,line width=1.pt, dashdotted] plot ({0.3*cos(\x-180)+1.8}, {0.3*sin(\x-180)+1.8});
    			
    			\draw (0.93,0.9)
    			node[color=asb] {$_{c_{123}}$};
    			
    			\draw [to-,asl,domain=0:360,line width=1.pt, dotted] plot ({0.3*cos(\x-180)+0.9}, {0.3*sin(\x-180)+0.9});
    			
    			\draw [-to,color=asl,line width=1pt, densely dashed] (6,2.25) -- (7,2.25);
    			\draw (7,2.25)
    			node[color=asl,anchor=west] {Vertex-cell template vectors};
    			\draw [-to,color=asl,line width=1pt, dashdotted] (6,1.5) -- (7,1.5);
    			\draw (7,1.5)
    			node[color=asl,anchor=west] {Edge-cell template vectors};
    			\draw [-to,color=asl,line width=1pt, dotted] (6,0.75) -- (7,0.75);
    			\draw (7,0.75)
    			node[color=asl,anchor=west] {Cell template vectors};
    		\end{tikzpicture}
    	\caption{}
    	\end{subfigure}
    	\begin{subfigure}{0.3\linewidth}
    		\centering
    		\includegraphics[width=0.7\linewidth]{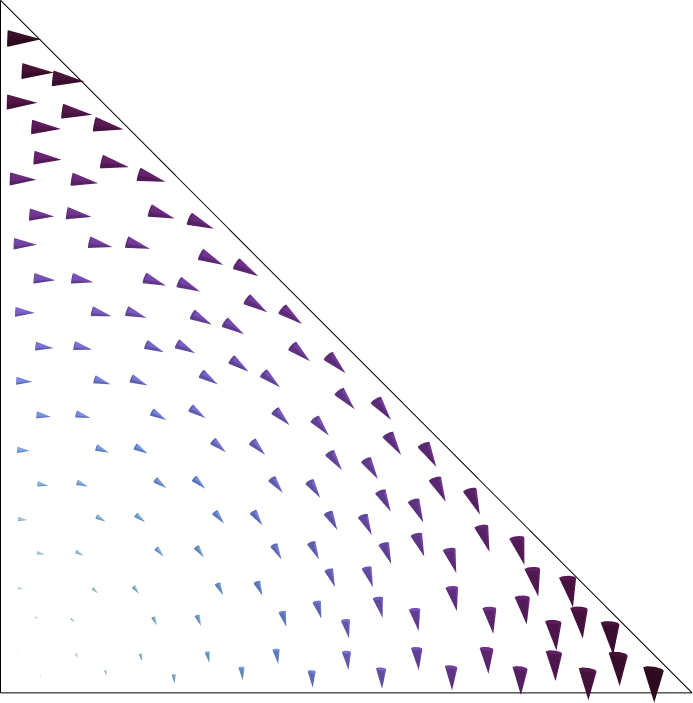}
    		\caption{}
    	\end{subfigure}
    	\begin{subfigure}{0.3\linewidth}
    		\centering
    		\includegraphics[width=0.7\linewidth]{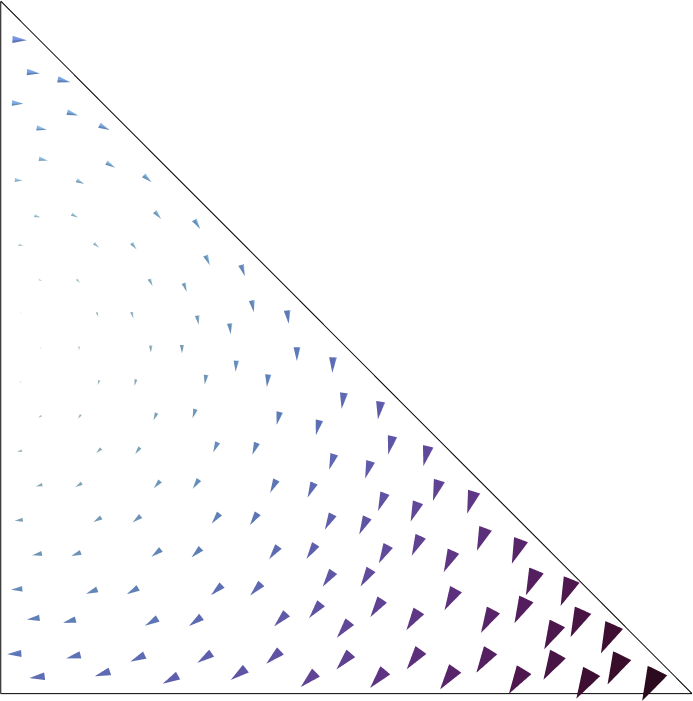}
    		\caption{}
    	\end{subfigure}
    	\begin{subfigure}{0.3\linewidth}
    		\centering
    		\includegraphics[width=0.7\linewidth]{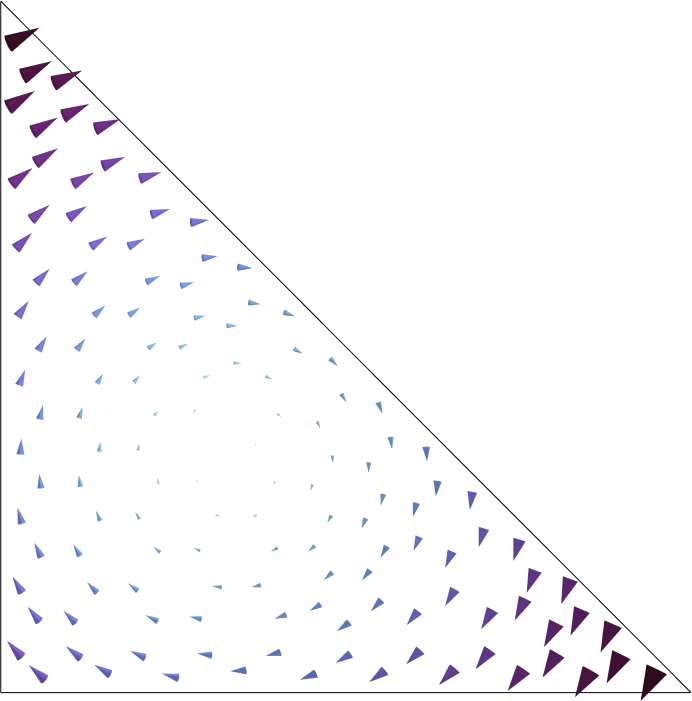}
    		\caption{}
    	\end{subfigure}
    	\caption{Template vectors on the reference triangle for base functions orthogonal to the kernel of the curl operator (a). The template of the first vertex (b), followed by the template of the first edge (c) and lastly, the cell template (d).}
    	\label{fig:rangetemptri}
    \end{figure}

\begin{figure}
		\centering
		\begin{subfigure}{0.3\linewidth}
			\centering
			\includegraphics[width=0.7\linewidth]{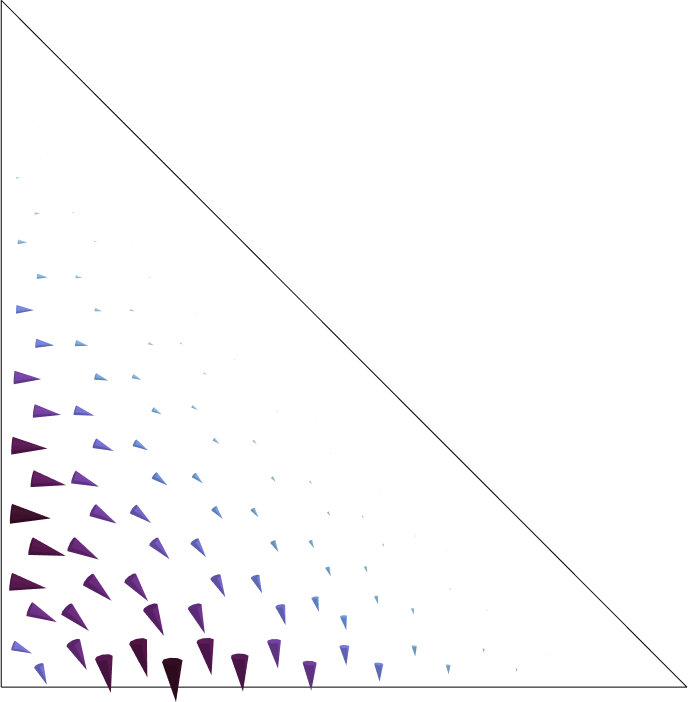}
			\caption{}
		\end{subfigure}
	    \begin{subfigure}{0.3\linewidth}
	    	\centering
	    	\includegraphics[width=0.7\linewidth]{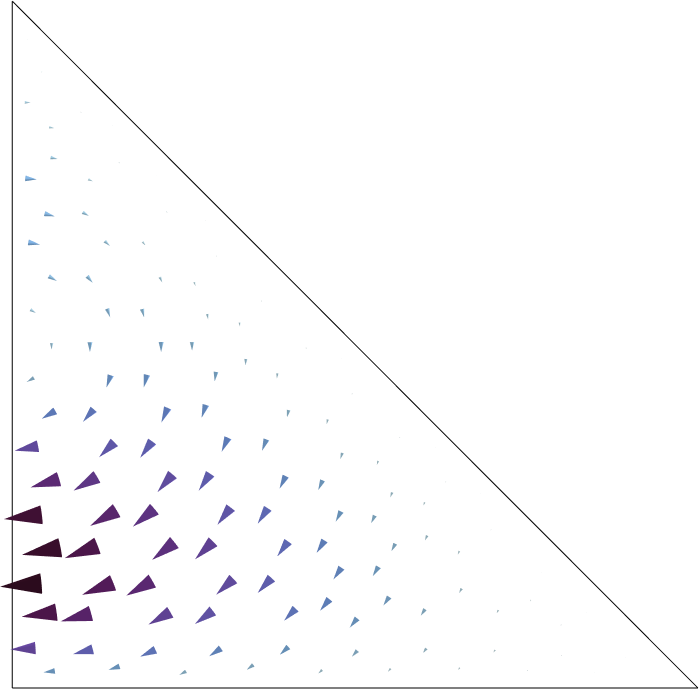}
	    	\caption{}
	    \end{subfigure}
        \begin{subfigure}{0.3\linewidth}
        	\centering
        	\includegraphics[width=0.7\linewidth]{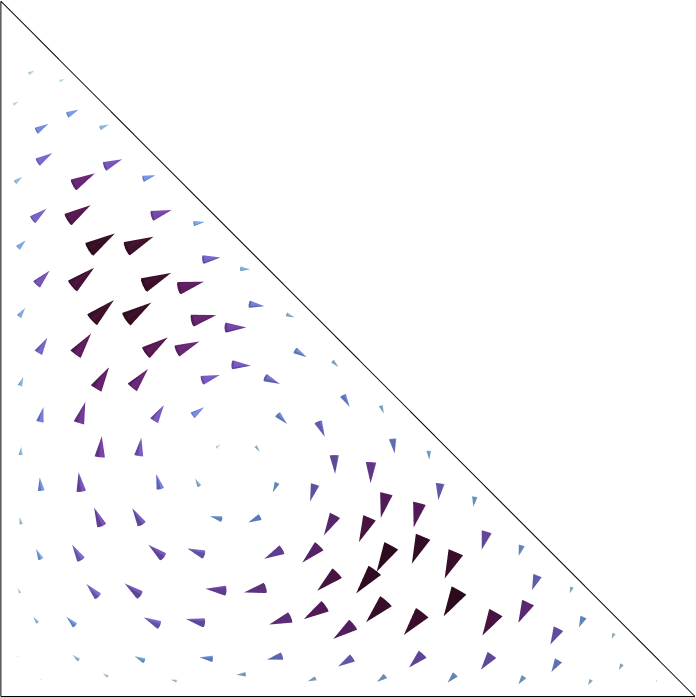}
        	\caption{}
        \end{subfigure}
		\begin{subfigure}{0.3\linewidth}
			\centering
			\includegraphics[width=0.7\linewidth]{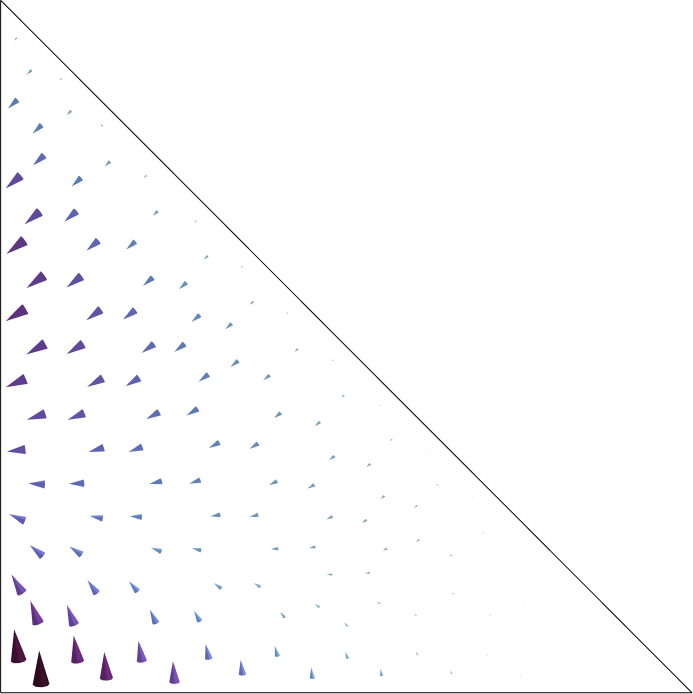}
			\caption{}
		\end{subfigure}
		\begin{subfigure}{0.3\linewidth}
			\centering
			\includegraphics[width=0.7\linewidth]{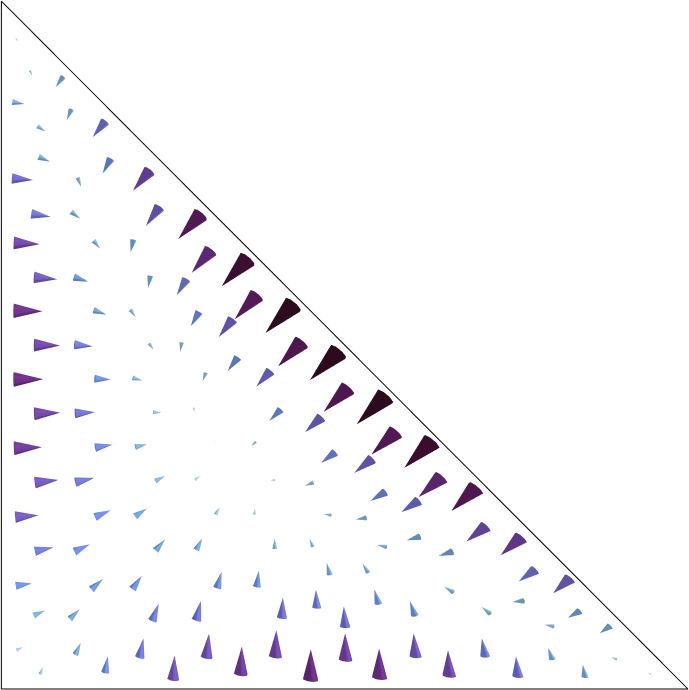}
			\caption{}
		\end{subfigure}
		\caption{Non-gradient vertex-cell (a), edge-cell (b) and pure cell (c) base functions. Gradient edge (d) and cell (e) base functions. 
		The base functions belong to the N\'ed\'elec element of the first type and are depicted on the reference triangle.}
		\label{fig:nedgrad}
	\end{figure}

We start with the kernel of $\Ned^p$ by taking gradients of the base functions of the $\U^{p+1}(\Gamma)$ space while excluding the vertex base functions
	\begin{align}
		&\bm{\vartheta}_{i}(\xi,\eta) = \nabla_\xi n_{i}^{p+1}  \, .
	\end{align}
This yields $(p+2)(p+1)/2 - 3$ base functions. We augment the space by adding the lowest order N\'ed\'elec ($\Ned^0$) base functions of the first type \cite{Anjam2015}, see \cref{fig:ned},
	\begin{align}
		&\bm{\vartheta}^I_1(\xi,\eta) = \begin{bmatrix} 
			\eta \\ 1 - \xi
		\end{bmatrix} \, , &&
		\bm{\vartheta}^I_2(\xi,\eta) = \begin{bmatrix} 
			1 - \eta \\ \xi
		\end{bmatrix} \, , &&
		\bm{\vartheta}^I_3(\xi,\eta) = \begin{bmatrix} 
			\eta \\ - \xi
		\end{bmatrix} \, .
		\label{eq:nedtri}
	\end{align}
Next we need to enhance the space with a minimal amount of base functions belonging to a higher degree polynomial space such that $\dim [\rot {\Ned^p}] = \dim\Po^p = (p+2)(p+1)/2$ base functions are found and the curl spans the next polynomial space in the sequence. 
In order to do so we introduce the polytopal template
	\begin{align}
		&\tem = \{\tem_{1}, \tem_{2}, \tem_{12} , \tem_{13}, \tem_{23}, \tem_{123}\} \, , 
	\end{align}
where the polytopal sets are derived from the lowest order N\'ed\'elec base functions of the first type (\cref{eq:nedtri})
    \begin{align}
    	&\tem_{1} = \{ \bm{\vartheta}_3^I \} \, , && \tem_{2} = \{ \bm{\vartheta}_2^I \} \, , && \tem_{12} = \{ \bm{\vartheta}_{3}^I - \bm{\vartheta}_{2}^I \} \, , \notag \\ 
    	& \tem_{13} = \{ \bm{\vartheta}_{1}^I + \bm{\vartheta}_{3}^I \} \, , && \tem_{23} = \{ \bm{\vartheta}_{1}^I - \bm{\vartheta}_{2}^I \} \, , && \tem_{123} = \{ \bm{\vartheta}_{1}^I - \bm{\vartheta}_{2}^I + \bm{\vartheta}_{3}^I \} \, .
    \end{align}
From the depiction in \cref{fig:rangetemptri} it is intuitively apparent that the template vectors represent the components needed to generate rotational flux. In fact, the template is composed of shifted and scaled linear vortex fields. 
The non-gradient base functions are given by the tensor product 
	\begin{align}
		&\left \{ \bigoplus_{i = 1}^2 \ver^{p}_i \otimes \tem_i  \right \} \oplus \left \{ \bigoplus_{j \in \mathcal{J}} \edge^{p}_j \otimes \tem_j  \right \} \oplus \{ \cell_{123}^{p} \otimes \tem_{123} \} \, ,  && \mathcal{J} = \{(1,2),(1,3),(2,3)\} \, . 
	\end{align}
This generates exactly $2p+(p-1)p/2$ base functions. Adding the constants from $\dim [\rot {\Ned^0}] = 1$ satisfies the dimensionality of the next polynomial space in the sequence $2p+(p-1)p/2 + 1 = (p+2)(p+1)/2 = \dim \Po^p$.
Further, since $\bm{\vartheta}_i^I$ belong to $[\Po^0]^2\oplus \bm{R} \widetilde{\Po}^0 \bm{\xi}$ the resulting base functions $n^p_i \bm{\vartheta}_j^I$ clearly belong to $[\Po^p]^2\oplus \bm{R} \widetilde{\Po}^p \bm{\xi}$.
The complete N\'ed\'elec space reads
	\begin{align}
		&\Ned^p = \Ned^0 \oplus  \left \{ \bigoplus_{j \in \mathcal{J} } \nabla \edge^{p+1}_j \right \} \oplus \nabla \cell^{p+1}_{123} \oplus \left \{ \bigoplus_{i = 1}^2 \ver^{p}_i \otimes \tem_i  \right \} \oplus \left \{ \bigoplus_{j \in \mathcal{J}} \edge^{p}_j \otimes \tem_j  \right \} \oplus \{ \cell^{p}_{123} \otimes \tem_{123} \} \, , \notag \\ & \mathcal{J} = \{(1,2),(1,3),(2,3)\} \, .
	\end{align}
Using the template and the exact sequence we can state the base functions for the N\'ed\'elec element of the first type.
\begin{definition}[Triangle $\Ned^p$ base functions]
	The base functions of the triangle N\'ed\'elec element of the first type are defined on their respective polytope as follows.
	\begin{itemize}
	    \item On each edge $e_{ij}$ the base functions read
	    \begin{subequations}
	        \begin{align}
	        \text{Edge:}&&\bm{\vartheta}(\xi, \eta) &=  \bm{\vartheta}^I \, ,   &\bm{\vartheta}^I &\in \left \{ \bm{\vartheta} \in \Ned^0 \; | \; \rtr \bm{\vartheta} \at_{e_{ij}}  \neq 0 \right \} \, , \\
	        &&\bm{\vartheta}(\xi, \eta) &= \nabla_\xi n \, , & n &\in \edge^{p+1}_{ij} \, , 
	    \end{align}
	    \end{subequations}
	    \item The cell base functions are given by 
	    \begin{subequations}
	        \begin{align}
	        \text{Vertex-cell:}&&\bm{\vartheta}(\xi, \eta) &= n \tv \, , & n &\in \ver^p_{1} \, , &\tv &\in \tem_{1} \, , \\
	        &&\bm{\vartheta}(\xi, \eta) &= n \tv \, , & n &\in \ver^p_{2} \, , &\tv &\in \tem_{2} \, , \\
	        \text{Edge-cell:}&&\bm{\vartheta}(\xi, \eta) &= n \tv \, , & n &\in \edge^p_{12} \, , &\tv &\in \tem_{12} \, , \\
	        &&\bm{\vartheta}(\xi, \eta) &= n \tv \, , & n &\in \edge^p_{13} \, , &\tv &\in \tem_{13} \, , \\
	        &&\bm{\vartheta}(\xi, \eta) &= n \tv \, , & n &\in \edge^p_{23} \, , &\tv &\in \tem_{23} \, , \\
	        \text{Cell:}&&\bm{\vartheta}(\xi, \eta) &= n \tv \, , & n &\in \cell^p_{123} \, , &\tv &\in \tem_{123} \, , \\
	        &&\bm{\vartheta}(\xi, \eta) &= \nabla_\xi n  \, , & n &\in \cell^{p+1}_{123} \, , 
	    \end{align}
	    \end{subequations}
	    such that their tangential trace vanishes on all edges. This holds true for the cell gradients due to the exactness of the de Rham sequence in \cref{fig:derham0}.
	\end{itemize}
	\end{definition}
	A depiction of the higher order base functions is given in \cref{fig:nedgrad}.
    \begin{theorem} [Linear independence]
    \label{th:linin}
		The set of base functions given by the lowest order N\'ed\'elec elements, the gradients of an $\Hone$-conforming polynomial subspace $\U^{p+1}$ (excluding vertex base functions) and the tensor product of the $\U^{p}$ base functions with the polytopal template yield a linearly independent polynomial basis for $\Ned^p$.
	\end{theorem}
	\begin{proof}
		We start by showing the gradients of the $\U^p$ base functions are linearly independent of each other by using contradiction.
		Assume the set of gradients is linearly dependent, then there holds 
		\begin{align}
			\sum_i c_i \nabla_\xi n_i = \nabla_\xi \sum_i c_i n_i = 0 \, ,
		\end{align}
		for some combination of constants $c_i$ where not all $c_i$ values are zero. However, the vertex base functions are not employed. Thus, if the basis satisfies the partition of unity property, then the kernel of the gradient operator, namely constants $\R$, is missing and the exact sequence property yields a contradiction (see \cref{fig:derham}).
		The same holds true for a hierarchical polynomial basis, since the vertex base functions are used to capture constants.
		
		The base functions of the lowest order N\'ed\'elec elements of the first type are linearly independent of the gradients since their tangential traces on the edges of the triangle are constant $\rtr \bm{\vartheta}^I_i |_{\mu_i} \in \R$, whereas the tangential traces of the edge gradients are at least linear and the tangential traces of the cell gradients vanish on the entire boundary by the exact sequence property, compare \cref{fig:derham0}.
		Together, the three lowest order base functions span the constant space $[\Po^0]^2 = \R^2$. Further, their curls span the constant space $\R$. 
		
		In order to complete the proof we must show that the remaining base functions are non-gradients and linearly independent of the lowest order base functions.
		Observe that the template vectors have the general form
		\begin{align}
			&\bm{\vartheta}^I = \begin{bmatrix}
						c_1 \eta - c_2 \\ c_3 - c_1 \xi 
			\end{bmatrix}  \, , && c_1 \in  \{1,2,3\}  \, , &&  c_2,  c_3 \in \{0,1\} \, ,
		\end{align}
		such that the curl of the base function reads
		\begin{align}
			\mathrm{div}_\xi (\bm{R} n \bm{\vartheta}^I) &= \langle \nabla_\xi n , \, \bm{R} \bm{\vartheta}^I \rangle - 2 c_1  n = (c_3 - c_1 \xi)n_{,\xi} + (c_2 - c_1 \eta) n_{,\eta} - 2c_1 n \, .
		\end{align}
	    Clearly, the polynomial order of the underlying scalar base function $n$ is maintained under the curl operator.
	    Therefore, there holds
	    \begin{align}
        	& \mathrm{div}_\xi(\bm{R}[c_1 \bm{R} \bm{\xi} + \vb{c}]\Po^p(\Gamma)) = \Po^p(\Gamma) \, , && \vb{c} = \begin{bmatrix}
        		-c_2 \\ c_3 
        	\end{bmatrix} \, ,
        \end{align}
        such that
	    \begin{align}
	         &\mathrm{div}_\xi (\bm{R}\sum_{i=0}^{p} c_i n_i \bm{\vartheta}^I)   = 0 \quad \iff \quad c_i = 0 \quad \forall \, i \in \{0,1,...,p\} \, ,
	    \end{align}
	    by dimensionality. The latter is readily observed when a hierarchical polynomial basis is used.
	    If the complete span of $n_i$ functions is employed, then the constant element is present and the curl operator maps also to the space of constants $\R$. However, this space is already obtained by employing the $\Ned^0$-basis in the construction and as such, leads to linear dependence. By removing the last vertex base function, the hierarchical basis no longer contains the space of constants, thus asserting linear independence of the total construction. The same holds true for any other basis that satisfies the partition of unity property since removing one base function cancels this characteristic and removes the constant element from the space.
	\end{proof}
    \begin{theorem} [$\Hc{,\surf}$-conformity]
    	
    	The resulting finite element is $\Hc{}$-conforming under covariant Piola mappings.
    \end{theorem}
    \begin{proof} 
    	The lowest order N\'ed\'elec base functions are $\Hr{}$-conforming by their degrees of freedom. The gradient base functions are conforming due to the exact sequence property $\nabla \Hone \subset \Hr{}$, see \cref{fig:derham}. Lastly, the remaining base functions are cell-type and do not affect the conformity of the finite element. 
 
    	The covariant Piola transformations maintain the conformity of the base functions across the mapping from the reference to the physical element.
    \end{proof}
	\begin{remark}
	        By employing the polynomial spaces $\ver_i^{p-1}$, $\edge_{ij}^{p-1}$ and $\cell_{123}^{p-1}$ in the construction of the non-gradient base functions, one finds the $\Nedtwo^p$-element with a split between the kernel and non-kernel functions.
	\end{remark}

\subsection{Raviart-Thomas}
The Raviart-Thomas triangle element \cite{Raviart} enhances the Brezzi-Douglas-Marini triangle element by adding base functions in $\ker^\perp(\mathrm{div})$ such that $\di \RT^p = \Po^p$, thus improving convergence estimates for divergence terms.
In order to construct Raviart-Thomas elements, one must split the polynomial basis between solenoidal and non-curl base functions. The construction follows the same lines as in \cref{sec:ned}.  

\begin{figure}
		\centering
		\begin{subfigure}{0.3\linewidth}
			\centering
			\includegraphics[width=0.7\linewidth]{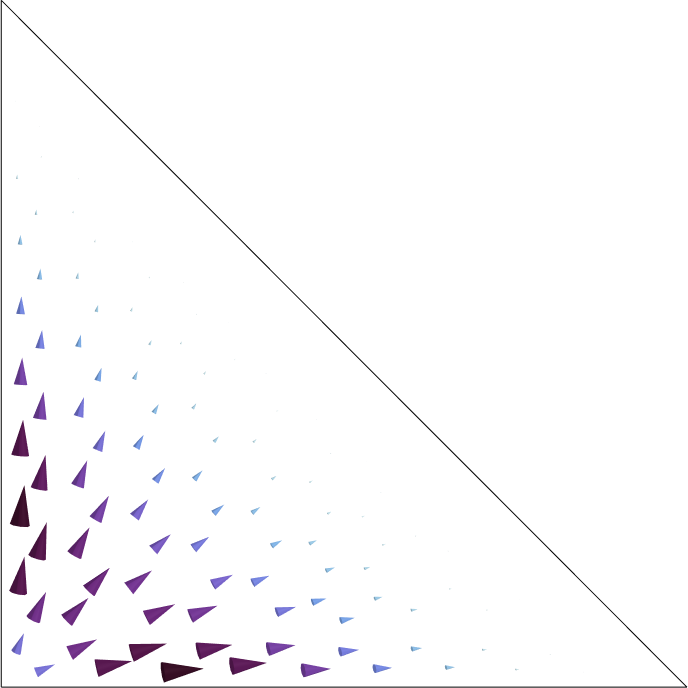}
			\caption{}
		\end{subfigure}
	    \begin{subfigure}{0.3\linewidth}
	    	\centering
	    	\includegraphics[width=0.7\linewidth]{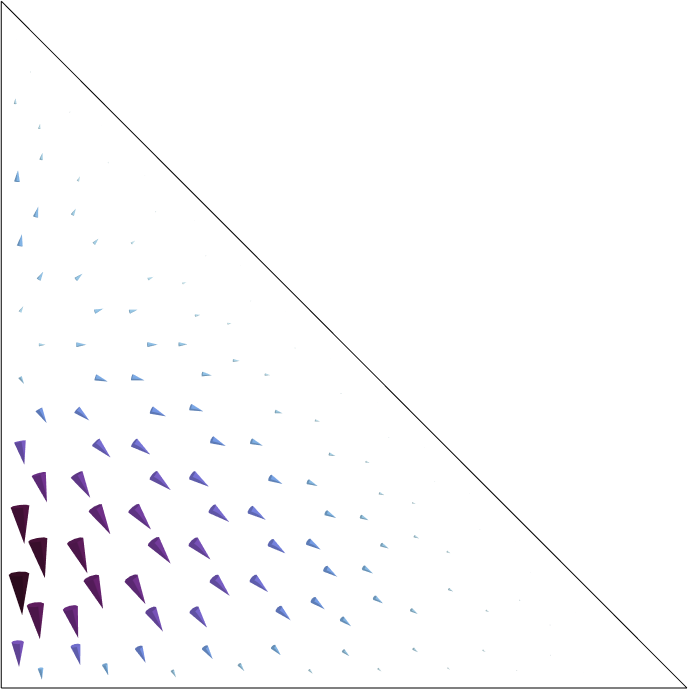}
	    	\caption{}
	    \end{subfigure}
        \begin{subfigure}{0.3\linewidth}
        	\centering
        	\includegraphics[width=0.7\linewidth]{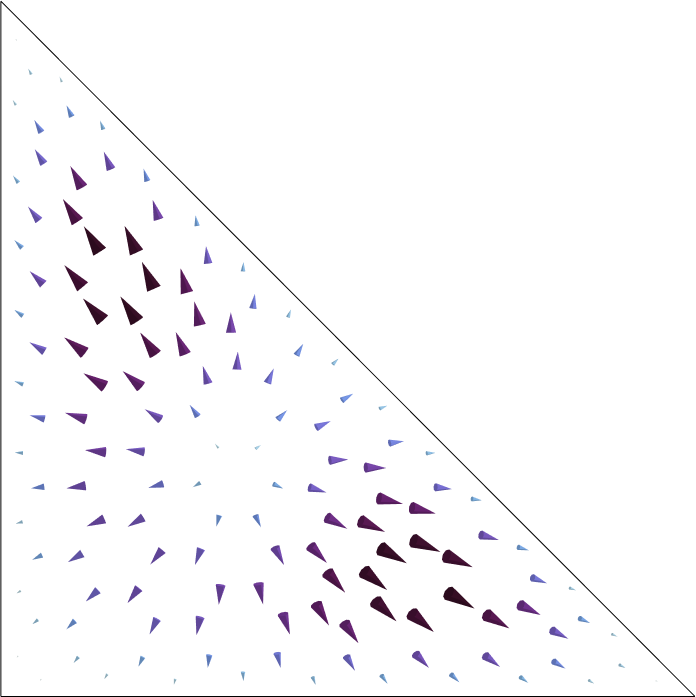}
        	\caption{}
        \end{subfigure}
		\begin{subfigure}{0.3\linewidth}
			\centering
			\includegraphics[width=0.7\linewidth]{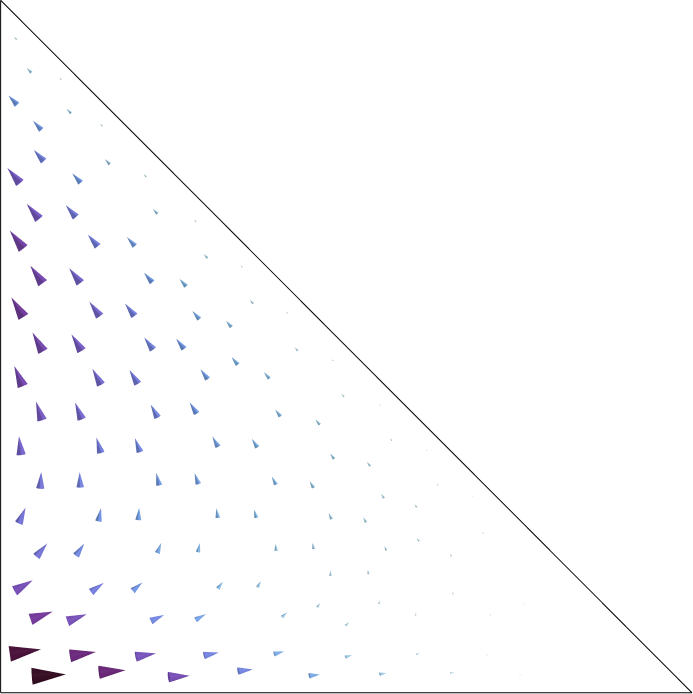}
			\caption{}
		\end{subfigure}
		\begin{subfigure}{0.3\linewidth}
			\centering
			\includegraphics[width=0.7\linewidth]{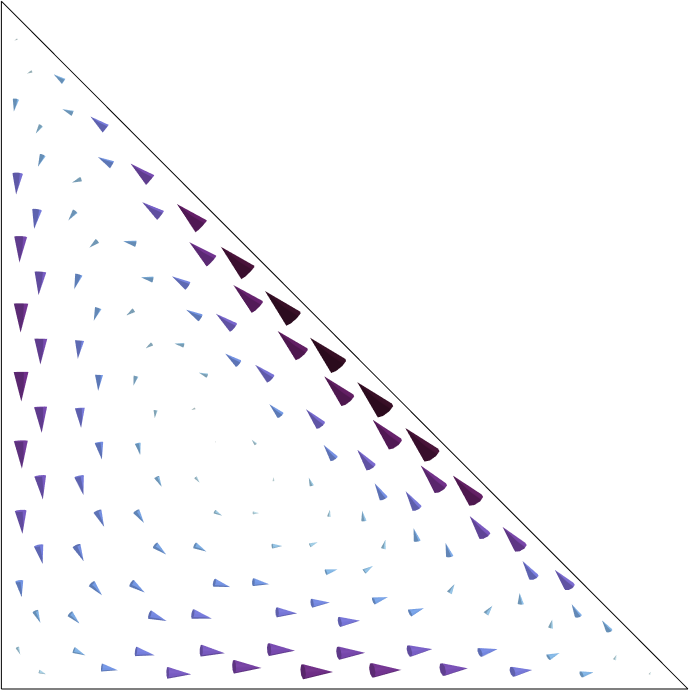}
			\caption{}
		\end{subfigure}
		\caption{Non-curl vertex-cell (a), edge-cell (b) and pure cell (c) base functions. Curl edge (d) and cell (e) base functions. The base functions belong to the triangle Raviart-Thomas element.}
		\label{fig:rtgrad}
	\end{figure}

The kernel of the space is partially given by curls of the $\U^{p+1}(\Gamma)$ space without vertex base functions
\begin{align}
    \bm{\phi}_i(\xi,\eta) = \rog_\xi n_i^{p+1} \, ,
\end{align}
amounting to $(p+2)(p+1)/2 - 3$ base functions. Next the lowest order Raviart-Thomas ($\RT^0$) base functions \cite{Anjam2015} are added 
\begin{align}
    &\bm{\phi}_1^I = \begin{bmatrix} 
    1 - \xi \\ -\eta
    \end{bmatrix} \, , &&
    \bm{\phi}_2^I = \begin{bmatrix} 
    \xi \\ \eta-1
    \end{bmatrix} \, &&
    \bm{\phi}_3^I = \begin{bmatrix} 
    - \xi \\ -\eta
    \end{bmatrix} \, ,
\end{align}
accounting for constant fields in $\R^2$ and mapping to constant fields in $\R$ via the divergence operator, $\di \RT^0 = \R$.
We complete the space by introducing polytopal template sets for non-curl base functions
\begin{align}
    	&\tem_{1} = \{ -\bm{\phi}_3^I \} \, , && \tem_{2} = \{ \bm{\phi}_2^I \} \, , && \tem_{12} = \{ \bm{\phi}_{2}^I - \bm{\phi}_{3}^I \} \, , \notag \\ 
    	& \tem_{13} = \{ -\bm{\phi}_{1}^I - \bm{\phi}_{3}^I \} \, , && \tem_{23} = \{ \bm{\phi}_{2}^I - \bm{\phi}_{1}^I \} \, , && \tem_{123} = \{ \bm{\phi}_{2}^I -  \bm{\phi}_{1}^I - \bm{\phi}_{3}^I \} \, ,
    \end{align}
such that the entire template reads
\begin{align}
		&\tem = \{\tem_{1}, \tem_{2}, \tem_{12} , \tem_{13}, \tem_{23}, \tem_{123}\} \, . 
	\end{align}
With the template at hand we can construct the Raviart-Thomas element 
\begin{align}
		&\RT^p = \RT^0 \oplus  \left \{ \bigoplus_{j \in \mathcal{J} } \rog \edge^{p+1}_j \right \} \oplus \rog \cell^{p+1}_{123} \oplus \left \{ \bigoplus_{i = 1}^2 \ver^{p}_i \otimes \tem_i  \right \} \oplus \left \{ \bigoplus_{j \in \mathcal{J}} \edge^{p}_j \otimes \tem_j  \right \} \oplus \{ \cell^{p}_{123} \otimes \tem_{123} \} \, , \notag \\ & \mathcal{J} = \{(1,2),(1,3),(2,3)\} \, .
	\end{align}
\begin{definition}[Triangle $\RT^p$ base functions]
	The base functions of the triangle Raviart-Thomas element are defined on their respective polytope as follows.
	\begin{itemize}
	    \item On each edge $e_{ij}$ the base functions read
	    \begin{subequations}
	        \begin{align}
	        \text{Edge:}&&\bm{\phi}(\xi, \eta) &=  \bm{\phi}^I \, ,   &\bm{\phi}^I &\in \left \{ \bm{\phi} \in \RT^0 \; | \; \ntr \bm{\phi} \at_{e_{ij}}  \neq 0 \right \} \, , \\
	        &&\bm{\phi}(\xi, \eta) &= \rog_\xi n \, , & n &\in \edge^{p+1}_{ij} \, , 
	    \end{align}
	    \end{subequations}
	    \item The cell base functions are given by 
	    \begin{subequations}
	        \begin{align}
	        \text{Vertex-cell:}&&\bm{\phi}(\xi, \eta) &= n \tv \, , & n &\in \ver^p_{1} \, , &\tv &\in \tem_{1} \, , \\
	        &&\bm{\phi}(\xi, \eta) &= n \tv \, , & n &\in \ver^p_{2} \, , &\tv &\in \tem_{2} \, , \\
	        \text{Edge-cell:}&&\bm{\phi}(\xi, \eta) &= n \tv \, , & n &\in \edge^p_{12} \, , &\tv &\in \tem_{12} \, , \\
	        &&\bm{\phi}(\xi, \eta) &= n \tv \, , & n &\in \edge^p_{13} \, , &\tv &\in \tem_{13} \, , \\
	        &&\bm{\phi}(\xi, \eta) &= n \tv \, , & n &\in \edge^p_{23} \, , &\tv &\in \tem_{23} \, , \\
	        \text{Cell:}&&\bm{\phi}(\xi, \eta) &= n \tv \, , & n &\in \cell^p_{123} \, , &\tv &\in \tem_{123} \, , \\
	        &&\bm{\phi}(\xi, \eta) &= \rog_\xi n  \, , & n &\in \cell^{p+1}_{123} \, , 
	    \end{align}
	    \end{subequations}
	    such that their tangential trace vanishes on all edges.
	\end{itemize}
	\end{definition}
	Several cubic base functions are depicted in \cref{fig:rtgrad}.
    \begin{theorem} [Linear independence]
		
		The set of base functions given by the lowest order Raviart-Thomas element, the curls of an $\Hone$-conforming polynomial subspace $\U^{p+1}$ (excluding vertex base functions) and the tensor product of the $\U^{p}$ base functions with the polytopal template yield a linearly independent polynomial basis for $\RT^p$.
	\end{theorem}
	\begin{proof}
	The proof follows the same lines as in \cref{th:linin}.
	\end{proof}
	\begin{theorem} [$\Hd{,\surf}$-conformity]
    	
    	The resulting finite element is $\Hd{}$-conforming under contravariant Piola mappings.
    \end{theorem}
    \begin{proof} 
    	The lowest order Raviart-Thomas base functions are $\Hd{}$-conforming due to their degrees of freedom. The curl base functions are conforming by exact sequence property $\rog \Hone \subset \Hd{}$, as depicted in \cref{fig:derham}. Finally, the remaining base functions are cell-type  and as such, do not affect the conformity of the finite element. 
 
    	The contravariant Piola transformations maintain the conformity of the base functions across the mapping from the reference to the physical element.
    \end{proof}
    
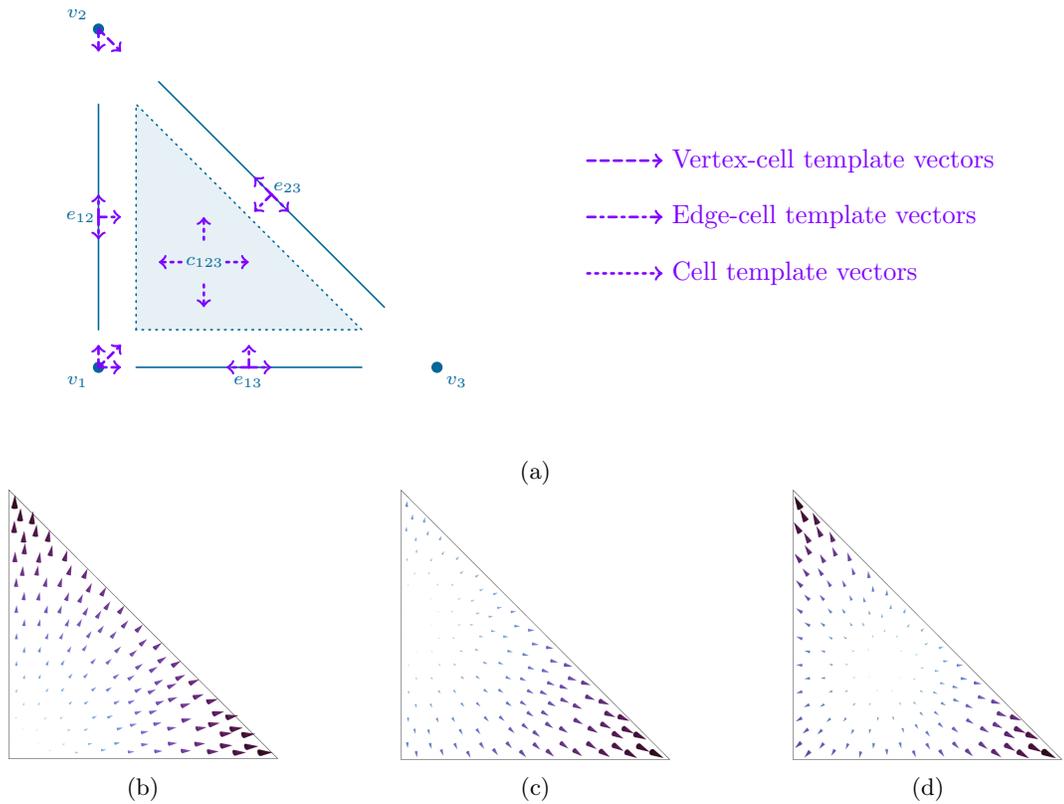
\begin{figure}
    	\centering
    	\begin{subfigure}{1\linewidth}
    		\centering
    		\definecolor{asl}{rgb}{0.4980392156862745,0.,1.}
    		\definecolor{asb}{rgb}{0.,0.4,0.6}
    		\begin{tikzpicture}[line cap=round,line join=round,>=triangle 45,x=1.0cm,y=1.0cm]
    			\clip(-2,-1.5) rectangle (12.5,4.5);
    			\draw (-0.5,-0.5) node[circle,fill=asb,inner sep=1.5pt] {};
    			\draw (-0.5,4) node[circle,fill=asb,inner sep=1.5pt] {};
    			\draw (4,-0.5) node[circle,fill=asb,inner sep=1.5pt] {};
    			\draw [color=asb,line width=.6pt] (-0.5,0) -- (-0.5,3);
    			\draw [color=asb,line width=.6pt] (0,-0.5) -- (3,-0.5);
    			\draw [color=asb,line width=.6pt] (0.3,3.3) -- (3.3,0.3);
    			\draw [dotted,color=asb,line width=.6pt] (0,0) -- (0,3) -- (3,0) -- (0,0);
    			\fill[opacity=0.1, asb] (0,0) -- (0,3) -- (3,0) -- cycle;
    			\draw (-0.5,-0.5) node[color=asb,anchor=north east] {$_{v_1}$};
    			
    			\draw [-to,color=asl,line width=1pt, densely dashed] (-0.5,-0.5) -- (-0.2,-0.2);
    			\draw [-to,color=asl,line width=1pt, densely dashed] (-0.5,-0.5) -- (-0.5,-0.2);
    			\draw [-to,color=asl,line width=1pt, densely dashed] (-0.5,-0.5) -- (-0.2,-0.5);
    			
    			\draw (4,-0.5) node[color=asb,anchor=north west] {$_{v_3}$};
    			\draw (-0.5,4) node[color=asb,anchor=south east] {$_{v_2}$};
    			
    			\draw [-to,color=asl,line width=1pt, densely dashed] (-0.5,4) -- (-0.2,3.7);
    			\draw [-to,color=asl,line width=1pt, densely dashed] (-0.5,4) -- (-0.5,3.7);
    			
    			\draw (-0.4,1.5) node[color=asb,anchor=east] {$_{e_{12}}$};
    			
    			\draw [-to,color=asl,line width=1pt, dashdotted] (-0.5,1.5) -- (-0.5,1.8);
    			\draw [-to,color=asl,line width=1pt, dashdotted] (-0.5,1.5) -- (-0.2,1.5);
    			\draw [-to,color=asl,line width=1pt, dashdotted] (-0.5,1.5) -- (-0.5,1.2);
    			
    			\draw (1.5,-.5) node[color=asb,anchor=north] {$_{e_{13}}$};
    			
    			\draw [-to,color=asl,line width=1pt, dashdotted] (1.5,-0.5) -- (1.5,-0.2);
    			\draw [-to,color=asl,line width=1pt, dashdotted] (1.5,-0.5) -- (1.8,-0.5);
    			\draw [-to,color=asl,line width=1pt, dashdotted] (1.5,-0.5) -- (1.2,-0.5);
    			
    			\draw (1.7,1.7) node[color=asb,anchor=south west] {$_{e_{23}}$};
    			
    			\draw [-to,color=asl,line width=1pt, dashdotted] (1.8,1.8) -- (1.57,1.57);
    			\draw [-to,color=asl,line width=1pt, dashdotted] (1.8,1.8) -- (2.03,1.57);
    			\draw [-to,color=asl,line width=1pt, dashdotted] (1.8,1.8) -- (1.57,2.03);
    			
    			\draw (0.92,0.9)
    			node[color=asb] {$_{c_{123}}$};
    			
    			\draw [-to,color=asl,line width=1pt, dotted] (1.2,0.9) -- (1.5,0.9);
    			\draw [-to,color=asl,line width=1pt, dotted] (0.6,0.9) -- (0.3,0.9);
    			\draw [-to,color=asl,line width=1pt, dotted] (0.9,1.2) -- (0.9,1.5);
    			\draw [-to,color=asl,line width=1pt, dotted] (0.9,0.6) -- (0.9,0.3);
    			
    			\draw [-to,color=asl,line width=1pt, densely dashed] (6,2.25) -- (7,2.25);
    			\draw (7,2.25)
    			node[color=asl,anchor=west] {Vertex-cell template vectors};
    			\draw [-to,color=asl,line width=1pt, dashdotted] (6,1.5) -- (7,1.5);
    			\draw (7,1.5)
    			node[color=asl,anchor=west] {Edge-cell template vectors};
    			\draw [-to,color=asl,line width=1pt, dotted] (6,0.75) -- (7,0.75);
    			\draw (7,0.75)
    			node[color=asl,anchor=west] {Cell template vectors};
    		\end{tikzpicture}
    	\caption{}
    	\end{subfigure}
    	\begin{subfigure}{0.3\linewidth}
    		\centering
    		\includegraphics[width=0.7\linewidth]{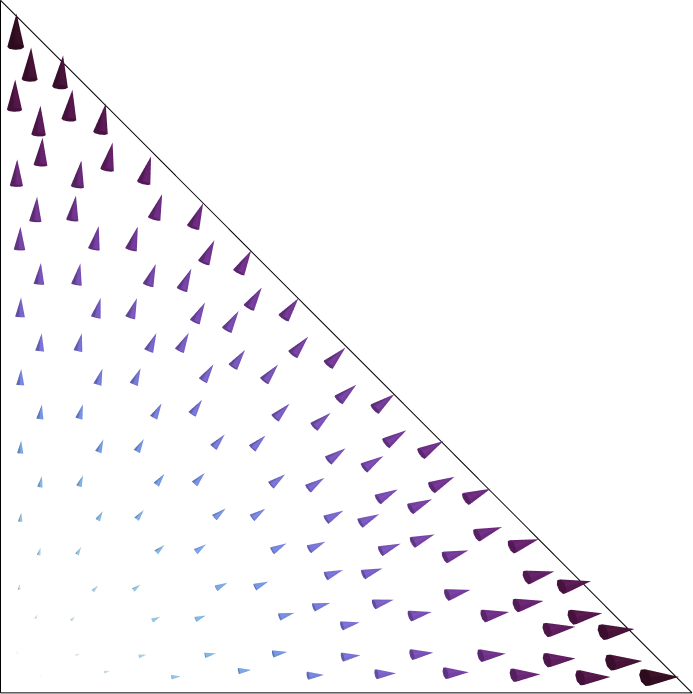}
    		\caption{}
    	\end{subfigure}
    	\begin{subfigure}{0.3\linewidth}
    		\centering
    		\includegraphics[width=0.7\linewidth]{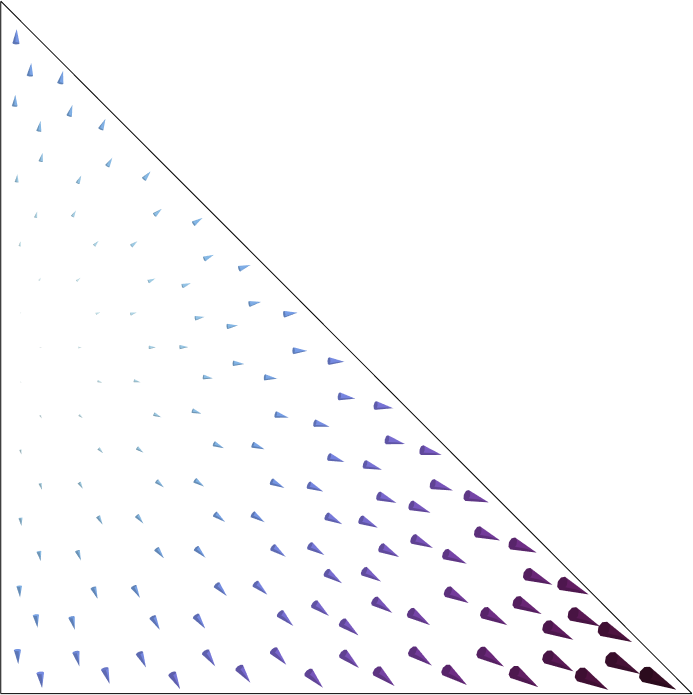}
    		\caption{}
    	\end{subfigure}
    	\begin{subfigure}{0.3\linewidth}
    		\centering
    		\includegraphics[width=0.7\linewidth]{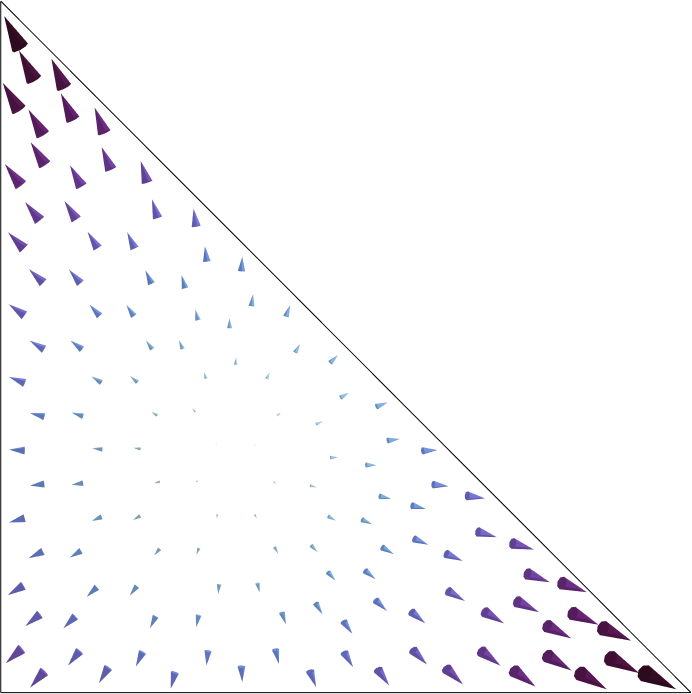}
    		\caption{}
    	\end{subfigure}
    	\caption{Template vectors on the reference triangle for base functions orthogonal to the kernel of the div operator (a). The template of the first vertex (b), followed by the template of the first edge (c) and lastly, the cell template (d).}
    	\label{fig:rt_tri}
    \end{figure} 
	\begin{remark}
	        By employing the polynomial spaces $\ver_i^{p-1}$, $\edge_{ij}^{p-1}$ and $\cell_{123}^{p-1}$ in the construction of the non-curl base functions one finds the $\BDM^p$-element with a split between the kernel and non-kernel functions.
	\end{remark}

\section{Three-dimensional templates}
\begin{figure}
    \centering
    \definecolor{asl}{rgb}{0.4980392156862745,0.,1.}
		\definecolor{asb}{rgb}{0.,0.4,0.6}
		\begin{tikzpicture}
			\begin{axis}
				[
				width=30cm,height=25cm,
				view={50}{15},
				enlargelimits=true,
				xmin=-1,xmax=2,
				ymin=-1,ymax=2,
				zmin=-1,zmax=2,
				domain=-10:10,
				axis equal,
				hide axis
				]
				\draw (-0.2, -0.2, -0.2) node[circle,fill=asb,inner sep=1.5pt] {};
				\draw (-0.2, -0.2, 1.2) node[circle,fill=asb,inner sep=1.5pt] {};
				\draw (1.2, -0.2, -0.2) node[circle,fill=asb,inner sep=1.5pt] {};
				\draw (-0.2, 1.2, -0.2) node[circle,fill=asb,inner sep=1.5pt] {};
				
				\addplot3[color=asb][line width=1pt,mark=o]
				coordinates {(0.1, 0.1, 0.1)};
				
				\addplot3[color=asb][line width=0.6pt,dotted]
				coordinates {(0,0,0)(0.5,0,0)(0,0.5,0)(0,0,0)};
				\addplot3[color=asb][line width=0.6pt,dotted]
				coordinates {(0,0,0)(0,0,0.5)};
				\addplot3[color=asb][line width=0.6pt,dotted]coordinates {(0.5,0,0)(0,0,0.5)};
				\addplot3[color=asb][line width=0.6pt,dotted]coordinates {(0,0.5,0)(0,0,0.5)};
				\fill[opacity=0.1, asb] (axis cs: 0,0,0) -- (axis cs: 0.5,0,0) -- (axis cs: 0,0.5,0) -- (axis cs: 0,0,0.5) -- cycle;
				
				\draw[color=asb] (-0.2, -0.2, -0.2) node[anchor=north east] {$_{v_{1}}$};
				\draw[color=asb] (-0.2, -0.2, 1.2) node[anchor=south east] {$_{v_{2}}$};
				\draw[color=asb] (-0.2, 1.2, -0.2) node[anchor=south west] {$_{v_{3}}$};
				\draw[color=asb] (1.2, -0.2, -0.2) node[anchor=north west] {$_{v_{4}}$};
				
				\draw[line width=.6pt, color=asb](-0.2, -0.2, 0)--(-0.2,-0.2,1);
				\draw[line width=.6pt, color=asb](0, -0.2, -0.2)--(1,-0.2,-0.2);
				\draw[line width=.6pt, color=asb](-0.2, 0, -0.2)--(-0.2,1,-0.2);
				\draw[line width=.6pt, color=asb](-0, 1.0, -0.2)--(1,0,-0.2);
				\draw[line width=.6pt, color=asb](0,-0.2,1)--(1,-0.2,0);
				\draw[line width=.6pt, color=asb](-0.2,0,1)--(-0.2,1,0);
				
				\draw[color=asb] (-0.2, -0.2, 0.5) node[anchor=east] {$_{e_{12}}$};
				\draw[color=asb] (0.5, -0.2, -0.2) node[anchor=north] {$_{e_{14}}$};
				\draw[color=asb] (-0.2,0.5,0.5) node[anchor=west] {$_{e_{23}}$};
				\draw[color=asb] (0.5,0.5,-0.2) node[anchor=west] {$_{e_{34}}$};
				\draw[color=asb] (0.5,-0.2,0.6) node[anchor=north west] {$_{e_{24}}$};
				\draw[color=asb] (0.1,0.1,-0.04) node[anchor=south] {$_{f_{124}}$};
				\draw[color=asb] (0.06,0.06,0.1) node[anchor=south] {$_{c_{1234}}$};
			\end{axis}
		\end{tikzpicture}
    \caption{Decomposition of the unit tetrahedron into vertices, edges, faces and the cell.}
    \label{fig:tet}
\end{figure}
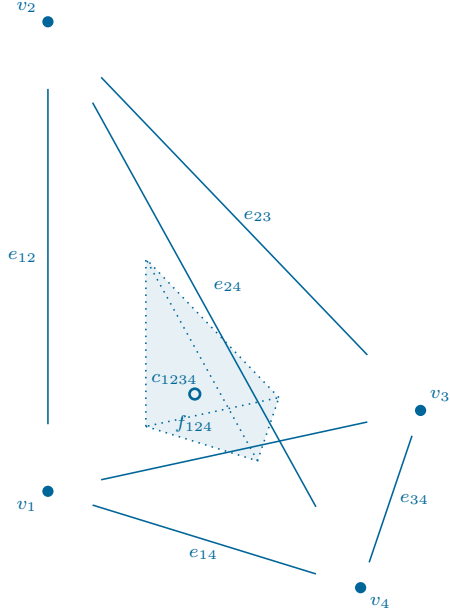

This section is dedicated to the introduction of polytopal templates on the reference tetrahedron 
\begin{align}
    \Omega = \{ (\xi, \eta, \zeta) \in [0,1]^3 \; | \; \xi + \eta + \zeta \leq 1 \} \, .
\end{align}
To that end, the tetrahedron is decomposed into its base polytopes given by its vertices $\{v_1, v_2, v_3, v_4\}$, its edges $\{e_{12}, e_{13}, e_{14}, e_{23}, e_{24}, e_{34}\}$, its faces $\{f_{123},f_{124},f_{134},f_{234}\}$, and its cell $c_{1234}$, see \cref{fig:tet}. Further, each polytope is associated with base functions belonging to an $\Hone$-conforming subspace $\U^p(\Omega)$ with $\dim \U^p(\Omega) = \dim \Po^p(\Omega) = (p+3)(p+2)(p+1)/6$.
\begin{definition}[Tetrahedron $\U^p(\Omega)$-polytopal spaces]
Each polytope is associated with a space of base functions as follows:
\begin{itemize}
    \item Each vertex is associated with the space of its respective base function $\ver^p_i$. As such, there are four spaces in total $i \in \{1,2,3,4\}$ and each one is of dimension one, $\dim \ver^p_i = 1 \quad \forall \, i \in \{1,2,3,4\}$.
    \item For each edge there exists a space of edge functions $\edge^p_{j}$ with  $j \in \mathcal{J} = \{(1,2),(1,3),(1,4),(2,3),(2,4),(3,4)\}$.
    The dimension of each edge space is given by $\dim \edge^p_j = p-1$.
    \item For each face there exists a space of of face base functions $\face_{k}^p$ with $k \in \mathcal{K} =  \{(1,2,3),(1,2,4),(1,3,4),(2,3,4)\}$, where the dimension of the spaces reads $\dim \face_{k} = (p-2)(p-1)/2$.
    \item Lastly, the space of cell base function is given by $\cell_{1234}$ with the dimensionality $\dim \cell_{1234} = (p-3)(p-2)(p-1)/6$.
\end{itemize}
The association with a respective polytope is according to \cref{def:poly} with respect to the trace operator for $\Hone$.
\end{definition}
We present polytopal templates for the construction of N\'ed\'elec elements of the second type and Brezzi-Douglas-Marini elements. For the formulation of N\'ed\'elec elements of the first type and Raviart-Thomas elements using Legendre polynomials see \cite{Zaglmayr2006,Joachim2005}, or \cite{Ainsworth2015,AINSWORTH2018178} for a Bernstein basis. 

\begin{figure}
		\centering
		\begin{subfigure}{0.24\linewidth}
			\centering
			\includegraphics[width=0.95\linewidth]{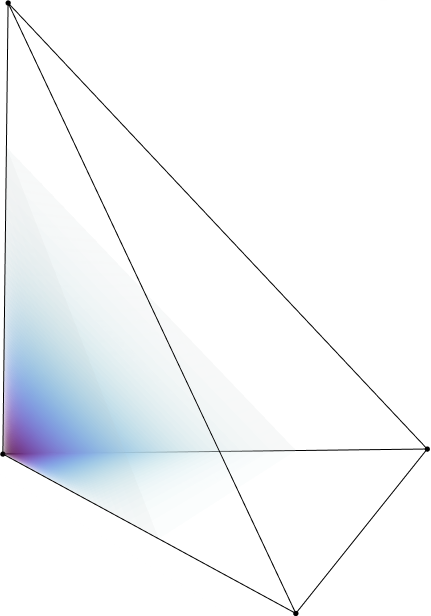}
			\caption{}
		\end{subfigure}
	    \begin{subfigure}{0.24\linewidth}
	    	\centering
	    	\includegraphics[width=0.95\linewidth]{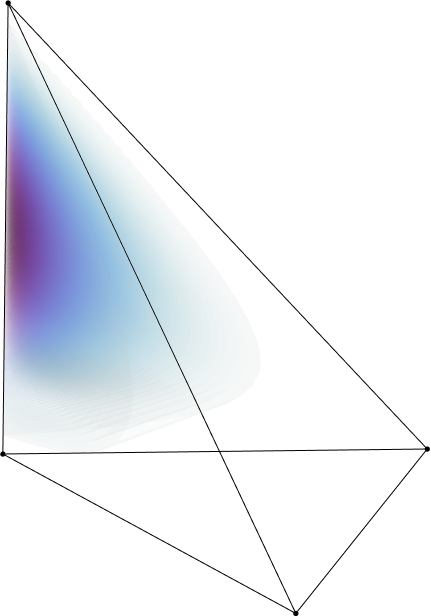}
	    	\caption{}
	    \end{subfigure}
        \begin{subfigure}{0.24\linewidth}
        	\centering
        	\includegraphics[width=0.95\linewidth]{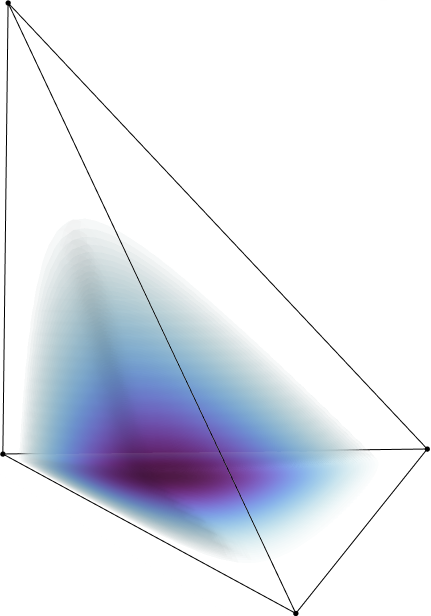}
        	\caption{}
        \end{subfigure}
        \begin{subfigure}{0.24\linewidth}
        	\centering
        	\includegraphics[width=0.95\linewidth]{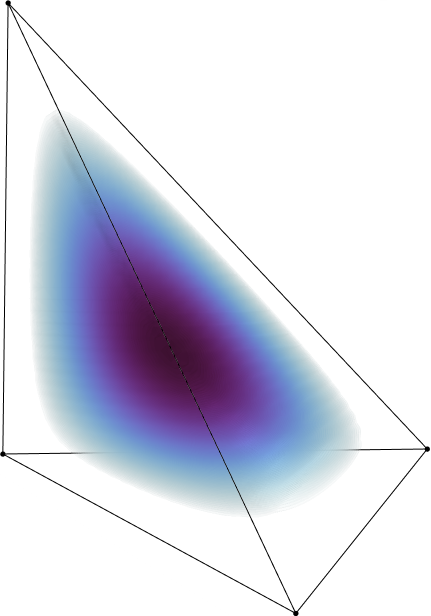}
        	\caption{}
        \end{subfigure}
		\caption{Vertex (a), edge (b), face (c), and cell (c) base functions  on the reference tetrahedron. Dark colours represent higher values of the field.}
		\label{fig:beziertet}
	\end{figure} 

\subsection{N\'ed\'elec II}
We proceed analogously to the definition of the N\'ed\'elec element of the second type \cite{Ned2} for triangles by constructing a polytopal template on the unit tetrahedron. The template is then used in conjunction with an $\Hone$-conforming polynomial basis to span the $\Nedtwo$-space on the unit tetrahedron. 
	
	We define the vertex-edge tangent vector $\vb{e}_3$ for $v_1$-$e_{12}$ such that its tangential projection is one on $e_{12}$ and zero on all other neighbouring edges. The same vector is the tangent vector associated with edge $e_{12}$. Next we define the edge-face vector $-\vb{e}_2$ for the edge $e_{12}$ and face $f_{123}$. On the face $f_{123}$ we employ the base vectors $\vb{e}_3$ and $\vb{e}_2$ which span the plane on the face. Lastly, we employ the full set of the Cartesian base vectors for the cell, namely $\vb{e}_3$, $\vb{e}_2$ and $\vb{e}_1$. The template vectors of the remaining polytopes are derived by covariant Piola transformations of the unit tetrahedron $c_{1234}$ to equivalent permutations $c_{ijkl}$ and by adjusting the sign of the vector to ensure a positive projection on the tangent vector, analogous to \cref{fig:permut}.  
\begin{figure}
		\centering
		\definecolor{asl}{rgb}{0.4980392156862745,0.,1.}
		\definecolor{asb}{rgb}{0.,0.4,0.6}
		\begin{tikzpicture}
			\begin{axis}
				[
				width=30cm,height=25cm,
				view={50}{15},
				enlargelimits=true,
				xmin=-1,xmax=2,
				ymin=-1,ymax=2,
				zmin=-1,zmax=2,
				domain=-10:10,
				axis equal,
				hide axis
				]
				\draw (-0.2, -0.2, -0.2) node[circle,fill=asb,inner sep=1.5pt] {};
				\draw (-0.2, -0.2, 1.2) node[circle,fill=asb,inner sep=1.5pt] {};
				\draw (1.2, -0.2, -0.2) node[circle,fill=asb,inner sep=1.5pt] {};
				\draw (-0.2, 1.2, -0.2) node[circle,fill=asb,inner sep=1.5pt] {};
				
				\addplot3[color=asl][line width=1pt,mark=o]
				coordinates {(0.1, 0.1, 0.1)};
				
				\addplot3[color=asb][line width=0.6pt,dotted]
				coordinates {(0,0,0)(0.5,0,0)(0,0.5,0)(0,0,0)};
				\addplot3[color=asb][line width=0.6pt,dotted]
				coordinates {(0,0,0)(0,0,0.5)};
				\addplot3[color=asb][line width=0.6pt,dotted]coordinates {(0.5,0,0)(0,0,0.5)};
				\addplot3[color=asb][line width=0.6pt,dotted]coordinates {(0,0.5,0)(0,0,0.5)};
				\fill[opacity=0.1, asb] (axis cs: 0,0,0) -- (axis cs: 0.5,0,0) -- (axis cs: 0,0.5,0) -- (axis cs: 0,0,0.5) -- cycle;
				
				\draw[color=asb] (-0.2, -0.2, -0.2) node[anchor=north east] {$_{v_{1}}$};
				\draw[color=asb] (-0.2, -0.2, 1.2) node[anchor=south east] {$_{v_{2}}$};
				\draw[color=asb] (-0.2, 1.2, -0.2) node[anchor=south west] {$_{v_{3}}$};
				\draw[color=asb] (1.2, -0.2, -0.2) node[anchor=north west] {$_{v_{4}}$};
				
				\draw[line width=.6pt, color=asb](-0.2, -0.2, 0)--(-0.2,-0.2,1);
				\draw[line width=.6pt, color=asb](0, -0.2, -0.2)--(1,-0.2,-0.2);
				\draw[line width=.6pt, color=asb](-0.2, 0, -0.2)--(-0.2,1,-0.2);
				\draw[line width=.6pt, color=asb](-0, 1.0, -0.2)--(1,0,-0.2);
				\draw[line width=.6pt, color=asb](0,-0.2,1)--(1,-0.2,0);
				\draw[line width=.6pt, color=asb](-0.2,0,1)--(-0.2,1,0);
				
				\draw[-to, line width=1.pt, color=asl](-0.2, -0.2, 0)--(-0.2,-0.2,0.2);
				\draw[-to, line width=1.pt, color=asl](0, -0.2, -0.2)--(0.2,-0.2,-0.2);
				\draw[-to, line width=1.pt, color=asl](-0.2, 0, -0.2)--(-0.2,0.2,-0.2);
				
				\draw[to-, line width=1.pt, color=asl](-0.2, -0.2, 1)--(-0.3,-0.3,0.9);
				\draw[to-, line width=1.pt, color=asl](1,-0.2,-0.2)--(0.9,-0.3,-0.3);
				\draw[to-, line width=1.pt, color=asl](-0.2,1,-0.2)--(-0.3,0.9,-0.3);
				
				\draw[-to, line width=1.pt, color=asl](0,-0.2,1)--(0.2,-0.2,1);
				\draw[-to, line width=1.pt, color=asl](-0.2,0,1)--(-0.2,0.2,1);
				\draw[-to, line width=1.pt, color=asl](1,-0.2,0.2)--(1,-0.2,0);
				\draw[-to, line width=1.pt, color=asl](-0.2,1,0.2)--(-0.2,1,0);
				
				\draw[-to, line width=1.pt, color=asl](1,0.2,-0.2)--(1,0,-0.2);
				\draw[-to, line width=1.pt, color=asl](0, 1.0, -0.2)--(0.2, 1.0, -0.2);
				
				\draw[-to, line width=1.pt, color=asl,densely dashed](-0.2,-0.2,0.4)--(-0.2,-0.2,0.6);
				\draw[-to, line width=1.pt, color=asl,densely dashed](0.4,-0.2,-0.2)--(0.6,-0.2,-0.2);
				\draw[-to, line width=1.pt, color=asl, densely dashed](0.5,-0.2,0.5)--(0.7,-0.2,0.5);
				\draw[-to, line width=1.pt, color=asl, densely dashed](-0.2,0.5,0.5)--(-0.2,0.7,0.5);
				\draw[-to, line width=1.pt, color=asl, densely dashed](0.5,0.5,-0.2)--(0.7,0.5,-0.2);
				
				\draw[-to, line width=1.pt, color=asl, dashdotted](-0.2,-0.2,0.5)--(-0.4,-0.2,0.5);
				\draw[-to, line width=1.pt, color=asl, dashdotted](-0.2,-0.2,0.5)--(-0.2,-0.4,0.5);
				\draw[to-, line width=1.pt, color=asl, dashdotted](0.5,-0.2,-0.2)--(0.5,-0.4,-0.2);
				\draw[to-, line width=1.pt, color=asl, dashdotted](0.5,-0.2,-0.2)--(0.5,-0.2,-0.4);
				\draw[to-, line width=1.pt, color=asl, dashdotted](-0.2,0.5,-0.2)--(-0.2,0.5,-0.4);
				
				\draw[-to, line width=1.pt, color=asl, dashdotted](0.5,-0.2,0.5)--(0.5,0,0.5);
				\draw[-to, line width=1.pt, color=asl, dashdotted](0.5,-0.2,0.5)--(0.6,-0.1,0.6);
				
				\draw[-to, line width=1.pt, color=asl, dashdotted](-0.2,0.5,0.5)--(-0.1,0.6,0.6);
				\draw[to-, line width=1.pt, color=asl, dashdotted](-0.2,0.5,0.5)--(0,0.5,0.5);
				
				\draw[-to, line width=1.pt, color=asl, dashdotted](0.5,0.5,-0.2)--(0.6,0.6,-0.1);
				\draw[to-, line width=1.pt, color=asl, dashdotted](0.5,0.5,-0.2)--(0.5,0.5,0);
				
				\draw[-to, line width=1.pt, color=asl, dotted](0.1,0,0.1)--(0.1,0,0.3);
				\draw[-to, line width=1.pt, color=asl, dotted](0.1,0,0.1)--(0.3,0,0.1);
				
				\draw[to-, line width=1.pt, color=asl, dashdotdotted](0.1,0,0.1)--(0.1,-0.2,0.1);
				
				\draw[-to, line width=1.pt, color=asl, dotted](0.2,0.2,0.1)--(0.2,0.4,0.1);
				\draw[-to, line width=1.pt, color=asl, dotted](0.2,0.2,0.1)--(0.4,0.2,0.1);
				
				\draw[-to, line width=1.pt, color=asl, dashdotdotted](0.2,0.2,0.1)--(0.3,0.3,0.2);
				
				\draw[-to,line width=1pt, color=asl](1, 1, 1)--(1.2, 1.2, 1.007);
				\draw[color=asl] (1.2, 1.2, 1.007) node[anchor=west] {Vertex-edge template vectors};
				\draw[-to,line width=1pt, color=asl, densely dashed](1, 1, 1-0.15)--(1.2, 1.2, 1.007-0.15);
				\draw[color=asl] (1.2, 1.2, 1.007-0.15) node[anchor=west] {Edge template vectors};
				\draw[-to,line width=1pt, color=asl, dashdotted](1, 1, 1-0.3)--(1.2, 1.2, 1.007-0.3);
				\draw[color=asl] (1.2, 1.2, 1.007-0.3) node[anchor=west] {Edge-face template vectors};
				\draw[-to,line width=1pt, color=asl, dotted](1, 1, 1-0.45)--(1.2, 1.2, 1.007-0.45);
				\draw[color=asl] (1.2, 1.2, 1.007-0.45) node[anchor=west] {Face template vectors};
				\draw[-to,line width=1pt, color=asl,dashdotdotted](1, 1, 1-0.6)--(1.2, 1.2, 1.007-0.6);
				\draw[color=asl] (1.2, 1.2, 1.007-0.6) node[anchor=west] {Face-cell template vectors};
				
				\addplot3[color=asl][line width=1pt,mark=o]
				coordinates {(1.1, 1.1, 1-0.75)};
				\draw[color=asl] (1.2, 1.2, 1.007-0.75) node[anchor=west] {Cell-Cartesian template vectors};
				
				\draw[color=asb] (-0.2, -0.2, 0.5) node[anchor=west] {$_{e_{12}}$};
				\draw[color=asb] (0.5, -0.2, -0.2) node[anchor=south] {$_{e_{14}}$};
				\draw[color=asb] (-0.2,0.5,0.5) node[anchor=east] {$_{e_{23}}$};
				\draw[color=asb] (0.5,0.5,-0.2) node[anchor=east] {$_{e_{34}}$};
				\draw[color=asb] (0.5,-0.2,0.6) node[anchor=south] {$_{e_{24}}$};
				\draw[color=asb] (0.1,0.1,-0.05) node[anchor=south] {$_{f_{124}}$};
			\end{axis}
		\end{tikzpicture}
		\caption{Template vectors for the reference tetrahedron on their corresponding polytopes. Only vectors on the visible sides of the tetrahedron are depicted. The template allows to construct N\'ed\'elec elements of the second type.}
		\label{fig:tet_nii}
	\end{figure}
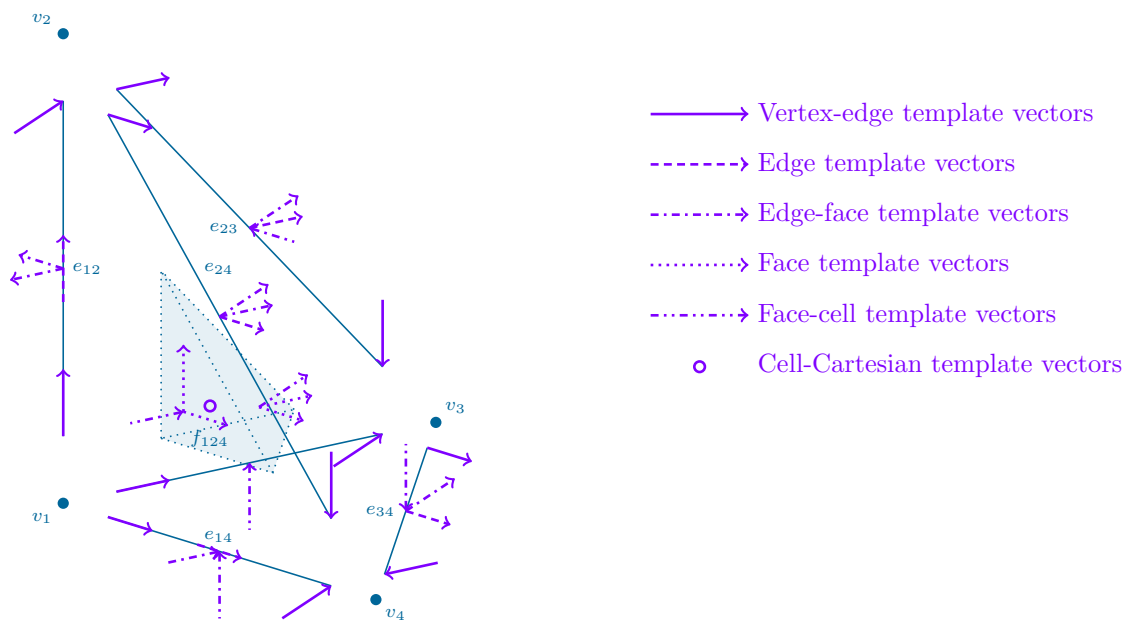
The resulting template is given by the super-set of the following polytopal sets depicted in \cref{fig:tet_nii}
    \begin{align}
    	\tem_1 &= \{ \vb{e}_3,\vb{e}_2,\vb{e}_1 \} \, , & \tem_2 &= \{ \vb{e}_1 + \vb{e}_2 + \vb{e}_3 , \vb{e}_2 , \vb{e}_1 \} \, , & \tem_3 &= \{ \vb{e}_1 + \vb{e}_2 + \vb{e}_3, -\vb{e}_3 ,\vb{e}_1 \} \, , \notag \\
    	\tem_4 &= \{ \vb{e}_1 + \vb{e}_2 + \vb{e}_3, -\vb{e}_3, -\vb{e}_2\} \, , & \tem_{12} &= \{ \vb{e}_3, -\vb{e}_2, -\vb{e}_1 \} \, , & \tem_{13} &= \{ \vb{e}_2, \vb{e}_3, -\vb{e}_1 \} \, , \notag \\
    	\tem_{14} &= \{ \vb{e}_1, \vb{e}_3, \vb{e}_2 \} \, , & \tem_{23} &= \{ \vb{e}_2, \vb{e}_1 + \vb{e}_2 + \vb{e}_3, -\vb{e}_1 \} \, , & \tem_{24} &= \{ \vb{e}_1, \vb{e}_1 + \vb{e}_2 + \vb{e}_3, \vb{e}_2 \} \, , \notag \\
    	\tem_{34} &= \{ \vb{e}_1, \vb{e}_1 + \vb{e}_2 + \vb{e}_3, -\vb{e}_3 \} \, , & \tem_{123} &= \{ \vb{e}_3, \vb{e}_2, -\vb{e}_1 \} \, , & \tem_{124} &= \{ \vb{e}_3, \vb{e}_1, \vb{e}_2 \} \, , \notag \\
    	\tem_{134} &= \{\vb{e}_2, \vb{e}_1, -\vb{e}_3\} \, , & \tem_{234}  &= \{ \vb{e}_2, \vb{e}_1,\vb{e}_1 + \vb{e}_2 + \vb{e}_3\} \, , & \tem_{1234} &= \{ \vb{e}_3, \vb{e}_2, \vb{e}_1 \} \, ,
    \end{align}
	and reads
	\begin{align}
		\tem = \{ \tem_1,\tem_2,\tem_3,\tem_4,\tem_{12},\tem_{13},\tem_{14},\tem_{23},\tem_{24},\tem_{34},\tem_{123},\tem_{124},\tem_{134},\tem_{234},\tem_{1234} \} \, .
	\end{align}
	As such, we can define the N\'ed\'elec element of the second type with polytopal template an underlying subspace $\U^p(\Omega)$
	\begin{align}
    		&\Nedtwo^p = \left\{ \bigoplus_{i=1}^4 \ver_i^p \otimes \tem_i \right\} \oplus \left\{ \bigoplus_{j \in \mathcal{J}  } \edge^p_j \otimes \tem_j \right\} \oplus \left\{ \bigoplus_{k \in \mathcal{K}} \face^p_k \otimes \tem_k  \right\} \oplus \{ \cell^p_{1234} \otimes \tem_{1234} \} \, , \notag \\
    		&\mathcal{J} = \{ (1,2),(1,3),(1,4),(2,3),(2,4),(3,4) \} \, , \qquad \mathcal{K} = \{ (1,2,3),(1,2,4),(1,3,4),(2,3,4) \} \, ,
    	\end{align}
    where $\ver^p_i$ are the sets of vertex base functions, $\edge^p_j$ are the sets of edge base functions, $\face^p_k$ are the sets of face base functions and $\cell^p_{1234}$ is the set of cell base functions.
    \begin{definition}[Tetrahedron $\Nedtwo^p$ base functions]
	The base functions of the tetrahedral N\'ed\'elec element of the second type are defined on their respective polytopes as follows.
	\begin{itemize}
	    \item On each edge $e_{ij}$ with vertices $v_i$ and $v_j$, the base functions read
	    \begin{subequations}
	        \begin{align}
	        \text{Vertex-edge:}& &\bm{\vartheta}(\xi, \eta, \zeta) &= n \tv \, , & n &\in \ver^p_i \, , &\tv &\in \left \{ \tv \in \tem_i \; | \; \rtr \tv \at_{e_{ij}}  \neq 0 \right \} \, , \\
	        &&\bm{\vartheta}(\xi, \eta, \zeta) &= n \tv \, , & n &\in \ver^p_j \, , &\tv &\in \left \{ \tv \in \tem_j \; | \; \rtr \tv \at_{e_{ij}}  \neq 0 \right \} \, , \\
	        \text{Edge:} &&\bm{\vartheta}(\xi, \eta, \zeta) &= n \tv \, , & n &\in \edge^p_{ij} \, , &\tv &\in \left \{ \tv \in \tem_{ij} \; | \; \rtr \tv \at_{e_{ij}} \neq 0  \right \} \, ,
	    \end{align}
	    \end{subequations}
	    such that $\tv$ is an element of the template sets, whose tangential trace does not vanish on the edge. 
	    \item For each face $f_{ijk}$ with edges $e_{ij}$, $e_{ik}$ and $e_{jk}$ the base functions read
	    \begin{subequations}
	        \begin{align}
	        \text{Edge-face:}& &\bm{\vartheta}(\xi, \eta, \zeta) &= n \tv \, , & n &\in \edge^p_{ij} \, , &\tv &\in \left \{ \tv \in \tem_{ij} \; | \; \rtr \tv \at_{e_{ij}}  = 0 \, , \quad \ttr \tv \at_{f_{ijk}}  \neq 0 \right \} \, , \\
	        & &\bm{\vartheta}(\xi, \eta, \zeta) &= n \tv \, , & n &\in \edge^p_{ik} \, , &\tv &\in \left \{ \tv \in \tem_{ik} \; | \; \rtr \tv \at_{e_{ik}} = 0 \, , \quad \ttr \tv \at_{f_{ijk}}  \neq 0 \right \} \, , \\
	        & &\bm{\vartheta}(\xi, \eta, \zeta) &= n \tv \, , & n &\in \edge^p_{jk} \, , &\tv &\in \left \{ \tv \in \tem_{jk} \; | \; \rtr \tv \at_{e_{jk}} = 0 \, , \quad \ttr \tv \at_{f_{ijk}}  \neq 0 \right \} \, , \\
	        \text{Face:}& &\bm{\vartheta}(\xi, \eta, \zeta) &= n \tv \, , & n &\in \cell^p_{ijk} \, , &\tv &\in \left \{ \tv \in \tem_{ijk} \; | \; \ttr \tv \at_{f_{ijk}} \neq 0  \right \} \, .
	    \end{align}
	    \end{subequations}
	    such that their tangential trace vanishes on all edges.
	    \item The cell base functions read
	    \begin{subequations}
	        \begin{align}
	        \text{Face-cell:}& &\bm{\vartheta}(\xi, \eta, \zeta) &= n \tv \, , & n &\in \cell^p_{123} \, , &\tv &\in \left \{ \tv \in \tem_{123} \; | \; \ttr \tv \at_{f_{123}} = 0  \right \} \, , \\
	        & &\bm{\vartheta}(\xi, \eta, \zeta) &= n \tv \, , & n &\in \cell^p_{124} \, , &\tv &\in \left \{ \tv \in \tem_{124} \; | \; \ttr \tv \at_{f_{124}} = 0  \right \} \, , \\
	        & &\bm{\vartheta}(\xi, \eta, \zeta) &= n \tv \, , & n &\in \cell^p_{134} \, , &\tv &\in \left \{ \tv \in \tem_{134} \; | \; \ttr \tv \at_{f_{134}} = 0  \right \} \, , \\
	        & &\bm{\vartheta}(\xi, \eta, \zeta) &= n \tv \, , & n &\in \cell^p_{234} \, , &\tv &\in \left \{ \tv \in \tem_{234} \; | \; \ttr \tv \at_{f_{234}} = 0  \right \} \, , \\
	        \text{Cell:}& &\bm{\vartheta}(\xi, \eta, \zeta) &= n \tv \, , & n &\in \cell^p_{1234} \, , &\tv & \in \tem_{1234}  \, , 
	    \end{align}
	    \end{subequations}
	\end{itemize}
	\end{definition}
	 \begin{theorem} [Linear independence] 
	 \label{th:linin-tet}
    	Let $\U^p$ be an $\Hone$-conforming basis on the unit tetrahedron, then its tensor product with the polytopal template yields a unisolvent N\'ed\'elec element of the second type.
    \end{theorem}
    \begin{proof}
    	Linear independence of the vectorial base functions is inherited from the linear independence of the underlying $\Hone$-conforming basis since each base function is multiplied with three linearly independent template vectors. Further, the dimension of the basis 
    	\begin{align}
    		\dim [\U^p(\Omega)]^3 = \dim[\Po^p(\Omega)]^3 = \dim \Nedtwo^p(\Omega) \, , 
    	\end{align}
        agrees with the dimension of the N\'ed\'elec element of the second type.
    \end{proof}
    \begin{theorem} [$\Hc{,\body}$-conformity]
    	\label{th:tet-conf}
    	Under covariant Piola transformations, the basis spans an $\Hc{}$-conforming subspace. 
    \end{theorem}
    \begin{proof}
    	The conformity of a grid composed solely of unit tetrahedrons is $\Hc{}$-conforming due to the methodology used to construct the polytopal template, namely by an initial definition of a minimal set of template vectors and permutations of the unit tetrahedron, such that the jump condition on the interface of neighbouring elements reduces from $\jump{\ttr \ud}|_{\Xi_{ij}} = 0$ to $\jump{\tr\langle \vb{t}, \, \ud \rangle}|_{\Xi_{ij}} = 0$ of the tangential components of the non-cell base functions. The latter is satisfied by the underlying $\U^p$-space. Conformity of a general grid is consequently achieved by employing covariant Piola transformations, as the tangential projection is maintained.
    \end{proof}
    \begin{figure}
		\centering
		\begin{subfigure}{0.3\linewidth}
			\centering
			\includegraphics[width=0.7\linewidth]{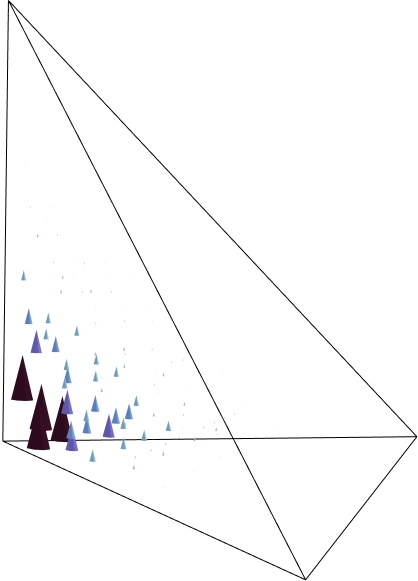}
			\caption{}
		\end{subfigure}
	    \begin{subfigure}{0.3\linewidth}
	    	\centering
	    	\includegraphics[width=0.7\linewidth]{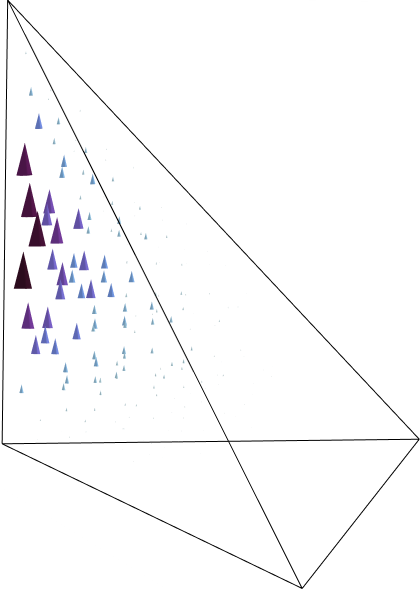}
	    	\caption{}
	    \end{subfigure}
        \begin{subfigure}{0.3\linewidth}
        	\centering
        	\includegraphics[width=0.7\linewidth]{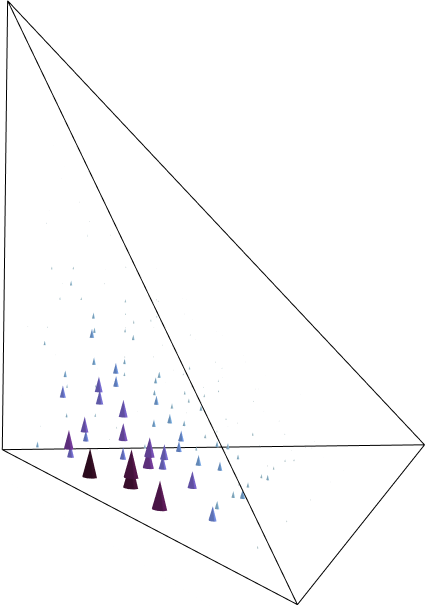}
        	\caption{}
        \end{subfigure}
        \begin{subfigure}{0.3\linewidth}
        	\centering
        	\includegraphics[width=0.7\linewidth]{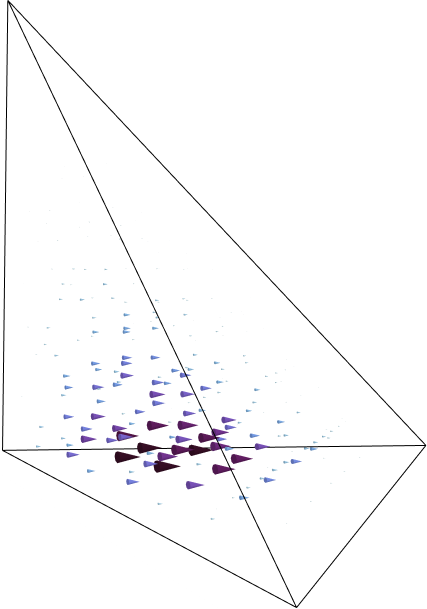}
        	\caption{}
        \end{subfigure}
        \begin{subfigure}{0.3\linewidth}
        	\centering
        	\includegraphics[width=0.7\linewidth]{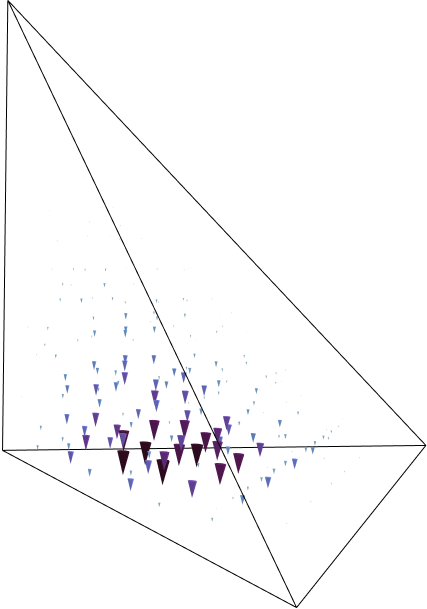}
        	\caption{}
        \end{subfigure}
        \begin{subfigure}{0.3\linewidth}
        	\centering
        	\includegraphics[width=0.7\linewidth]{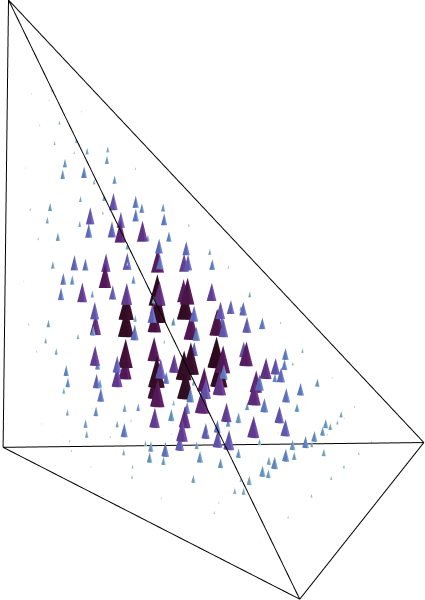}
        	\caption{}
        \end{subfigure}
		\caption{Quartic vertex-edge (a) and edge (b) base functions on $e_{12}$, edge-face (c) base function on $e_{14}$ for $f_{124}$, face (d) and face-cell (e) base functions on on $f_{134}$, and pure cell (f) base functions of the N\'ed\'elec element of the second type on the reference tetrahedron.}
		\label{fig:nedpiitet}
	\end{figure} 

\subsection{Brezzi-Douglas-Marini}
\begin{figure}
		\centering
		\definecolor{asl}{rgb}{0.4980392156862745,0.,1.}
		\definecolor{asb}{rgb}{0.,0.4,0.6}
		\begin{tikzpicture}
			\begin{axis}
				[
				width=30cm,height=25cm,
				view={50}{15},
				enlargelimits=true,
				xmin=-1,xmax=2,
				ymin=-1,ymax=2,
				zmin=-1,zmax=2,
				domain=-10:10,
				axis equal,
				hide axis
				]
				\draw (-0.2, -0.2, -0.2) node[circle,fill=asb,inner sep=1.5pt] {};
				\draw (-0.2, -0.2, 1.2) node[circle,fill=asb,inner sep=1.5pt] {};
				\draw (1.2, -0.2, -0.2) node[circle,fill=asb,inner sep=1.5pt] {};
				\draw (-0.2, 1.2, -0.2) node[circle,fill=asb,inner sep=1.5pt] {};
				
				\addplot3[color=asl][line width=1pt,mark=o]
				coordinates {(0.1, 0.1, 0.1)};
				
				\addplot3[color=asb][line width=0.6pt,dotted]
				coordinates {(0,0,0)(0.5,0,0)(0,0.5,0)(0,0,0)};
				\addplot3[color=asb][line width=0.6pt,dotted]
				coordinates {(0,0,0)(0,0,0.5)};
				\addplot3[color=asb][line width=0.6pt,dotted]coordinates {(0.5,0,0)(0,0,0.5)};
				\addplot3[color=asb][line width=0.6pt,dotted]coordinates {(0,0.5,0)(0,0,0.5)};
				\fill[opacity=0.1, asb] (axis cs: 0,0,0) -- (axis cs: 0.5,0,0) -- (axis cs: 0,0.5,0) -- (axis cs: 0,0,0.5) -- cycle;
				
				\draw[color=asb] (-0.2, -0.2, -0.2) node[anchor=north east] {$_{v_{1}}$};
				\draw[color=asb] (-0.2, -0.2, 1.2) node[anchor=east] {$_{v_{2}}$};
				\draw[color=asb] (-0.2, 1.2, -0.2) node[anchor=west] {$_{v_{3}}$};
				\draw[color=asb] (1.2, -0.2, -0.2) node[anchor=west] {$_{v_{4}}$};
				
				\draw[line width=.6pt, color=asb](-0.2, -0.2, 0)--(-0.2,-0.2,1);
				\draw[line width=.6pt, color=asb](0, -0.2, -0.2)--(1,-0.2,-0.2);
				\draw[line width=.6pt, color=asb](-0.2, 0, -0.2)--(-0.2,1,-0.2);
				\draw[line width=.6pt, color=asb](-0, 1.0, -0.2)--(1,0,-0.2);
				\draw[line width=.6pt, color=asb](0,-0.2,1)--(1,-0.2,0);
				\draw[line width=.6pt, color=asb](-0.2,0,1)--(-0.2,1,0);
				
				\draw[-to, line width=1.pt, color=asl](-0.2, -0.2, -0.2)--(-0.4,-0.2,-0.2);
				\draw[-to, line width=1.pt, color=asl](-0.2, -0.2, -0.2)--(-0.2,0,-0.2);
				\draw[-to, line width=1.pt, color=asl](-0.2, -0.2, -0.2)--(-0.2,-0.2,-0.4);
				
				\draw[-to, line width=1.pt, color=asl](-0.2, -0.2, 1.2)--(-0.4,-0.2,1.4);
				\draw[-to, line width=1.pt, color=asl](-0.2,-0.2,1.2)--(-0.2,0,1);
				\draw[-to, line width=1.pt, color=asl](-0.2,-0.2,1.2)--(-0.2,-0.2,1);
				
				\draw[-to, line width=1.pt, color=asl](-0.2, 1.2, -0.2)--(-0.4, 1.4, -0.2);
				\draw[-to, line width=1.pt, color=asl](-0.2, 1.2, -0.2)--(-0.2, 1.4, -0.4);
				\draw[-to, line width=1.pt, color=asl](-0.2, 1.2, -0.2)--(-0.2, 1, -0.2);
				
				\draw[-to, line width=1.pt, color=asl](1.2, -0.2, -0.2)--(1., -0., -0.2);
				\draw[-to, line width=1.pt, color=asl](1.2, -0.2, -0.2)--(1.4, -0.2, -0.4);
				\draw[-to, line width=1.pt, color=asl](1.2, -0.2, -0.2)--(1., -0.2, -0.2);
				
				\draw[-to, line width=1.pt, color=asl,dashdotted](-0.2,-0.2,0.4)--(-0.2,-0.2,0.6);
				\draw[-to, line width=1.pt, color=asl,dashdotted](0.4,-0.2,-0.2)--(0.6,-0.2,-0.2);
				\draw[-to, line width=1.pt, color=asl, dashdotted](0.4,-0.2,0.6)--(0.6,-0.2,0.4);
				\draw[-to, line width=1.pt, color=asl, dashdotted](-0.2,0.4,0.6)--(-0.2,0.6,0.4);
				\draw[-to, line width=1.pt, color=asl, dashdotted](0.4,0.6,-0.2)--(0.6,0.4,-0.2);
				
				\draw[-to, line width=1.pt, color=asl, densely dashed](-0.2,-0.2,0.5)--(-0.4,-0.2,0.5);
				\draw[to-, line width=1.pt, color=asl, densely dashed](-0.2,-0.2,0.5)--(-0.2,-0.4,0.5);
				\draw[to-, line width=1.pt, color=asl, densely dashed](0.5,-0.2,-0.2)--(0.5,-0.4,-0.2);
				\draw[-to, line width=1.pt, color=asl, densely dashed](0.5,-0.2,-0.2)--(0.5,-0.2,-0.4);
				\draw[-to, line width=1.pt, color=asl, densely dashed](-0.2,0.5,-0.2)--(-0.2,0.5,-0.4);
				
				\draw[to-, line width=1.pt, color=asl, densely dashed](0.5,-0.2,0.5)--(0.5,-0.4,0.7);
				\draw[to-, line width=1.pt, color=asl, densely dashed](0.5,-0.2,0.5)--(0.5,-0.2,0.7);
				
				\draw[to-, line width=1.pt, color=asl, densely dashed](-0.2,0.5,0.5)--(-0.2,0.5,0.7);
				\draw[-to, line width=1.pt, color=asl, densely dashed](-0.2,0.5,0.5)--(-0.4,0.5,0.7);
				
				\draw[-to, line width=1.pt, color=asl, densely dashed](0.5,0.5,-0.2)--(0.5,0.7,-0.4);
				\draw[-to, line width=1.pt, color=asl, densely dashed](0.5,0.7,-0.2)--(0.5,0.5,-0.2);
				
				\draw[-to, line width=1.pt, color=asl, dashdotdotted](0.1,0,0.1)--(0.1,0,0.3);
				\draw[-to, line width=1.pt, color=asl, dashdotdotted](0.1,0,0.1)--(0.3,0,0.1);
				
				\draw[to-, line width=1.pt, color=asl, dotted](0.1,0,0.1)--(0.1,-0.2,0.1);

				\draw[-to,line width=1pt, color=asl](1, 1, 1)--(1.2, 1.2, 1.007);
				\draw[color=asl] (1.2, 1.2, 1.007) node[anchor=west] {Vertex-face template vectors};
				\draw[-to,line width=1pt, color=asl, densely dashed](1, 1, 1-0.15)--(1.2, 1.2, 1.007-0.15);
				\draw[color=asl] (1.2, 1.2, 1.007-0.15) node[anchor=west] {Edge-face template vectors};
				\draw[-to,line width=1pt, color=asl, dashdotted](1, 1, 1-0.3)--(1.2, 1.2, 1.007-0.3);
				\draw[color=asl] (1.2, 1.2, 1.007-0.3) node[anchor=west] {Edge-cell template vectors};
				\draw[-to,line width=1pt, color=asl, dotted](1, 1, 1-0.45)--(1.2, 1.2, 1.007-0.45);
				\draw[color=asl] (1.2, 1.2, 1.007-0.45) node[anchor=west] {Face template vectors};
				\draw[-to,line width=1pt, color=asl,dashdotdotted](1, 1, 1-0.6)--(1.2, 1.2, 1.007-0.6);
				\draw[color=asl] (1.2, 1.2, 1.007-0.6) node[anchor=west] {Face-cell template vectors};
				
				\addplot3[color=asl][line width=1pt,mark=o]
				coordinates {(1.1, 1.1, 1-0.75)};
				\draw[color=asl] (1.2, 1.2, 1.007-0.75) node[anchor=west] {Cell-Cartesian template vectors};
				
				\draw[color=asb] (-0.2, -0.2, 0.5) node[anchor=west] {$_{e_{12}}$};
				\draw[color=asb] (0.5, -0.2, -0.2) node[anchor=south] {$_{e_{14}}$};
				\draw[color=asb] (-0.2,0.5,0.5) node[anchor=east] {$_{e_{23}}$};
				\draw[color=asb] (0.5,0.5,-0.2) node[anchor=east] {$_{e_{34}}$};
				\draw[color=asb] (0.5,-0.2,0.6) node[anchor=south] {$_{e_{24}}$};
				\draw[color=asb] (0.1,0.1,-0.05) node[anchor=south] {$_{f_{124}}$};
			\end{axis}
		\end{tikzpicture}
		\caption{Template vectors for the reference tetrahedron on their corresponding polytopes. Only vectors on the visible sides of the tetrahedron are depicted. The template is associated with the Brezzi-Douglas-Marini element.}
		\label{fig:tet_bdm}
	\end{figure}
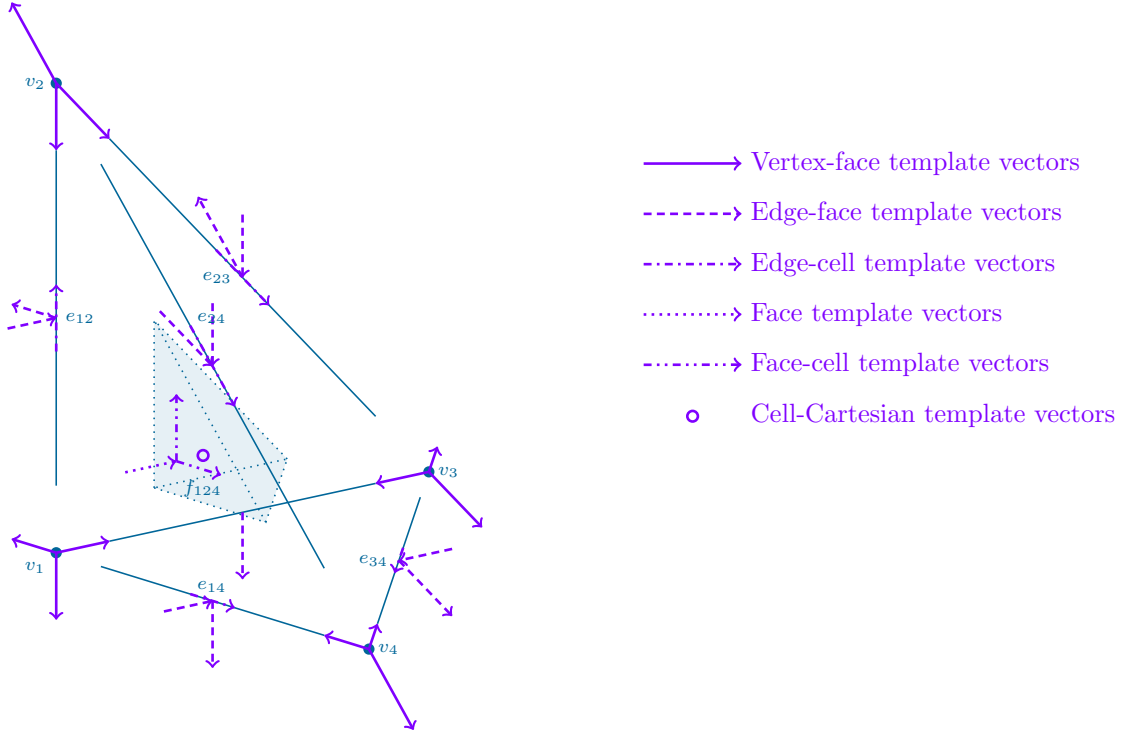

The last element we consider is the Brezzi-Douglas-Marini element \cite{BDM} on the unit tetrahedron. We start the construction by defining a single vertex-face vector for the vertex $v_1$ and the face $f_{123}$, namely $-\vb{e}_1$. Next we define two template vectors on the edge $e_{12}$. The vector $-\vb{e}_1$ is associated with the face $f_{123}$ and the vector $\vb{e}_3$ is an edge-cell vector. On the face introduce the template $\{-\vb{e}_1, \vb{e}_3,\vb{e}_2\}$, where the first vector is associated with the face and the last two are face-cell vectors. Lastly, we employ the Cartesian basis in the cell $c_{1234}$. The template vectors for the remaining polytopes are computed by permutations of the unit tetrahedron $c_{ijkl}$ with the contravariant Piola transformation, compare with \cref{fig:permut}.
The complete template is depicted in \cref{fig:tet_bdm},
where the polytopal template sets read
\begin{align}
    	\tem_1 &= \{ -\vb{e}_1,\vb{e}_2,-\vb{e}_3 \} \, , & \tem_2 &= \{ \vb{e}_3 - \vb{e}_1,  \vb{e}_2 - \vb{e}_3, -\vb{e}_3 \} \, , & \tem_3 &= \{ \vb{e}_2 - \vb{e}_1, \vb{e}_2 - \vb{e}_3, -\vb{e}_2 \} \, , \notag \\
    	\tem_4 &= \{ \vb{e}_2 - \vb{e}_1, \vb{e}_1 - \vb{e}_3 , -\vb{e}_1 \} \, , & \tem_{12} &= \{ -\vb{e}_1, \vb{e}_2, \vb{e}_3 \} \, , & \tem_{13} &= \{ -\vb{e}_1, -\vb{e}_3, \vb{e}_2 \} \, , \notag \\
    	\tem_{14} &= \{ \vb{e}_2, -\vb{e}_3 ,\vb{e}_1 \} \, , & \tem_{23} &= \{ \vb{e}_3 - \vb{e}_1 , -\vb{e}_3, \vb{e}_2 - \vb{e}_3  \} \, , & \tem_{24} &= \{ \vb{e}_2 - \vb{e}_3, -\vb{e}_3, \vb{e}_1 - \vb{e}_3 \} \, , \notag \\
    	\tem_{34} &= \{ \vb{e}_2 - \vb{e}_3, -\vb{e}_2, \vb{e}_1 - \vb{e}_2 \} \, , & \tem_{123} &= \{ -\vb{e}_1, \vb{e}_3, \vb{e}_2 \} \, , & \tem_{124} &= \{ \vb{e}_2, \vb{e}_3, \vb{e}_1 \} \, , \notag \\
    	\tem_{134} &= \{ -\vb{e}_3, \vb{e}_2, \vb{e}_1 \} \, , & \tem_{234}  &= \{ -\vb{e}_3, \vb{e}_2 - \vb{e}_3 , \vb{e}_1 - \vb{e}_3 \} \, , & \tem_{1234} &= \{\vb{e}_3, \vb{e}_2, \vb{e}_1 \} \, .
    \end{align}
	The super-set is given by
	\begin{align}
		\tem = \{ \tem_1,\tem_2,\tem_3,\tem_4,\tem_{12},\tem_{13},\tem_{14},\tem_{23},\tem_{24},\tem_{34},\tem_{123},\tem_{124},\tem_{134},\tem_{234},\tem_{1234} \} \, .
	\end{align}
	With template at hand, we define the $\BDM^p$-element using an underlying $\U^p$ subspace by tensor products
	\begin{align}
    		&\BDM^p = \left\{ \bigoplus_{i=1}^4 \ver_i^p \otimes \tem_i \right\} \oplus \left\{ \bigoplus_{j \in \mathcal{J}  } \edge^p_j \otimes \tem_j \right\} \oplus \left\{ \bigoplus_{k \in \mathcal{K}} \face^p_k \otimes \tem_k  \right\} \oplus \{ \cell^p_{1234} \otimes \tem_{1234} \} \, , \notag \\
    		&\mathcal{J} = \{ (1,2),(1,3),(1,4),(2,3),(2,4),(3,4) \} \, , \qquad \mathcal{K} = \{ (1,2,3),(1,2,4),(1,3,4),(2,3,4) \} \, ,
    	\end{align}
    where $\ver^p_i$ are the sets of vertex base functions, $\edge^p_j$ are the sets of edge base functions, $\face^p_k$ are the sets of face base functions and $\cell^p_{1234}$ is the set of cell base functions.
    \begin{definition}[Tetrahedron $\BDM^p$ base functions]
	The base functions of the tetrahedral Brezzi-Douglas-Marini element are defined on their respective polytope as follows.
	\begin{itemize}
	    \item For each face $f_{ijk}$ with vertices $v_i$, $v_j$ and $v_k$, and edges $e_{ij}$, $e_{ik}$ and $e_{jk}$ the base functions read
	    \begin{subequations}
	        \begin{align}
	        \text{Vertex-face:}& &\bm{\phi}(\xi, \eta, \zeta) &= n \tv \, , & n &\in \ver^p_i \, , &\tv &\in \left \{ \tv \in \tem_i \; | \; \ntr \tv \at_{f_{ijk}}  \neq 0 \right \} \, , \\
	        &&\bm{\phi}(\xi, \eta, \zeta) &= n \tv \, , & n &\in \ver^p_j \, , &\tv &\in \left \{ \tv \in \tem_j \; | \; \ntr \tv \at_{f_{ijk}}  \neq 0 \right \} \, , \\
	        &&\bm{\phi}(\xi, \eta, \zeta) &= n \tv \, , & n &\in \ver^p_k \, , &\tv &\in \left \{ \tv \in \tem_k \; | \; \ntr \tv \at_{f_{ijk}}  \neq 0 \right \} \, , \\
	        \text{Edge-face:} &&\bm{\phi}(\xi, \eta, \zeta) &= n \tv \, , & n &\in \edge^p_{ij} \, , &\tv &\in \left \{ \tv \in \tem_{ij} \; | \; \ntr \tv \at_{f_{ijk}} \neq 0  \right \} \, , \\
	        & &\bm{\phi}(\xi, \eta, \zeta) &= n \tv \, , & n &\in \edge^p_{ik} \, , &\tv &\in \left \{ \tv \in \tem_{ik} \; | \; \ntr \tv \at_{f_{ijk}} \neq 0  \right \} \, , \\
	        & &\bm{\phi}(\xi, \eta, \zeta) &= n \tv \, , & n &\in \edge^p_{jk} \, , &\tv &\in \left \{ \tv \in \tem_{jk} \; | \; \ntr \tv \at_{f_{ijk}} \neq 0  \right \} \, , \\
	        \text{Face:}& &\bm{\phi}(\xi, \eta, \zeta) &= n \tv \, , & n &\in \cell^p_{ijk} \, , &\tv &\in \left \{ \tv \in \tem_{ijk} \; | \; \ntr \tv \at_{f_{ijk}} \neq 0  \right \} \, .
	    \end{align}
	    \end{subequations}
	    such that their tangential trace vanishes on all edges.
	    \item The cell base functions read
	    \begin{subequations}
	        \begin{align}
	        \text{Edge-cell:} &&\bm{\phi}(\xi, \eta, \zeta) &= n \tv \, , & n &\in \edge^p_{12} \, , &\tv &\in \left \{ \tv \in \tem_{12} \; | \; \ntr \tv \at_{f_{123}} = 0 \, , \quad \ntr \tv \at_{f_{124}} = 0 \right \} \, , \\
	        &&\bm{\phi}(\xi, \eta, \zeta) &= n \tv \, , & n &\in \edge^p_{13} \, , &\tv &\in \left \{ \tv \in \tem_{13} \; | \; \ntr \tv \at_{f_{123}} = 0 \, , \quad \ntr \tv \at_{f_{134}} = 0 \right \} \, , \\
	        &&\bm{\phi}(\xi, \eta, \zeta) &= n \tv \, , & n &\in \edge^p_{14} \, , &\tv &\in \left \{ \tv \in \tem_{14} \; | \; \ntr \tv \at_{f_{124}} = 0 \, , \quad \ntr \tv \at_{f_{134}} = 0 \right \} \, , \\
	        &&\bm{\phi}(\xi, \eta, \zeta) &= n \tv \, , & n &\in \edge^p_{23} \, , &\tv &\in \left \{ \tv \in \tem_{23} \; | \; \ntr \tv \at_{f_{123}} = 0 \, , \quad \ntr \tv \at_{f_{234}} = 0 \right \} \, , \\
	        &&\bm{\phi}(\xi, \eta, \zeta) &= n \tv \, , & n &\in \edge^p_{24} \, , &\tv &\in \left \{ \tv \in \tem_{24} \; | \; \ntr \tv \at_{f_{124}} = 0 \, , \quad \ntr \tv \at_{f_{234}} = 0 \right \} \, , \\
	        &&\bm{\phi}(\xi, \eta, \zeta) &= n \tv \, , & n &\in \edge^p_{34} \, , &\tv &\in \left \{ \tv \in \tem_{34} \; | \; \ntr \tv \at_{f_{134}} = 0 \, , \quad \ntr \tv \at_{f_{234}} = 0 \right \} \, , \\
	        \text{Face-cell:}& &\bm{\phi}(\xi, \eta, \zeta) &= n \tv \, , & n &\in \cell^p_{123} \, , &\tv &\in \left \{ \tv \in \tem_{123} \; | \; \ntr \tv \at_{f_{123}} = 0  \right \} \, , \\
	        & &\bm{\phi}(\xi, \eta, \zeta) &= n \tv \, , & n &\in \cell^p_{124} \, , &\tv &\in \left \{ \tv \in \tem_{124} \; | \; \ntr \tv \at_{f_{124}} = 0  \right \} \, , \\
	        & &\bm{\phi}(\xi, \eta, \zeta) &= n \tv \, , & n &\in \cell^p_{134} \, , &\tv &\in \left \{ \tv \in \tem_{134} \; | \; \ntr \tv \at_{f_{134}} = 0  \right \} \, , \\
	        & &\bm{\phi}(\xi, \eta, \zeta) &= n \tv \, , & n &\in \cell^p_{234} \, , &\tv &\in \left \{ \tv \in \tem_{234} \; | \; \ntr \tv \at_{f_{234}} = 0  \right \} \, , \\
	        \text{Cell:}& &\bm{\phi}(\xi, \eta, \zeta) &= n \tv \, , & n &\in \cell^p_{1234} \, , &\tv & \in \tem_{1234}  \, , 
	    \end{align}
	    \end{subequations}
	\end{itemize}
	\end{definition}
A depiction of the base functions is given in \cref{fig:bdmtet}.
	\begin{theorem} [Linear independence] 
    	Let $\U^p$ be a polynomial $\Hone$-conforming subspace on the unit tetrahedron, then its tensor product with the polytopal template yields a unisolvent Brezzi-Douglas-Marini.
    \end{theorem}
    \begin{proof} 
    Each base function from $\U^p$ is multiplied with three linearly independent vectors, such that linear Independence is inherited on the vectorial level. 
    \end{proof}
    \begin{theorem} [$\Hd{,\body}$-conformity]
    	Under contravariant Piola transformations, the basis spans an $\Hd{}$-conforming subspace. 
    \end{theorem}
    \begin{proof}
    The proof follows the same lines as in \cref{th:tet-conf}. The base functions are designed with the property
    \begin{align}
        \langle \bm{\nu} , \, \bm{\phi}_i \rangle = n_i \, , 
    \end{align}
    for the normal components of the non-cell base functions, such that the trace interface condition is reduced to $\jump{\tr \langle \vb{n} , \, \ud \rangle}|_{\Xi_{ij}} = 0$ from $\jump{\ntr  \ud}|_{\Xi_{ij}} = 0$, which is satisfied by the underlying $\U^p$ subspace. 
    \end{proof}
\begin{figure}
		\centering
		\begin{subfigure}{0.3\linewidth}
			\centering
			\includegraphics[width=0.7\linewidth]{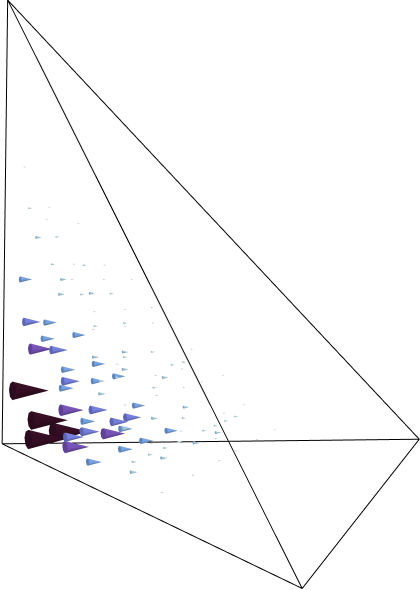}
			\caption{}
		\end{subfigure}
	    \begin{subfigure}{0.3\linewidth}
	    	\centering
	    	\includegraphics[width=0.7\linewidth]{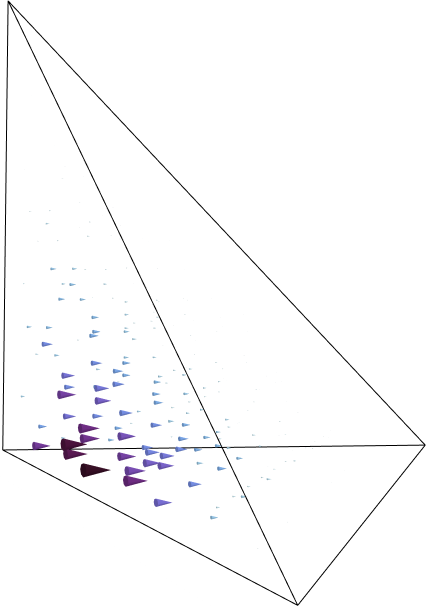}
	    	\caption{}
	    \end{subfigure}
        \begin{subfigure}{0.3\linewidth}
        	\centering
        	\includegraphics[width=0.7\linewidth]{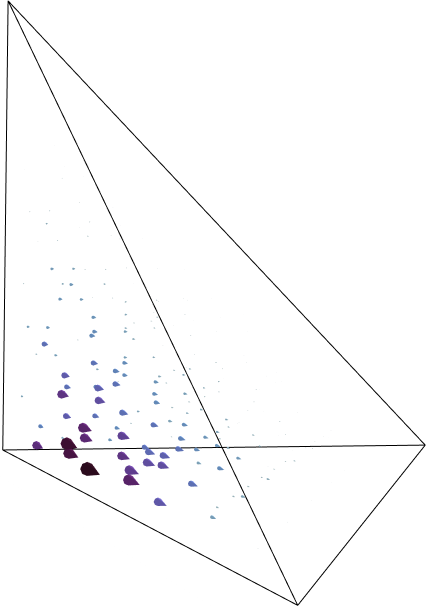}
        	\caption{}
        \end{subfigure}
        \begin{subfigure}{0.3\linewidth}
        	\centering
        	\includegraphics[width=0.7\linewidth]{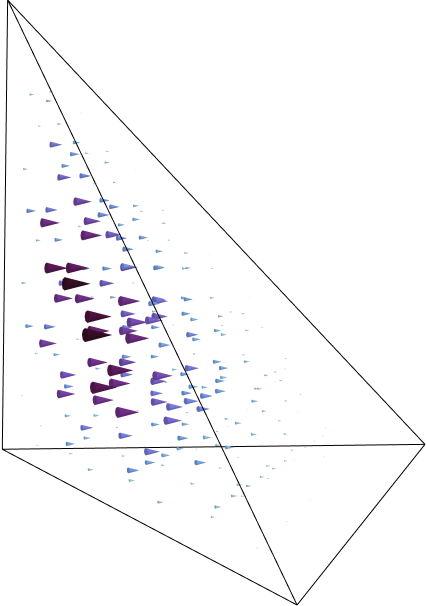}
        	\caption{}
        \end{subfigure}
        \begin{subfigure}{0.3\linewidth}
        	\centering
        	\includegraphics[width=0.7\linewidth]{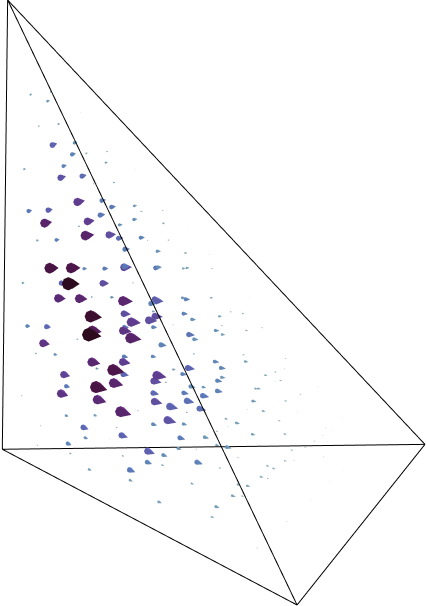}
        	\caption{}
        \end{subfigure}
        \begin{subfigure}{0.3\linewidth}
        	\centering
        	\includegraphics[width=0.7\linewidth]{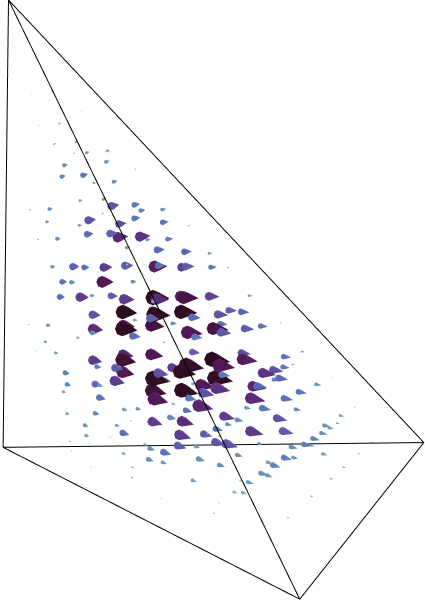}
        	\caption{}
        \end{subfigure}
		\caption{Quartic vertex-face (a), edge-face (b), edge-cell (c), face (d), face-cell (e) and pure cell (f) base functions of the Brezzi-Douglas-Marini element on the reference tetrahedron. The face functions are associated with the face $f_{124}$.}
		\label{fig:bdmtet}
	\end{figure}

\section{Examples}
In the following we present two examples using N\'ed\'elec elements and the relaxed micromorphic model. The decay in the error is measured in the Lebesgue norm 
\begin{align}
		&\| \widetilde{u} - u^h \|_{\Le} = \sqrt{\int_\body \| \widetilde{u} - u^h \|^2 \dd \body} \, ,
	\end{align}
in which context $\widetilde{u}$ is the analytical solution and $u^h$ is the finite element approximation.

\subsection{Relaxed micromorphic model of antiplane shear}
The strong from of the relaxed micromorphic model of antiplane shear \cite{Sky2021} reads
\begin{subequations}
		\begin{align}
			-\mue \di(\nabla u - \vb{p}) &= f  && \text{in} \quad \surf \, , \label{eq:strong_anti_p}  \\
			-\mue(\nabla u - \vb{p}) + \mumi \, \vb{p} + \muma \Lc^2 \rog \rot{\vb{p}} &= \vb{m}  && \text{in} \quad \surf \, , \label{eq:strong_anti_u} \\
			u &= \widetilde{u} && \text{on} \quad \curv_D^u \, , \\
			\langle \vb{p} , \, \vb{t} \rangle &= \langle \widetilde{\vb{p}} , \, \vb{t} \rangle && \text{on} \quad \curv_D^P \, , \\
			\langle \nabla u , \, \vb{n} \rangle &= \langle \vb{p} , \, \vb{n} \rangle && \text{on} \quad \curv_N^u \, ,\\
			\rot{\vb{p}} &= 0  && \text{on} \quad \curv_N^P \, ,
		\end{align}
	\end{subequations}
where $u$ is the antiplane displacement field, $\vb{p}$ is the microdistortion, $\mue, \mumi, \muma$ are shear material parameters and $\Lc$ is the characteristic length scale parameter. Body forces and micro-moments are given by $f$ and $\bm{M}$, respectively.
The corresponding bilinear and linear forms are given by 
\begin{subequations}
		\begin{align}
			a(\{\delta u, \delta \vb{p}\},\{u,\vb{p}\}) &= \int_\surf \mue \langle \nabla \delta u - \delta \vb{p} , \, \nabla u - \vb{p} \rangle + \mumi \langle \delta \vb{p} , \, \vb{p} \rangle + \muma \Lc^2 \rot{\delta \vb{p}} \rot{\vb{p}} \, \dd \surf \, , \\
			l(\{\delta u, \delta \vb{p}\}) &= \int_\surf \delta u \, f + \langle \delta \vb{p} , \, \vb{m} \rangle \, \dd \surf \, .
		\end{align}
		\label{eq:bili_anti}
	\end{subequations}
The problem is uniquely solvable for the space $\X(\surf) = \Hone(\surf) \times \Hc{,\surf}$ \cite{Sky2021}.

We set the material constants to $\mue = \mumi = \muma = \Lc = 1$ and prescribe an analytical solution by inserting predefined displacement and microdistortion fields 
\begin{align}
    	&\widetilde{u} =  \left \{ \begin{matrix}
    		(1 - y ^ 2)  (e^{x + 1} - 1) & \text{for} & x \leq 0  \\[2ex]
    		(1 - y ^ 2)  (e^{1-x} - 1) & \text{for} & x > 0 
    	\end{matrix} \right . \, , && \widetilde{\vb{p}} = \nabla \widetilde{u} = \left \{ \begin{matrix}
    		\begin{bmatrix}\left(1 - y^{2}\right) e^{x + 1}\\2 y \left(1 - e^{x + 1}\right)\end{bmatrix} & \text{for} & x \leq 0 \\[4ex]
    		\begin{bmatrix} \left(y^{2}- 1 \right) e^{1 - x}\\ 2 y \left(1-e^{1 - x} \right)\end{bmatrix} & \text{for} & x > 0 
    	\end{matrix} \right . \, .
    	\label{eq:arti_sol}
    \end{align}
into the strong form \cref{eq:strong_anti_u,eq:strong_anti_p}, such that the corresponding right-hand-side reads
\begin{align}
    	&f = 0 \, , && \vb{m} = \left \{ \begin{matrix}
    		\begin{bmatrix}\left(1 - y^{2}\right) e^{x + 1}\\2 y \left(1 - e^{x + 1}\right)\end{bmatrix} & \text{for} & x \leq 0 \\[4ex]
    		\begin{bmatrix}\left(y^{2} - 1\right) e^{1 - x}\\2 y \left(e^{x} - e\right) e^{- x}\end{bmatrix}& \text{for} & x > 0 
    	\end{matrix} \right . \, .
    \end{align}
The prescribed solution \cref{eq:arti_sol} is the analytical solution for the boundary conditions 
\begin{align}
    &u \at_{\partial \surf} = \widetilde{u} \, , && \langle \vb{p} , \, \vb{t} \rangle\at_{\partial \surf} = \langle \widetilde{\vb{p}} , \, \vb{t} \rangle \, ,
\end{align}
such that $\curv_D = \partial \surf$, since the solution is unique.
Clearly, the displacement field $u$ is $\C^0$-continuous and its gradient is solely tangentially continuous, such that N\'ed\'elec elements are required for optimal convergence results. 

\begin{figure}
    	\centering
    	\begin{subfigure}{0.3\linewidth}
    		\centering
    		\includegraphics[width=0.9\linewidth]{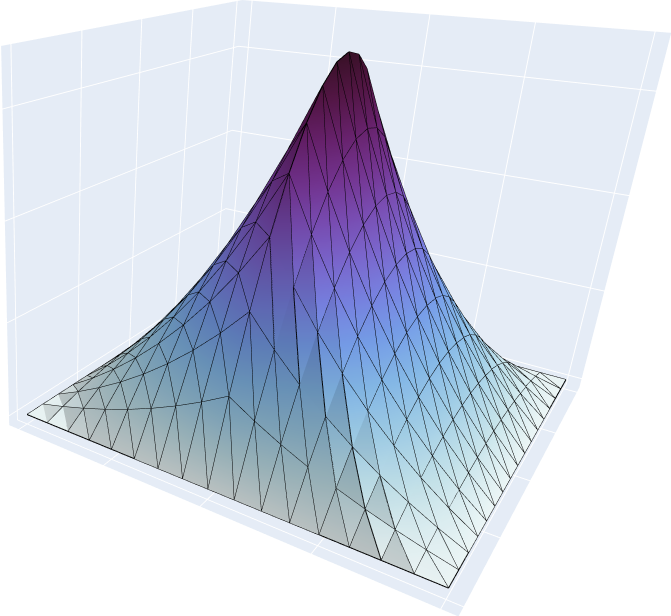}
    		\caption{}
    	\end{subfigure}
    	\begin{subfigure}{0.3\linewidth}
    		\centering
    		\includegraphics[width=0.9\linewidth]{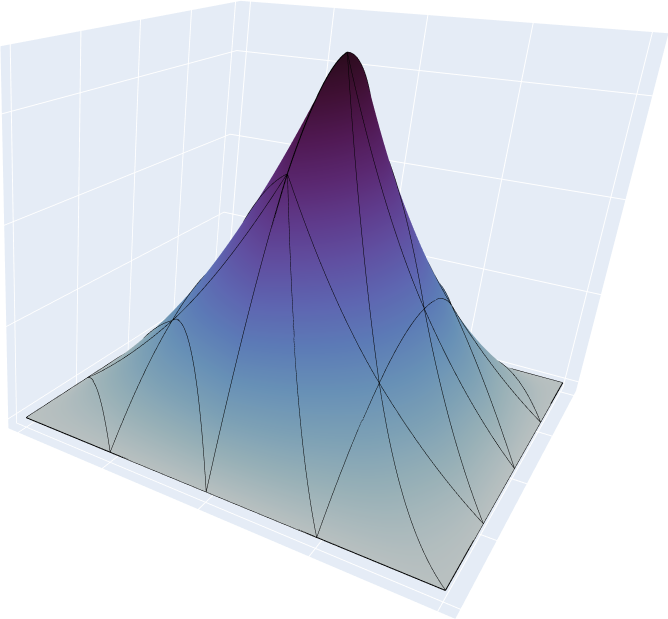}
    		\caption{}
    	\end{subfigure}
    	\begin{subfigure}{0.3\linewidth}
    		\centering
    		\includegraphics[width=0.9\linewidth]{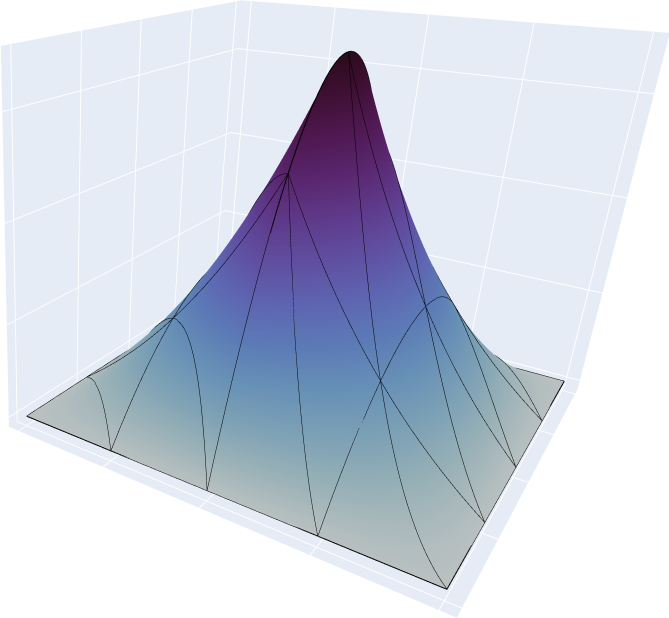}
    		\caption{}
    	\end{subfigure}
    	\begin{subfigure}{0.3\linewidth}
    		\centering
    		\includegraphics[width=0.9\linewidth]{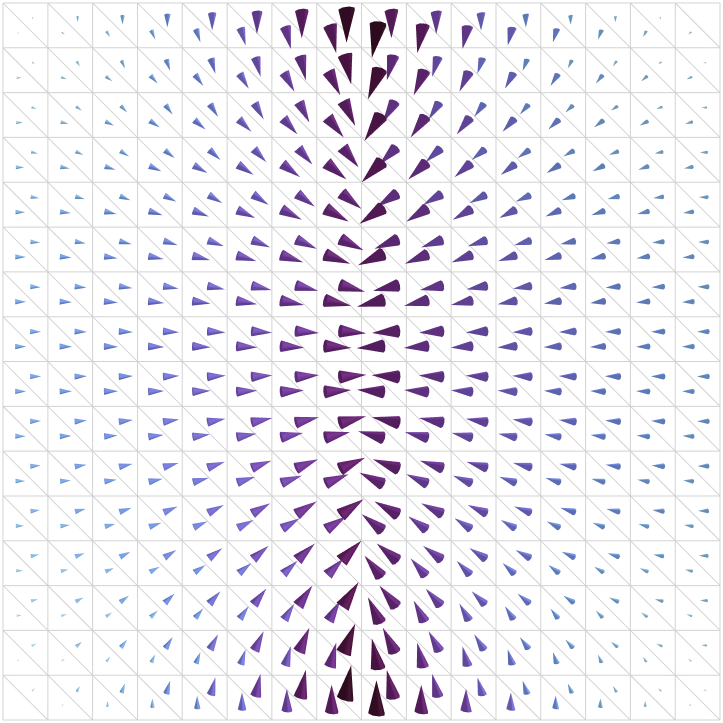}
    		\caption{}
    	\end{subfigure}
    	\begin{subfigure}{0.3\linewidth}
    		\centering
    		\includegraphics[width=0.9\linewidth]{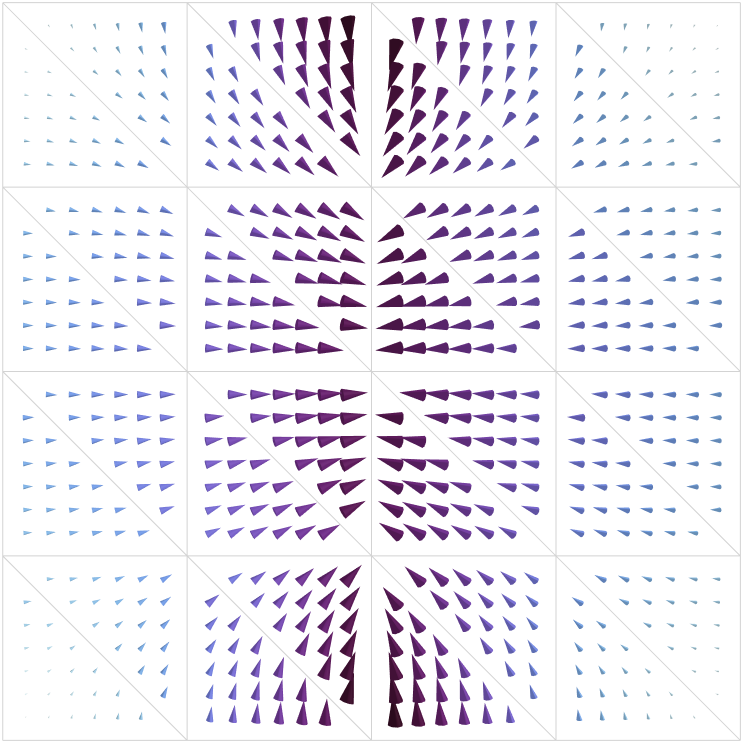}
    		\caption{}
    	\end{subfigure}
    	\begin{subfigure}{0.3\linewidth}
    		\centering
    		\includegraphics[width=0.9\linewidth]{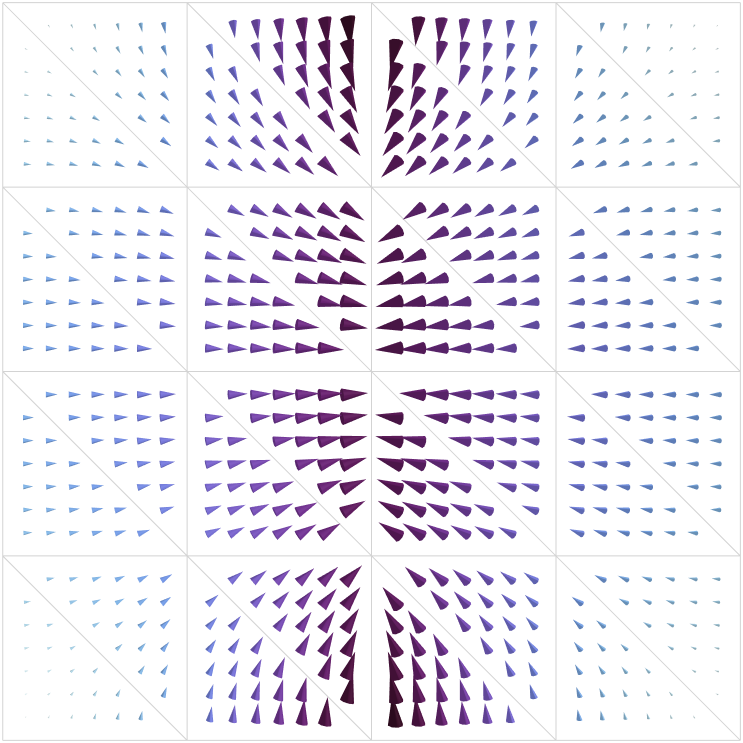}
    		\caption{}
    	\end{subfigure}
    	\caption{Displacement field using 1090 linear elements (a), 32 quadratic elements (b) and 32 cubic elements (c), corresponding to 1090, 258, and 530 degrees of freedom, respectively. The related microdistortion fields are depicted in (d)-(f).}
    	\label{fig:normal1}
    \end{figure}

\begin{figure}
    	\centering
    	\begin{subfigure}{0.48\linewidth}
    		\centering
    		\begin{tikzpicture}
	\definecolor{asl}{rgb}{0.4980392156862745,0.,1.}
	\definecolor{asb}{rgb}{0.,0.4,0.6}
	\begin{loglogaxis}[
		/pgf/number format/1000 sep={},
		axis lines = left,
		xlabel={Degrees of freedom},
		ylabel={$\| \widetilde{u} - u^h \|$},
		xmin=10, xmax=1e5,
		ymin=1e-6, ymax=1e3,
		xtick={1e1, 1e3, 1e5},
		ytick={1e-4, 1e-2, 1e0},
		legend pos= north east,
		ymajorgrids=true,
		grid style=dotted,
		]
		\addplot[color=asl, mark=pentagon] coordinates {
			( 26 ,  0.7841634360892709 )
			( 82 ,  0.24408098746290924 )
			( 290 ,  0.06499728682503178 )
			( 1090 ,  0.016523023999240258 )
			( 4226 ,  0.00414849192515152 )
			( 16642 ,  0.0010382458823288178 )
		};
		\addlegendentry{$\Lag^1 \times \Ned^0$}
		
		\addplot[color=violet, mark=triangle] coordinates {
			( 58 ,  0.10482972623626953 )
			( 194 ,  0.01344564736739643 )
			( 706 ,  0.0017116101417605969 )
			( 2690 ,  0.0002168273923953385 )
			( 10498 ,  2.734823359851182e-05 )
			( 41474 ,  3.4370707743099257e-06 )
		};
		\addlegendentry{$\Lag^2 \times \Nedtwo^1$}
		
		\addplot[color=asb, mark=o] coordinates {
			( 146 ,  0.00966153399844631 )
			( 530 ,  0.0005911476191940268 )
			( 2018 ,  3.5802107871026504e-05 )
			( 7874 ,  2.190350139820566e-06 )
			( 31106 ,  1.3526692952756594e-07 )
		};
		\addlegendentry{$\Ber^3 \times \Ned^2$}
		
		\addplot[color=blue, mark=diamond] coordinates {
			( 122 ,  0.009661533998708435 )
			( 434 ,  0.0005911473308038146 )
			( 1634 ,  3.580210757838578e-05 )
			( 6338 ,  2.1903500497684865e-06 )
			( 24962 ,  1.352669230655464e-07 )
		};
		\addlegendentry{$\Ber^3 \times \Nedtwo^2$}
		
		\addplot[dashed,color=black, mark=none]
		coordinates {
			(20, 4)
			(1e3, 0.08)
		};
	
	    \addplot[dashed,color=black, mark=none]
	    coordinates {
	    	(1e3, 3e-3)
	    	(5e4, 8.485281374238571e-06)
	    };
        
        \addplot[dashed,color=black, mark=none]
        coordinates {
        	(8e2, 7e-5)
        	(5e3, 1.7919999999999997e-06)
        };
		
	\end{loglogaxis}
	\draw (2.2,3.6) node[anchor=south]{$\mathcal{O}(h^{2})$};
	\draw (3.65,0.05) node[anchor=south]{$\mathcal{O}(h^{4})$};
	\draw (5.2,1.3) node[anchor=south]{$\mathcal{O}(h^{3})$};
\end{tikzpicture}
    		\caption{}
    	\end{subfigure}
    	\begin{subfigure}{0.48\linewidth}
    		\centering
    		\begin{tikzpicture}
	\definecolor{asl}{rgb}{0.4980392156862745,0.,1.}
	\definecolor{asb}{rgb}{0.,0.4,0.6}
	\begin{loglogaxis}[
		/pgf/number format/1000 sep={},
		axis lines = left,
		xlabel={Degrees of freedom},
		ylabel={$\| \widetilde{\vb{p}} - \vb{p}^h \|$},
		xmin=10, xmax=1e5,
        ymin=1e-5, ymax=1e5,
        xtick={1e1, 1e3, 1e5},
        ytick={1e-4, 1e-2, 1e0},
		legend pos= north east,
		ymajorgrids=true,
		grid style=dotted,
		]
		\addplot[color=asl, mark=pentagon] coordinates {
			( 26 ,  3.114764908738634 )
			( 82 ,  1.788526579488874 )
			( 290 ,  0.9362489808704159 )
			( 1090 ,  0.47413612476977995 )
			( 4226 ,  0.2378534772422181 )
			( 16642 ,  0.11902627255743328 )
		};
		\addlegendentry{$\Lag^1 \times \Ned^0$}
		
		\addplot[color=violet, mark=triangle] coordinates {
			( 58 ,  0.8223291420936867 )
			( 194 ,  0.22135058660792784 )
			( 706 ,  0.05750981739089367 )
			( 2690 ,  0.014666672042813498 )
			( 10498 ,  0.0037042456642980707 )
			( 41474 ,  0.0009308607827571709 )
		};
		\addlegendentry{$\Lag^2 \times \Nedtwo^1$}
		
		\addplot[color=asb, mark=o] coordinates {
			( 146 ,  0.09985668624512833 )
			( 530 ,  0.012787363436420847 )
			( 2018 ,  0.0015956028726516554 )
			( 7874 ,  0.00019851473839186068 )
			( 31106 ,  2.4732016431791918e-05 )
		};
		\addlegendentry{$\Ber^3 \times \Ned^2$}
		
		\addplot[color=blue, mark=diamond] coordinates {
			( 122 ,  0.09989931669179952 )
			( 434 ,  0.012788632010245177 )
			( 1634 ,  0.0015956329154497993 )
			( 6338 ,  0.0001985154557129703 )
			( 24962 ,  2.473203480226435e-05 )
		};
		\addlegendentry{$\Ber^3 \times \Nedtwo^2$}
		
		\addplot[dashed,color=black, mark=none]
		coordinates {
			(5e1, 1e1)
			(1e3, 2.2360679774997894)
		};
	    
	    \addplot[dashed,color=black, mark=none]
	    coordinates {
	    	(1e3, 8e-2)
	    	(5e4, 0.0016)
	    };
    
        \addplot[dashed,color=black, mark=none]
        coordinates {
        	(6e2, 4e-3)
        	(1e4, 5.878775382679627e-05)
        };
		
	\end{loglogaxis}
	\draw (2.4,3.15) node[anchor=south]{$\mathcal{O}(h)$};
	\draw (5.1,1.7) node[anchor=south]{$\mathcal{O}(h^{2})$};
	\draw (3.9,0.3) node[anchor=south]{$\mathcal{O}(h^{3})$};
\end{tikzpicture}
    		\caption{}
    	\end{subfigure}
    	\caption{Displacement (a) and microdistortion (b) fields of the discontinuous-normal problem and convergence rates under h-refinement (b).}
    	\label{fig:normal2}
    \end{figure}
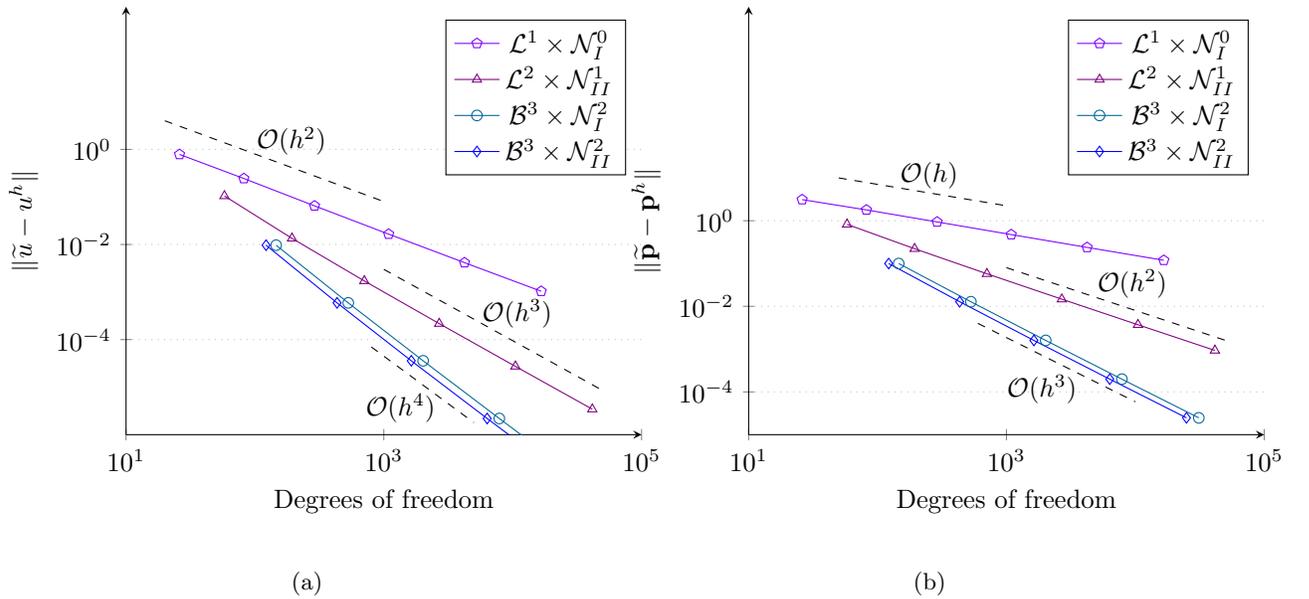

We demonstrate the behaviour of the elements constructed using linear and quadratic Lagrangian $\Lag$, and cubic Bernstein $\Ber$ polynomials. 
The solution of the displacement and microdistortion fields is depicted in \cref{fig:normal1}.
The difference between the quadratic and cubic displacement fields is subtle but can be observed by the sharp form of the solution at $(x,y) = (0,0)$.
As shown in \cref{fig:normal2}, the elements yield optimal convergence rates and correctly capture the discontinuity of the normal component of the microdistortion. This is also visible in the depiction of the fields in \cref{fig:normal1}.
In comparison, a formulation with $\vb{p} \in [\Hone(\surf)]^2$ would impose the higher regularity $\C^0$ on the microdistortion. This can be done using $\Hone$-conforming finite elements for the formulation of the microdistortion $\vb{p}$ and leads to sub-optimal convergence rates, see \cite{Sky2021}. 

\subsection{Three-dimensional relaxed micromorphic continuum}
The balance equations of the relaxed micromorphic continuum \cite{Neff2014} are given by
\begin{subequations}
		\begin{align}
			-\Di[\Ce \sym (\D \vb{u} - \bm{P}) + \Cc \skw (\D \vb{u} - \bm{P})] &= \vb{f} && \text{in} \quad \body \, , \label{eq:strong_u} \\
			-\Ce  \sym (\D \vb{u} - \Pm) - \Cc  \skw(\D \vb{u} - \Pm) + \Cmic \sym \Pm + \muma \, \Lc ^ 2  \Curl\Curl\Pm &= \bm{M} && \text{in} \quad \body \, , \label{eq:strong_p}  \\
			\vb{u} &= \widetilde{\vb{u}} && \text{on} \quad \surf_D^u \, , \\
			\Pm \times \, \vb{n} &= \widetilde{\Pm} \times \vb{n} && \text{on} \quad \surf_D^P \, , \label{eq:pdir} \\
			[\Ce \sym (\D \vb{u}- \bm{P}) + \Cc \skw (\D \vb{u} - \bm{P})] \, \vb{n} &= 0 && \text{on} \quad \surf_N^u \, ,\\
			\Curl \Pm \times \, \vb{n}  &= 0  && \text{on} \quad \surf_N^P \, .
		\end{align}
	    \label[Problem]{eq:full_relaxed}
	\end{subequations}
where $\ud$ is a three-dimensional displacement field, $\Pm$ is second order tensor field for the microdistortion, $\Ce, \Cmic$ are standard fourth order elasticity tensors, $\Cc$ is a rotational coupling tensor for infinitesimal rotations, $\muma$ is the macroscopic shear modulus and $\Lc$ is the characteristic length scale parameter. The body forces $\vb{f}$ and micro-moments $\bm{M}$ are now given by a three-dimensional vector a second order tensor, respectively. Further, the gradient $\D(\cdot)$ and $\Curl(\cdot)$ operators are constructed by applying the standard operators row-wise.
The bilinear form reads
\begin{align}
		a(\{\delta \ud , \delta \Pm\},\{\ud, \Pm\}) = \int_\body & \langle \sym(\D \delta \ud - \delta \Pm) , \, \Ce \sym(\D \ud - \Pm) \rangle + \langle \sym \delta \Pm, \, \Cmic \sym \Pm \rangle \notag \\
		& + \langle \skw(\D \delta \ud - \delta \Pm) , \, \Cc \skw(\D \ud - \Pm) \rangle + \muma \Lc^2 \langle \Curl \delta \Pm , \, \Curl \Pm \rangle \, \dd \body \, ,
		\label{eq:bi_full}
	\end{align}
and the corresponding linear form is given by
\begin{align}
		l(\{\delta \ud , \delta \Pm\}) = \int_\body \langle \delta \ud , \, \vb{f} \rangle + \langle \delta \Pm , \, \bm{M} \rangle \, \dd \body \, .
		\label{eq:li_full}
	\end{align}
The problem is well-posed in $\X(\body) = [\Hone(\body)]^3 \times [\Hc{,\body}]^3$, compare with \cite{SKY2022115298,GNMPR15,Neff_existence,Neff2015}.

We define the material constants $\lame = \mue = \lammi = \mumi= \muma = \Lc = 1$ and $\muc = 0$ and prescribe the solution
\begin{align}
		&\widetilde{\ud} = \left[\begin{matrix}0\\0\\\sin ( \pi x )\end{matrix}\right] \, , && \widetilde{\Pm} = \D \widetilde{\ud} + 10 (1-y^2) (1-z^2) \sin(\pi x) \begin{bmatrix}
			0 & 0 & 0 \\
			0 & 0 & 0 \\
			0 & -z & y 
		\end{bmatrix} \, ,
	\end{align}
for which the right-hand-side reads
\begin{align}
    	\vb{f} &=  \begin{bmatrix}10 \pi y (y^{2} - 1) (z^{2} - 1) \cos(\pi x )\\20 (- y^{2} + z^{2}) \sin(\pi x )\\y z (60 y^{2} - 20 z^{2} - 40) \sin(\pi x )\end{bmatrix} \, , \notag \\ 
    	\bm{M} &= \left [ \begin{matrix}20 y \left(y^{2} - 1\right) \left(z^{2} - 1\right) \sin{\left(\pi x \right)} & 0 & \\0 & 20 y \left(y^{2} - 1\right) \left(z^{2} - 1\right) \sin{\left(\pi x \right)} & \cdots \\ \pi \left(20 y^{3} z - 20 y z^{3} + 1\right) \cos{\left(\pi x \right)} & z \left(120 y^{2} - 10 \pi^{2} \left(y^{2} - 1\right) \left(z^{2} - 1\right) - 20 \left(y^{2} - 1\right) \left(z^{2} - 1\right) - 80\right) \sin{\left(\pi x \right)} & \end{matrix} \right . \notag \\
    	& \qquad \left . \begin{matrix} 
    		&\pi \cos{\left(\pi x \right)} \\
    		\cdots &- 20 z \left(y^{2} - 1\right) \left(z^{2} - 1\right) \sin{\left(\pi x \right)} \\
    		&y \left(- 120 z^{2} + 60 \left(y^{2} - 1\right) \left(z^{2} - 1\right) + 10 \pi^{2} \left(y^{2} - 1\right) \left(z^{2} - 1\right) + 80\right) \sin{\left(\pi x \right)}
    	\end{matrix} \right ] \, .
    \end{align}
We set the boundary $\surf_D = \partial \body$, such that the prescribed solution corresponds with the analytical solution due to uniqueness. 

\begin{figure}
    	\centering
    	\begin{subfigure}{0.48\linewidth}
    		\centering
    		\includegraphics[width=1\linewidth]{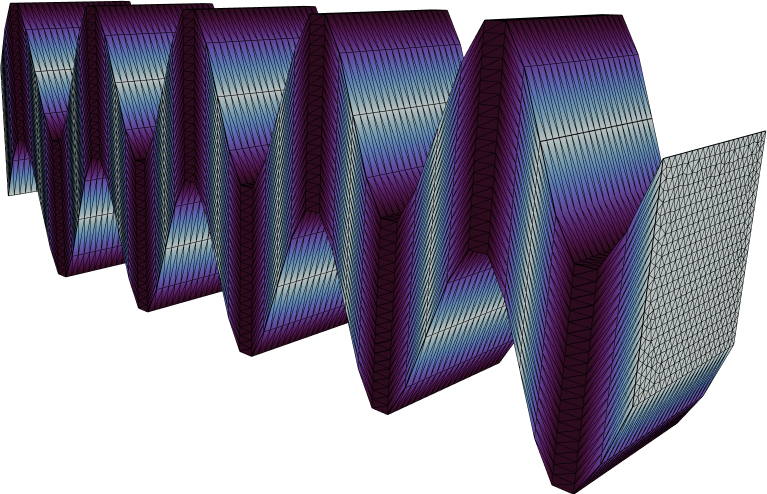}
    		\caption{}
    	\end{subfigure}
    	\begin{subfigure}{0.48\linewidth}
    		\centering
    		\includegraphics[width=1\linewidth]{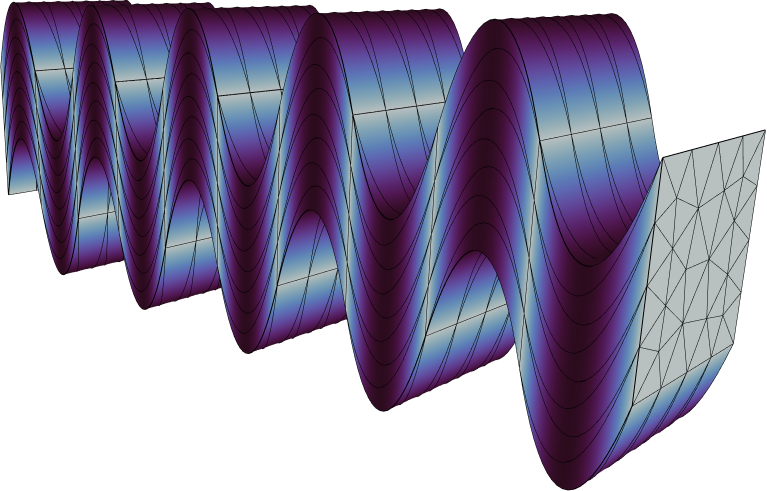}
    		\caption{}
    	\end{subfigure}
    	\begin{subfigure}{0.48\linewidth}
    		\centering
    		\includegraphics[width=1\linewidth]{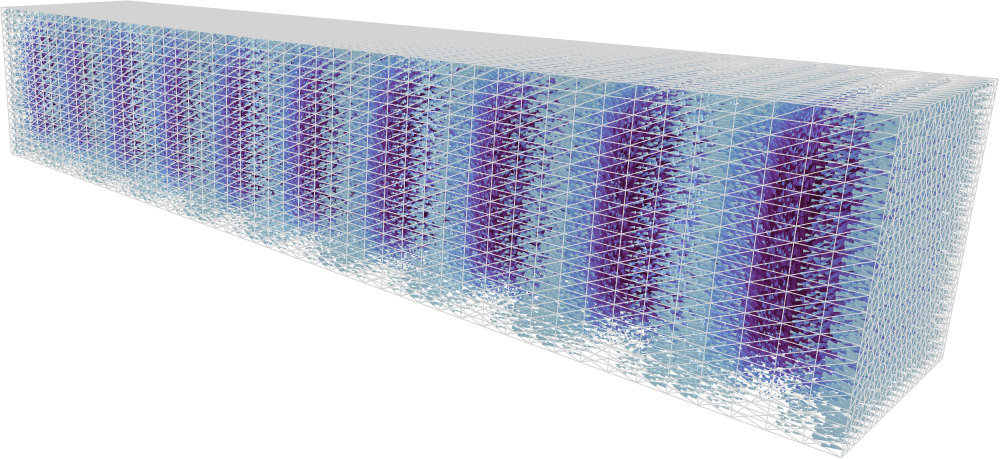}
    		\caption{}
    	\end{subfigure}
    	\begin{subfigure}{0.48\linewidth}
    		\centering
    		\includegraphics[width=1\linewidth]{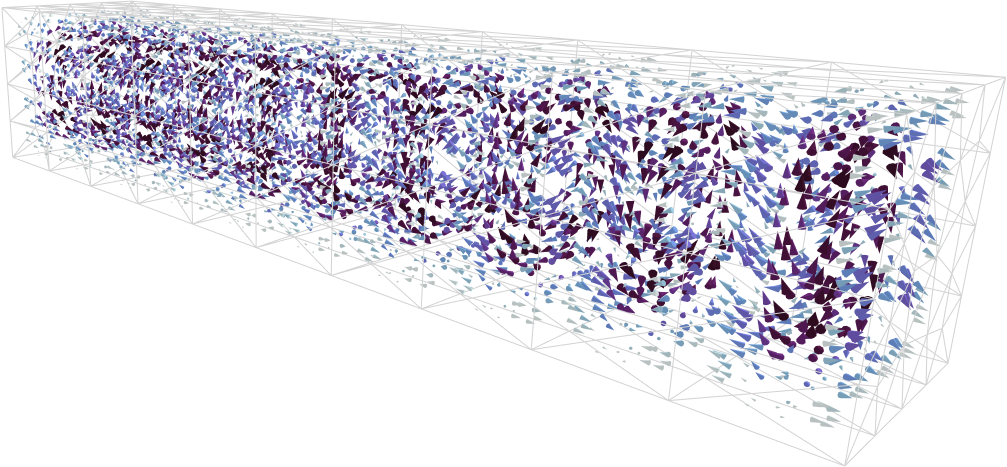}
    		\caption{}
    	\end{subfigure}
    	\caption{Depiction of the displacement field  (a),(c) and the last row of the microdistortion field (b),(d) for  $141900$ linear elements corresponding to $597820$ degrees of freedom and $1260$ cubic elements corresponding to $96358$ degrees of freedom.}
    	\label{fig:conv3d2}
    \end{figure}

\begin{figure}
    	\centering
    	\begin{subfigure}{0.48\linewidth}
    		\centering
    		\begin{tikzpicture}
	\definecolor{asl}{rgb}{0.4980392156862745,0.,1.}
	\definecolor{asb}{rgb}{0.,0.4,0.6}
	\begin{loglogaxis}[
		/pgf/number format/1000 sep={},
		axis lines = left,
		xlabel={Degrees of freedom},
		ylabel={$\| \widetilde{\ud} - \ud^h \|_{\Le}$},
		xmin=5000, xmax=1.5e6,
		ymin=1e-3, ymax=1e5,
		xtick={1e4,1e5,1e6},
		ytick={1e-1, 1e1},
		legend pos= north east,
		ymajorgrids=true,
		grid style=dotted,
		]
		\addplot[color=asl, mark=triangle] coordinates {
			(6364, 4.106011728956118)
			(11476,2.2474442972200017) 
			(26632, 1.3611198038372418) 
			(81550, 0.5900232773323868) 
			(195574, 0.3440619257596154) 
			(597820, 0.15193166252301096)
		};
		\addlegendentry{$\Lag^1 \times \Ned^0$}
		
		\addplot[color=violet, mark=o] coordinates {
			(17110, 0.1811877778079692)
			(30994, 0.401411432437721)
			(72244, 0.1749067651896356)
			(222184, 0.04680374735766063)
			( 534088 , 0.0203953698418525 )
		};
		\addlegendentry{$\Lag^2 \times \Nedtwo^1$}
		
		\addplot[color=cyan, mark=diamond] coordinates {
			( 10930 , 0.790109909735829 )
			(60520, 0.15372327516231174)
			(111208, 0.05572218846142478)
			(262912, 0.018061987659246076)
			( 437764 , 0.007494786894627259 )
		};
		\addlegendentry{$\Ber^{3} \times \Nedtwo^{2}$}
		
		\addplot[dashed,color=black, mark=none]
		coordinates {
			(1e4, 1e1)
			(1e5, 2.1544346900318843)
		};
	
	    \addplot[dashed,color=black, mark=none]
	    coordinates {
	    	(8e4, 5e-2)
	    	(3e5, 0.008582127863161156)
	    };
		
	\end{loglogaxis}
	\draw (2.3,2.6) node[anchor=south]{$\mathcal{O}(h^{2})$};
	\draw (4.05,0.27) node[anchor=south]{$\mathcal{O}(h^{4})$};
\end{tikzpicture}
    		\caption{}
    	\end{subfigure}
    	\begin{subfigure}{0.48\linewidth}
    		\centering
    		\begin{tikzpicture}
	\definecolor{asl}{rgb}{0.4980392156862745,0.,1.}
	\definecolor{asb}{rgb}{0.,0.4,0.6}
	\begin{loglogaxis}[
		/pgf/number format/1000 sep={},
		axis lines = left,
		xlabel={Degrees of freedom},
		ylabel={$\| \widetilde{\Pm} - \Pm^h \|_{\Le}$},
		xmin=5000, xmax=1.5e6,
		ymin=1e-2, ymax=1e5,
		xtick={1e4,1e5,1e6},
		ytick={1e-1, 1e1},
		legend pos= north east,
		ymajorgrids=true,
		grid style=dotted,
		]
		\addplot[color=asl, mark=triangle] coordinates {
			(6364, 16.98453410833475) 
			(11476, 12.450876564592395) 
			(26632, 9.8526437445637) 
			(81550, 6.681327992847365) 
			(195574, 5.046808088038453) 
			(597820, 3.367682097548589)
		};
		\addlegendentry{$\Lag^1 \times \Ned^0$}
		
		\addplot[color=violet, mark=o] coordinates {
			(17110, 9.75194508934138)
			(30994, 7.311486249234136)
			(72244, 4.563095731419553)
			(222184, 2.1100134568385815)
			( 534088 , 1.2107611533054206 )
		};
		\addlegendentry{$\Lag^2 \times \Nedtwo^1$}

		\addplot[color=cyan, mark=diamond] coordinates {
			( 10930 , 10.873495610584424 )
			(60520, 2.6847626548950787)
			(111208, 1.3761377373621508)
			(262912, 0.5921110207956943)
			( 437764 , 0.33107433517242935 )
		};
		\addlegendentry{$\Ber^{3} \times \Nedtwo^{2}$}
		
		\addplot[dashed,color=black, mark=none]
		coordinates {
			(1e4, 5e1)
			(1e5, 23.207944168063896)
		};
	
	    \addplot[dashed,color=black, mark=none]
	    coordinates {
	    	(5e4, 7e-1)
	    	(5e5, 0.07)
	    };
		
	\end{loglogaxis}
	\draw (1.8,2.83) node[anchor=south west]{$\mathcal{O}(h)$};
	\draw (3.6,0.4) node[anchor=south west]{$\mathcal{O}(h^3)$};
\end{tikzpicture}
    		\caption{}
    	\end{subfigure}
    	\caption{Convergence under h-adaption for various polynomial orders.}
    	\label{fig:conv3d}
    \end{figure}
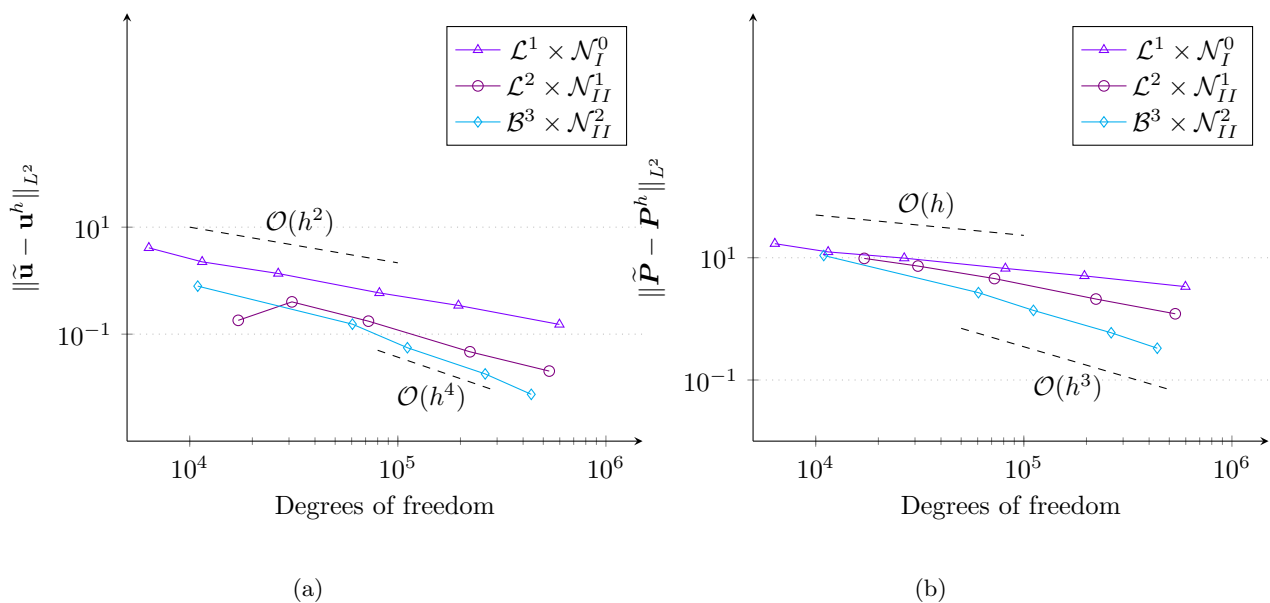

A depiction of the approximation is given in \cref{fig:conv3d2}. In order to capture the wave-shaped solution, either a fine grid or higher order elements are needed. The convergence rates shown in \cref{fig:conv3d} are optimal. Although the initial solution of the quadratic element on the coarse grid seems well enough, it is simply due to the fortunate placing of element interfaces at the peaks and cannot be expected in general discretizations. All discretizations converge at an optimal rate, compare \cite{SKY2022115298}.


\section{Conclusions and outlook}
In this work we presented a novel method of constructing base functions for the spaces $\Hc{}$ and $\Hd{}$ on the reference triangle and tetrahedron. The ability of the method to generate unisolvent bases is proven by linear independence and conformity theorems. Further, the validity of the construction is demonstrated by two examples of the relaxed micromorphic model using Lagrange and Bernstein base functions for the formulation. As such, we conclude that the method can be used to generate $\Hc{}$- and $\Hd{}$-conforming subspaces for a variety of $\Hone$-conforming polynomial subspaces. The simplicity of the method makes its application straight-forward in the formulation of arbitrary order vectorial finite elements, as one must simply employ a higher order $\Hone$-conforming subspace in the construction. Further, by construction, the vectorial basis can inherit characteristics of the underlying scalar basis. For example, the Kronecker delta property of Lagrange polynomials as demonstrated in the embedding of boundary conditions in \cite{SKY2022115298}, or the optimal complexity of Bernstein polynomials \cite{AinsworthOpt}.

This work did not discuss alternative reference elements, such as quads, prisms, pyramids or hexahedra. Further, the definition of Raviart-Thomas and N\'ed\'elec elements of the first type on the reference tetrahedron has not been addressed. The latter are topics for future works.

\section*{Acknowledgements}

\bibliographystyle{spmpsci}   

\bibliography{Ref}   

\appendix

\section{Piola transformations} \label{ap:a}
Consistent transformations are employed to map the base functions from the reference element to the physical element \cite{Mon03}. Effectively, every element in the physical domain is mapped by the same reference domain. If the mapping of the physical space is achieved via barycentric functions, the polynomial degree is maintained across transformations. 
	
	Scalar base functions transform according to
	\begin{align}
		&n(\vb{x}) = n \circ [\vb{x}^{-1}(\bm{\xi})] \, ,  &&
		\nabla_x n = \bm{J}^{-T} \nabla_\xi n \, ,
	\end{align}
	where the result concerning the Jacobi matrix is a direct consequence of the chain rule.
	
	N\'ed\'elec elements are defined via their action on the tangent vectors of the element. Consequently, a consistent transformation is given by the equality 
	\begin{align}
		\langle \bm{\theta}, \, \vb{t} \rangle \dd \curv = \langle \bm{\theta}, \, \dd \vb{s} \rangle = \langle \bm{\theta},  \, \bm{J} \dd \bm{\mu} \rangle = 
		\langle \bm{\vartheta}, \, \dd \bm{\mu} \rangle = \langle \bm{\vartheta}, \, \bm{\tau} \rangle \dd \mu \quad \iff \quad \bm{\theta} = \bm{J}^{-T} \bm{\vartheta} \, ,
	\end{align}
	known as the covariant Piola transformation.
	This is the same transformation as for gradients, thus respecting the commuting property \cref{fig:derham}.
	Further, vectors undergoing the latter transformation exhibit the following transformation of the curl operator 
	\begin{align}
		\mathrm{curl}_x \bm{\theta} = \nabla_x \times \bm{\theta} = (\bm{J}^{-T} \nabla_\xi) \times (\bm{J}^{-T} \bm{\vartheta}) = \cof(\bm{J}^{-T}) (\nabla_\xi \times \bm{\vartheta}) = \dfrac{1}{\det\bm{J}} \bm{J} \mathrm{curl}_\xi \bm{\vartheta} \, ,
	\end{align}
	being the so called contravariant Piola transformation.
	The result is won by observing that
	\begin{align}
		\nabla_x \times \bm{J}^{-T} = \nabla_x \times \nabla_x \bm{\xi} = 0 \, .
	\end{align}
    For two-dimensional domains the formula reduces to
    \begin{align}
    	\mathrm{div}_x (\bm{R} \bm{\theta}) = \dfrac{1}{\det \bm{J}} \mathrm{div}_\xi (\bm{R} \bm{\vartheta}) \, ,
    \end{align}
    since the curl operator produces a scalar.
	The contravariant Piola transformation is compatible with the commuting diagram and preserves normal projections on the element's boundary. To see this characteristic define the base function $\bm{\phi}$ in the reference domain and $\bm{\varphi}$ in the physical domain and equate their normal projections on the outer surface of both domains
	\begin{align}
		\langle \bm{\varphi}, \, \vb{n} \rangle \dd \surf = \langle \bm{\varphi} , \, \dd \vb{A} \rangle = \langle \bm{\varphi} , \, \cof(\bm{J}) \dd \bm{\Gamma} \rangle = \langle \bm{\phi}, \, \dd \bm{\Gamma} \rangle = \langle \bm{\phi} , \, \bm{\nu} \rangle \dd \Gamma \quad \iff \quad \bm{\varphi} = \dfrac{1}{\det \bm{J}} \bm{J} \bm{\phi} \, .  
	\end{align}
The divergence of functions mapped by a contravariant Piola transformation is given by
\begin{align}
	\int_\body q \, \mathrm{div}_x\bm{\varphi} \, \dd  \body &= 
	\oint_{\partial \body} q \, \langle \bm{\varphi} , \, \vb{n} \rangle \, \dd \surf - \int_\body \langle \nabla_x q , \, \bm{\varphi} \rangle \, \dd \body \notag \\
	&= \oint_{\partial \Omega} \hat{q} \, \langle \dfrac{1}{\det \bm{J}} \bm{J} \, \bm{\phi} , \, \det(\bm{J}) \, \bm{J}^{-T} \bm{\nu} \rangle \, \dd \Gamma - \int_\Omega \langle \bm{J}^{-T} \nabla_\xi \hat{q} , \, \dfrac{1}{\det \bm{J}} \bm{J} \, \bm{\phi} \rangle \, \det \bm{J} \, \dd \Omega \notag \\
	&= \oint_{\partial \Omega} \hat{q} \, \langle \bm{\phi} , \, \bm{\nu} \rangle \, \dd \Gamma - \int_\Omega \langle \nabla_\xi \hat{q} , \, \bm{\phi} \rangle \, \dd \Omega \notag \\
	&= \int_\Omega \hat{q} \, \mathrm{div}_\xi \bm{\phi} \, \dd \Omega = \int_\body q \, \mathrm{div}_\xi(\bm{\phi}) \, \dfrac{1}{\det \bm{J}} \dd \body \qquad \forall\, q \, \in \mathit{C}^\infty(\overline{\body}) \, ,
\end{align}
where $\hat{q} = q \circ \vb{x}$. Consequently, there holds
\begin{align}
	\mathrm{div}_x\bm{\varphi}= 
	\dfrac{1}{\det \bm{J}}\, \mathrm{div}_\xi \bm{\phi} \, .
\end{align}

\end{document}